\documentclass[a4paper]{article}

\pdfoutput=1
\usepackage{latexsym}
\usepackage{graphicx}
\usepackage[a4paper]{geometry}
\usepackage{amsmath}
\usepackage{amssymb}
\usepackage{hyperref}
\usepackage{amsfonts}
\usepackage{commath}
\usepackage{mdframed}
\usepackage[usenames,dvipsnames]{xcolor}
\usepackage{float}
\usepackage{enumitem}
\usepackage[natbib,style=numeric,defernumbers]{biblatex}
\newcommand\Tstrut{\rule{0pt}{2.6ex}}
\newcommand\Bstrut{\rule[-0.9ex]{0pt}{0pt}}

\newcommand{\notes}[1]{{
    \sf
}}

\newcommand{\solution}{\vspace{1.5ex}\textbf{Solution and comments:} }

\title{Teaching Resources for Embedding Ethics in Mathematics: \\ Exercises, Projects, and Handouts}

\author{Maurice Chiodo\footnote{Centre for the Study of Existential Risk, University of Cambridge, United Kingdom. \texttt{\href{mailto:mcc56@cam.ac.uk}{mcc56@cam.ac.uk}}.} \\
Dennis M\"uller\footnote{RWTH Aachen University, Germany. \texttt{\href{mailto:dennis.mueller3@rwth-aachen.de}{dennis.mueller3@rwth-aachen.de}}} \\
Rehan Shah \footnote{School of Engineering and Materials Science, Queen Mary University of London, United Kingdom. \texttt{\href{mailto:rehan.shah@qmul.ac.uk}{rehan.shah@qmul.ac.uk}}.}}   
\date{\today}

\begin{document}

\maketitle
\begin{abstract}
    The resources compiled in this document provide an approach to embed and teach Ethics in Mathematics at the undergraduate level. We provide mathematical exercises and homework problems that teach students ethical awareness and transferable skills, for many of the standard courses in the first and second years of a university degree in mathematics or related courses with significant mathematical content (e.g., physics, engineering, computer science, economics, etc). In addition to the exercises, this document also contains a list of projects, essay topics, and handouts for use as final projects and in seminars. This is a living document, and additional contributions are welcome.

\end{abstract}
\vspace{+60pt}
\underline{\textit{Acknowledgements:}} Some of these questions were designed during a project with undergraduate students at the University of Cambridge during the summer of 2018. We are very grateful for their work and continued support in the initial years of the  Ethics in Mathematics Project\footnote{The Ethics in Mathematics Project: \href{https://www.ethics-in-mathematics.com}{https://www.ethics-in-mathematics.com}}. We also wish to thank the many colleagues who have provided invaluable commentary while we were compiling this document, in particular Tony Gardiner and Thomas King for their valuable comments and feedback. 
\\Some of the mathematical content in these questions was inspired by questions on the first and second year example sheets provided by the Faculty of Mathematics at the University of Cambridge.  
\\ \\ 
\underline{\textit{Work in Progress:}} This document is a work in progress; a \textit{living document}. If you want to contribute to it by adding additional content, such as questions, projects, or reading list items, please get in touch. If you find any errors, please also let us know.

\newpage 
\setcounter{tocdepth}{1}
\tableofcontents
\newpage
\section{Introduction}

\subsection{What do we mean by ``Embedded Ethics in mathematics''?}

Over the last few decades, we have seen a sociopolitical turn in mathematics education that introduced the study of equity, social justice, identity and power into the apprenticeship of future mathematicians and mathematics teachers, thereby highlighting the impact of existing political and moral discourses and the contingent institutional contexts on both mathematics and its education\footnote{R. Gutiérrez (2013). The Sociopolitical Turn in Mathematics Education. \textit{Journal for Research in Mathematics Education} 44:1, pp. 37–-68. \href{https://doi.org/10.5951/jresematheduc.44.1.0037                 }{https://doi.org/10.5951/jresematheduc.44.1.0037  }}. More recently, this has been extended to include approaches that try to teach subject-specific ethical awareness to mathematicians.

What is Ethics in Mathematics?
Ethics in Mathematics tries to abstract from the subject-specific ethical discourse that is happening in many mathematical sciences (e.g., data science, AI, finance, and cryptography) to build a holistic approach towards researching, understanding and teaching the potential ethical questions that arise when doing \textit{any} form of mathematics\footnote{D. M\"uller (2022). Situating ``Ethics in Mathematics'' as a Philosophy of Mathematics Ethics Education, p. 5. \href{https://arxiv.org/abs/2202.00705}{https://arxiv.org/abs/2202.00705} \\ D. M\"uller, M. Chiodo, \& J. Franklin (2022). A Hippocratic Oath for mathematicians? Mapping
the landscape of ethics in mathematics. \textit{Science and Engineering Ethics} 28:5, article 41. \href{https://doi.org/10.1007/s11948-022-00389-y                                                                  }{https://doi.org/10.1007/s11948-022-00389-y              }  }. This approach is in line with a broader trend in the philosophy of mathematics. For many years, we have seen an increasing number of philosophers focus their work on mathematical practice\footnote{R. Wagner (2017). \textit{Making and Breaking Mathematical Sense: Histories and Philosophies of Mathematical Practice}. Princeton: Princeton University Press. \\
R. Wagner (2023). The Ethical Charge of Articulating Mathematics. \textit{Global Philosophy} 33:35.\href{https://doi.org/10.1007/s10516-023-09686-y                                                                                              }{https://doi.org/10.1007/s10516-023-09686-y                        }},  the potential for injustices\footnote{C.J. Rittberg (2023). Justified Epistemic Exclusions in Mathematics. \textit{Philosophia Mathematica}. \href{https://doi.org/10.1093/philmat/nkad008                                                                                                         }{https://doi.org/10.1093/philmat/nkad008   }
\\ C.J. Rittberg, F.S. Tanswell, \& J.P. van Bendegem (2020). Epistemic Injustice in Mathematics. \textit{Synthese} 197:9, pp. 3875--3904. \href{ https://doi.org/10.1007/s11229-018-01981-1                                                                                                                               }{ https://doi.org/10.1007/s11229-018-01981-1                     } 
\\E. Hunsicker, \& Colin Jakob Rittberg (2022). On the Epistemological Relevance of Social Power and Justice in Mathematics. \textit{Axiomathes} 32, pp. 1147–-68. \href{https://doi.org/10.1007/s10516-022-09629-z                                                                                                          }{https://doi.org/10.1007/s10516-022-09629-z                          }} and harm\footnote{P. Ernest (2018). The Ethics of Mathematics: Is Mathematics Harmful? \textit{The Philosophy of Mathematics Education Today}, pp. 187–-216: Springer, Cham. \href{https://link.springer.com/chapter/10.1007/978-3-319-77760-3_12                                                                                                                                                                                                           }{https://link.springer.com/chapter/10.1007/978-3-319-77760-3                                                        \_12                                                              }}, as well as potential institutional solutions such as codes of conduct\footnote{C. Buell, V. I. Piercey, and R. E. Tractenberg (2022). Leveraging guidelines for ethical practice of statistics and computing to develop standards for ethical mathematical practice: A White Paper. \href{https://arxiv.org/abs/2209.09311}{https://arxiv.org/abs/2209.09311}} or even Hippocratic oaths\footnote{D. M\"uller, M. Chiodo, \& J. Franklin (2022). A Hippocratic Oath for mathematicians? Mapping
the landscape of ethics in mathematics. \textit{Science and Engineering Ethics} 28:5, article 41. \href{https://doi.org/10.1007/s11948-022-00389-y                                                                                                                           }{https://doi.org/10.1007/s11948-022-00389-y                              }
\\ C. J. Rittberg (2023). Hippocratic Oaths for Mathematicians? \textit{Philosophia} 51:3, pp. 1579–-1603. \href{https://doi.org/10.1007/s11406-022-00588-8                                                                           }{https://doi.org/10.1007/s11406-022-00588-8                                                        } }. But unlike the philosophical discourse, the mathematicians' discourse on Ethics in Mathematics almost always also contains a practical question: what is the ethics of \textit{my} mathematics? 

For future mathematicians to be able to answer this question, they will first need to learn how to see ``ethics'' in the context of their own mathematical practice. Mathematics is hard and ethics is hard, and one of our biggest lessons over the last few years has been that to answer these questions, mathematicians will need to be adequately \textit{prepared} to answer them. They don't just need to be primed to spot them, but also need the skills to navigate the plurality of arguments and reasons that are part of every ethical discourse. Thus, teaching ethics in mathematics is a bit different to teaching a theorem, as it cannot be succinctly packaged: it is a \textit{skill} to be \textit{developed} and \textit{maintained}, just like other forms of ``mathematical reasoning'', rather than a fact or well-defined process to be learned. It requires a varied and diverse pedagogical approach, making use of different tools at different points in the learning process.

We, mathematicians, are only just beginning to grapple with this, but over the years, two complementary approaches have been developed to address this need: standalone ethics courses and, more recently, embedded ethics\footnote{Embedded ethics has been primarily pioneered as an approach to teach ethics in computer science. The approach aims to include the teaching of ethical awareness into the entire computer science curriculum, including homework problems, computational problems and standalone courses. One such approach can be found in the \textit{Embedded EthiCS} programme at Harvard: \href{https://embeddedethics.seas.harvard.edu}{https://embeddedethics.seas.harvard.edu}.}. This document is a first attempt to show how properly embedding ethics into the learning of mathematics can be done at the undergraduate level. Here, \textit{embedded ethics in mathematics} means embedding the teaching of subject-specific ethical awareness into all areas of mathematical education, and thus in particular into the exercises and projects through which students learn mathematics. Only through constant teaching and practice can the necessary skills be developed and maintained.

\subsection{What resources we provide: Exercises, Projects, and Handouts}
This document is divided into three types of resources for you to use and integrate into your teaching: exercises, projects, and handouts. Each category is important and useful, and serves a different purpose in the learning process. In this document, we have chosen to provide resources that are as easy as possible to ``plug and play'' into your existing teaching, without the need for any course or syllabus redesigns. The main component of this collection, in terms of length, detail, and depth, is the exercises, and we give a thorough explanation of how to make best use of those.

\subsubsection*{1. Exercises}
These have been written as closely as possible to standard undergraduate mathematics exercises, as would be typically found on an exercise sheet. The difficulty of the mathematics in each exercise is roughly equivalent to the difficulty of standard first or second year mathematics questions. That is to say, there is \textit{real} mathematics that has to be done in each question, of a level comparable to what one might expect in a first or second year course. This is not ``watered down mathematics'' in any sense.

Each of these exercises contains some ethical component to the problem. In order to solve the problem fully, the student will need to take into account this ethical aspect and consider it as part of their solution.

Each exercise is a standalone question, and you could could easily add one, or several, of these to your existing exercise sheets. You could also make use of some of these as worked examples in lectures, and discuss and explain both the mathematical, and societal, aspects of the question in that context. What we think would work best is if you made use of some of these as examples in lectures, to normalise the idea to students of ethical considerations in mathematical work, and then include some in exercise sheets for students to work through themselves.

We have provided ``solutions''  to as many of these exercises as possible. Our solutions include both a full exposition of the mathematical component of the question, as well as a discussion and incorporation of the societal and ethical issues that are embedded in the question. What we provide is a minimal coverage of the solution, and no doubt for many of these you may be able to add extra context, perspective, and insight, which we invite you to do. You should also feel free to modify and alter any part of any exercise to make it fit into your course more appropriately; we can only guess at what standard courses might contain, and what associated questions might be appropriate and relevant. Please note that, while every effort has been made to ensure the (mathematical) accuracy of these solutions, we cannot guarantee they are free from errors. So please check these thoroughly as you incorporate them into your teaching.

We have tried to cover as many courses as we could. What we have so far is a collection of exercises that cover \textit{most} standard first year undergraduate courses, and \textit{some} second year courses. We did not have the time, expertise, or resources, to go beyond that. But we encourage you to come up with more such questions if you can, and moreover to share them with us to add to this resource for others to make use of.   

\subsubsection*{2. Projects}
Some of the questions and problems that we came up with are more substantial than may be appropriate on an exercise sheet. So we split these off as self-contained projects for you to set students. Typically, these might take between one whole evening, and several days, to complete. They invite students to think about, and to explore a problem more deeply - both the mathematics, and the societal considerations. These will usually involve a degree of investigation, further reading, and research by the student, to complete fully.

These can be used as part of a course/module assessment. They can be used as optional additional enhancement work. They can be used for winter/summer break project work. They can be used as group work. Or in any other suitable context that may come up. We have not provided any solutions to these, as there are typically no ``right'' answers. The aim is for students to carry out their own investigations, find out what they can, and come to their own conclusions.

\subsubsection*{3. Handouts}

We have provided an Ethics in Mathematics reading list for students to go through. Feel free to use these in any combination you see fit, and also to modify them to better-suit your audience.

\subsubsection*{A note on copyright}
This material has been released under Creative Commons license CC BY-SA 4.0. Provided you comply with the requirements of the license\footnote{Full details available at \href{https://creativecommons.org/licenses/by-sa/4.0}{https://creativecommons.org/licenses/by-sa/4.0}}, you should be fine to include and adapt this material in your teaching resources, including all worked solutions and associated discussion text.

\subsection{Purpose and methodology behind the exercises}

The exercises are designed to teach mathematics, transferable skills and ethical  awareness in parallel. They cover many of the core topics of the first and second years of mathematical training at many universities. We tried to cover as many courses and areas as possible, both pure and applied. In addition, we tried to cover as many non-technical skills as possible. Each exercise has a tag that explains the topic, allied disciplines of study, approximate year levels of study, and which \textit{ethical pillar} (i.e., non-technical skill and/or ethical area) it is related to. The tags that we use for allied disciplines include \textit{Computer Science, Economics, Engineering, Natural Sciences, Probability, and Statistics}. The ethical pillars mentioned in the exercises correspond to those outlined in the \textit{Manifesto for the Responsible Development of Mathematical Works}\footnote{M. Chiodo \& D. M\"uller (2023). Manifesto for the Responsible Development of Mathematical Works -- A Tool for Practitioners and Management. \href{https://arxiv.org/abs/2306.09131}{https://arxiv.org/abs/2306.09131}}, which we repeat here for completeness:

\newpage
\label{list:pillars}
\subsubsection*{Ethical pillars:}
\begin{enumerate}
\item \textbf{Deciding whether to begin:} Why are you providing this mathematical product or service, and should you even do so?
\item \textbf{Diversity and perspectives:} Do your co-workers, superiors, and you have sufficient perspective, and do you understand the limitations and biases in your thinking?
\item \textbf{Handling data and information:} Are you using authorised and morally obtained datasets, in a responsible manner?
\item \textbf{Data manipulation and inference:} Do you have the expertise to properly manipulate data ensuring quality and ethics?
\item \textbf{The mathematisation of the problem:} What optimisation objectives and constraints have you chosen, and what are their real-life consequences? Who might the other impacted parties be?
\item \textbf{Communicating and documenting your work:} Are you properly considering how to comment and document your work and communicate the results to those who need them?
\item \textbf{Falsifiability and feedback loops:} Is your work falsifiable, and can you handle its large-scale impact and any feedback loops that arise?
\item \textbf{Explainable and safe mathematics:} Is your mathematical output explainable, and followed up with proper monitoring and maintenance?
\item \textbf{Mathematical artefacts have politics:} Are you aware of other non-mathematical aspects and the political nature of your work? What do you do to earn trust in yourself and your product?
\item \textbf{Emergency response strategies:} Do you have a non-technical response strategy for when things go wrong? Do you have a support network, including peers who support you and with whom you can talk freely?
\end{enumerate}

\begin{figure}[ht!]
    \centering
    \includegraphics[width=0.8\textwidth]{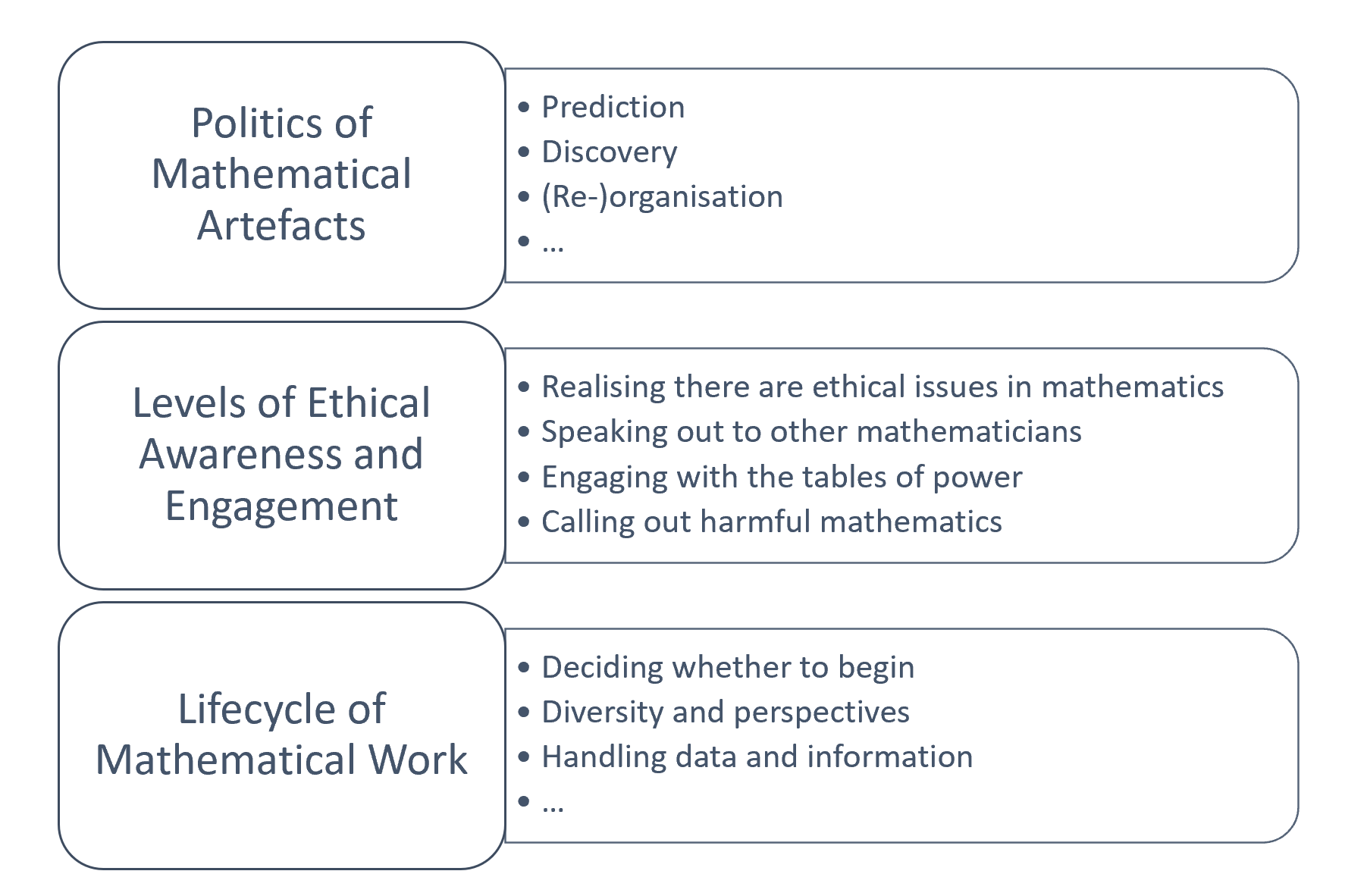}
    \caption{Elements of Embedded Ethics in Mathematics \protect\footnotemark}
    \label{fig:embedded_ethics}
\end{figure}

The exercises we provide  are designed with the insight we have gained from nine years of experience teaching ethics in mathematics at Cambridge. In particular, the educational philosophy of designing exercises by choosing one or more areas of mathematics, one or more pillars and an additional political dimension\footnote{D. M\"uller \& M. Chiodo (2023). Mathematical Artifacts Have Politics: The Journey from Examples to Embedded Ethics. \href{https://arxiv.org/abs/2308.04871}{https://arxiv.org/abs/2308.04871}} (see Figure \ref{fig:embedded_ethics}).

\footnotetext{ibid., p. 13.}

\subsection{How to use the exercises}
You will need to approach these exercises slightly differently from your usual mathematical exercises. Besides the mathematical knowledge, students need to understand that 
\begin{enumerate}
 \setlength\itemsep{-0.3em}
\item Such exercises exist,
\item There are good reasons to complete these exercises,
\item All exercises provide them with valuable skills,
\item These are interesting and challenging,
\item These will teach them things they will have to learn at some point anyway.
\end{enumerate}
In our experience, most students will be enthusiastic about exercises that include ethics, social justice and politics. But some may be slightly reluctant. Throughout their years at school and in other mathematical courses, these students might have formed a very rigid picture of what counts as proper and good mathematics\footnote{cf. D. M\"uller (2022). Situating ``Ethics in Mathematics'' as a Philosophy of Mathematics Ethics Education, p.3. \href{https://arxiv.org/abs/2202.00705}{https://arxiv.org/abs/2202.00705}}. It is essential that you, as the lecturer or tutorial section leader, are aware of this and teach  the new insight contained in these exercises in an unapologetic and enthusiastic way. The exercises contain exciting mathematics, exciting transferable skills and exciting ethics. They are meant to challenge the students in a deep way without dictating to them any right or wrong ``way of looking at the world''. In our experience, students who understand that they'll have fun learning something important that interests them will be most motivated and will get the most out of these exercises. It is therefore important that you not only motivate your students to do the mathematics, but also motivate and encourage them to see the bigger picture of their work and output. We all like to tell our students how widely mathematics is used: now is the time to also help them understand that not all uses of mathematics are automatically good. 

Similarly, some of your colleagues  may be uncomfortable and express concern when they hear that you plan to incorporate such tools and aspects into your teaching. We have tried to design exercises that do not sacrifice mathematical content: we even dare to suggest that one can produce \textit{better} mathematical exercises by incorporating an ethical lesson or a transferable skill. We encourage you to talk to your colleagues. When they understand that these are genuine, and often hard, mathematical exercises, they'll usually be more open to the idea of embedding ethical training into courses.

Occasionally, you will need to talk about some of the necessary ethical or transferable skills in your lectures or during your tutorial sections. There are some exercises that presume background knowledge that isn't typically presented in undergraduate courses (or high school). For example, one usually cannot expect students to magically know about p-hacking in statistics if it's their first statistics course and they've only just learned the concept of p-values. Similarly, the pillars, and many of the transferable skills taught in these exercises, will be new to most students. It might make sense to either briefly talk about the pillar(s) in class,  or to provide the necessary pillar(s) as an additional handout in order to guide students. We found that a good option is to use some of the exercises as worked examples in lectures, so that students understand some of the core ideas behind how these might be tackled. This also helps you to teach, and moreover normalise, the ethics, transferable skills and background knowledge throughout the entire course, instead of trying to teach all of it upfront. In our experience, it is better if students see 5-10 minutes here and there embedded into worked examples or as a standalone discussion, rather than in a continuous 60 minute block at one single point in the course.

Tutorial leaders will also need some preparation for this sort of teaching. For starters, they will need to be informed of the presence, and purpose, of these questions on the exercise sheets; many tutorial leaders will have never come across the concept of Ethics in Mathematics, let alone taught it. In addition, they will need to be prepared for the fact that many of these questions give rise to a degree of open discussion. This is different from your average mathematical exercises where there is often only one correct answer, and/or it is sufficient to provide only one proof. Ethics is very much governed by a plurality of competing claims\footnote{cf. A. Sen (2010). \textit{The Idea of Justice}, p. 10. London: Penguin}. Students will suggest different arguments, produce different solutions, and come to different conclusions. Indeed, two quite different ethical arguments may contradict each other, and yet both be defensible. It is important to be ready to discuss these questions with an open mind: they are designed to teach ethical awareness, and they explicitly do not try to teach a fixed picture of what is ``right'' or ``wrong''. All this should be outlined as part of the ``tutor induction pack'' that tutors receive, be it printed solutions to the exercise sheets, induction sessions before their teaching begins, teaching guidance notes, etc. 

\subsection{Further Remarks: Teaching Ethics in Mathematics}

\subsubsection{How do these resources fit into teaching Ethics in Mathematics?}

Teaching Ethics in Mathematics is rather challenging, and probably a very new concept and endeavour for most academic mathematicians.

For a quick introduction to Ethics in Mathematics and why is is worth considering, you can consult Ernest\footnote{P. Ernest (2020). The ethics of mathematical practice. \textit{Handbook of the history and philosophy of mathematical practice}, Springer, pp. 1--38. \href{https://doi.org/10.1007/978-3-030-19071-2_9-1}{https://doi.org/10.1007/978-3-030-19071-2\_9-1}} and Chiodo and Clifton\footnote{M. Chiodo \& T. Clifton (2019). The importance of ethics in mathematics. \textit{European Mathematical Society Magazine} 114.\href{https://ems.press/journals/mag/articles/16606}{https://ems.press/journals/mag/articles/16606}}. The discussion paper by Chiodo and Vyas\footnote{M. Chiodo \& R. Vyas (2019). The role of ethics in a mathematical education. \textit{Ethics in Mathematics Discussion Papers} \href{https://www.ethics-in-mathematics.com/pub}{https://www.ethics-in-mathematics.com/pub}} contains further information about why it might be necessary to teach ethics to undergraduate mathematicians. We also want to refer to the special issue on \textit{Ethics in Mathematics} of the Journal of Humanistic Mathematics\footnote{Ethics in Mathematics, special issue of \textit{Journal of Humanistic Mathematics}, volume 12, issue 2 (2022). \href{https://scholarship.claremont.edu/jhm/vol12/iss2}{https://scholarship.claremont.edu/jhm/vol12/iss2}} and the PRIMUS special issue on \textit{Mathematics for Social Justice}\footnote{Mathematics for Social Justice, special issue of \textit{Problems, Resources, and Issues in Mathematics Undergraduate Studies}, volume 29, issue 3--4 (2019).  \href{https://www.tandfonline.com/toc/upri20/29/3-4}{https://www.tandfonline.com/toc/upri20/29/3-4}}, both of which contain further information on the subject matter.

To understand how this resource fits into a more holistic approach to teaching Ethics in Mathematics, you can consult Chiodo and Bursill-Hall\footnote{M. Chiodo \& P. Bursill-Hall (2019). Teaching ethics in mathematics. \textit{European Mathematical Society Magazine} 114. \href{https://doi.org/10.4171/news/114/10}{https://doi.org/10.4171/news/114/10}}, which outlines the three main components needed to successfully carrying out such teaching:
\begin{enumerate}
    \item \textbf{Constructing a standalone course on ethics in mathematics}, to go into the specific details on what it is about mathematical work that leads to ethical issues.
    \item \textbf{Integrating ethics exercises across all courses}, to normalise the perception that there are such issues in mathematical work, and to train such thinking and make it habitual. 
    \item \textbf{Securing support and approval from colleagues}, to create a positive environment for students so that they understand that this forms part of the embedded values of their community.
\end{enumerate}

And to read more about why teaching Ethics to Mathematics students is only \textit{part} of a larger task of raising meaningful ethical awareness and engagement by mathematicians, you can consult M\"uller et al.\footnote{D. M\"uller, M. Chiodo, \& J. Franklin (2022). A Hippocratic Oath for mathematicians? Mapping
the landscape of ethics in mathematics. \textit{Science and Engineering Ethics} 28:5, article 41. \href{https://doi.org/10.1007/s11948-022-00389-y              }{https://doi.org/10.1007/s11948-022-00389-y                        }}.

\subsubsection{Where to go from here?}
If you are intrigued by this and want to go even further, then there are several options for you to take. For example, you might:
\begin{itemize}
    \item design your own exercises using our (or someone else's) methodology,
    \item encourage others to incorporate more ethical awareness into their teaching and courses,
    \item consider teaching a standalone course on ethics in mathematics\footnote{There are many resources available that can help you to design a standalone course on Ethics in Mathematics. A brief course description of each lecture and the full video recording of Dr. Maurice Chiodo's 8 hour seminar course is available on the  Ethics in Mathematics Project website (\href{https://www.ethics-in-mathematics.com/course/lectures}{https://www.ethics-in-mathematics.com/course/lectures}). Other course descriptions and syllabi are also readily available on the internet, for example, Dr. Matthew Cordes's Ethics in Mathematics course at ETH Zurich (\href{https://www.mcordes.com/F23_ethics_syllabus.pdf}{https://www.mcordes.com/F23\_ethics\_syllabus.pdf}) and Dr. Allison N. Miller's mini-seminar on Ethics in Mathematics (\href{https://scholarship.claremont.edu/jhm/vol12/iss2/11}{https://scholarship.claremont.edu/jhm/vol12/iss2/11}). Allison Miller also maintains a substantive reading list on the topic that can serve as a starter for those wishing to design a standalone course: \href{https://docs.google.com/document/d/1P4xAUa9qfEbeCjurGOdJH3dW2JErVRJZ1K9DLAT8uYE/edit}{https://docs.google.com/document/d/1P4xAUa9qfEbeCjurGOdJH3dW2JErVRJZ1K9DLAT8uYE/edit}.}.
    \item Or if you're very ambitious, even try to (re-)design the mathematics curriculum at your institution with ethics in mathematics in mind.
    \item Ethics in mathematics also integrates well with other professional and skill-based courses (e.g., on scientific integrity and publishing ethics), as well as research projects and research experiences for undergraduate students (REUs).
\end{itemize}
If you want to do more, we encourage you to have a look at the reading material in Section \ref{sec:further_reading}. The further reading contains a lot of information and inspiration about what to do next.

\subsubsection{Updates to this document}
This is intended to be a \textit{living document}, with  updates and additions over time. As we develop more of these resources (questions, projects, reading list), we will add them here. If you have any suggestions on content to add to any of these three categories, please feel free to get in touch. To make this a truly universal teaching resource requires contribution from many different mathematicians, with different skills, experience, and perspective.
\newpage
\section{Exercises: Calculus}
\begin{enumerate}\item \textbf{Topic: \textit{Calculus}
    \\Allied disciplines: \textit{Economics, Engineering, Statistics, Natural Sciences}
    \\Year level: 1
    \\\hyperref[list:pillars]{Ethical pillars}: 1,5,9}

    In a large population, the proportion with income between $x$ and $x+dx$ is $f(x) dx$. Express the mean income $\mu$ as an integral, assuming that any positive income is possible. Let $p = F(x)$ be the proportion of the population with income less than $x$, and $G(x)$ be the mean income earned by people with income less than $x$. Further, let $\theta (p)$ be the proportion of the total income which is earned by people with income less than $x$ as a function of the proportion $p$ of the population which has income less than $x$. Express $F(x)$ and $G(x)$ as integrals and thence derive an expression for $\theta(p)$, showing that
	\begin{equation*}
	\theta(0) = 0 ~, ~~ \theta(1) = 1 ~,
	\end{equation*}
	and
	\begin{equation*}
	\theta^\prime (p) = \frac{F^{-1}(p)}{\mu} ~, ~~ \theta^{\prime\prime}(p) = \frac{1}{\mu f(F^{-1}(p))} > 0~.
	\end{equation*}
	Sketch the graph of a function $\theta(p)$ with these properties. Deduce that, if there is any variation in income, the bottom (when ordered in terms of income) proportion $p$ of the population receive less than $p$ of the total income, for all positive values of $p$. Just how much less is quantified by the (in)famous ``Gini index'' beloved of economists, which is twice the area between the curve $\theta(p)$ and the diagonal line connecting $(0,0)$ and $(1,1)$.
	\linebreak
	\linebreak
	A particular population's income is described by the exponential distribution:
	\begin{equation*}
	f(x) = \lambda e^{-\lambda x} ~~~ \text{for} ~~x>0
	\end{equation*}
	for some constant $\lambda > 0$, using your expression for $\theta(p)$, compute the Gini index in this case. Do you think the Gini index is a good measure of inequality of income?

\ifdefined\solution
    \solution \textit{Do you think the Gini index is a good measure of inequality of income? The point of the example given at the end is that the answer is independent of $\lambda$. Do students agree with the index, that the family of exponential distributions all describe the same degree of inequality? This suggests the question of ``how do we define inequality?'' or ``how \textit{should} we define inequality?''}

	We have
	\begin{equation*}
	\mu = \int_{0}^{\infty} t f(t) dt ~, ~~ F(x) = \int_{0}^{x} f(t) dt ~, ~~ G(x) = \frac{1}{\mu F(x)} \int_{0}^{x} t f(t) dt ~.
	\end{equation*}
	Then
	\begin{align*}
	\theta(p) &= F(x)G(x) = \frac{1}{\mu} \int_{0}^{x} tf(t) dt \\
	&= \frac{1}{\mu} \int_{0}^{F^{-1}(p)} tf(t) dt ~.
	\end{align*}
	Noting $F^{-1}(0) = 0$ and $F^{-1}(p) \rightarrow \infty$ as $p \rightarrow 1$, we get
	\begin{equation*}
	\theta(0) = \frac{1}{\mu} \int_{0}^{0} tf(t) dt = 0
	\end{equation*}
	and 
	\begin{equation*}
	\theta(1) = \frac{1}{\mu} \int_{0}^{\infty} tf(t) dt = \frac{\mu}{\mu} = 1 ~.
	\end{equation*}
	Differentiating with respect to $p$, we get
	\begin{align*}
	\theta^\prime (p) &= \left( \frac{dp}{dx}\right)^{-1} \frac{d\theta}{dx} \\
	&= \frac{1}{F^\prime (x)} \frac{1}{\mu} \frac{d}{dx} \left(\int_{0}^{x} tf(t) dt\right) \\
	&= \frac{1}{\mu f(x)} xf(x) = \frac{x}{\mu} = \frac{F^{-1}(p)}{\mu} ~,
	\end{align*}
	and
	\begin{align*}
	\theta^{\prime\prime} (p) &= \left( \frac{dp}{dx}\right)^{-1} \frac{d}{dx} \left(\frac{x}{\mu}\right) \\
	&= \frac{1}{f(x)} \frac{1}{\mu} = \frac{1}{\mu f(F^{-1} (p))} ~,
	\end{align*}
	where each is $>0$ since both numerators and denominators are positive.
	\linebreak
	\linebreak
	\{Sketch: any strictly convex function $\theta (p)$ with $\theta (0) = 0$, $\theta^\prime (0) = 0$, $\theta (1) = 1$, and $\theta^\prime (p) \rightarrow \infty$ as $p \rightarrow 1$.\}
	
	In the case $f(x) = \lambda e^{-\lambda x}$, we have
	\begin{align*}
	\mu &= \int_{0}^{\infty} \lambda t e^{-\lambda t} dt \\
	&= \left.\left[ -t e^{-\lambda t} \right]\right._{0}^{\infty} + \int_{0}^{\infty} e^{-\lambda t} dt \\
	&= \left.\left[ -\frac{1}{\lambda} e^{-\lambda t} \right]\right._{0}^{\infty} = \frac{1}{\lambda} ~.
	\end{align*}
	Thus we find the Gini index $g$ as
	\begin{align*}
	g &= 2\left( \frac{1}{2} - \int_{0}^{1} \theta(p) dp
	\right) = 1 - 2 \int_{0}^{\infty} \theta (F(x)) \frac{dp}{dx} dx \\
	&= 1 - \int_{0}^{\infty} \left( 2\lambda \int_{0}^{x} \lambda t e^{-\lambda t} dt \right) \lambda e^{-\lambda x} dx \\
	&= 1 - \int_{0}^{\infty} \left( 2\lambda \left. \left[ -t e^{-\lambda t} \right]\right._{0}^{x} + \int_{0}^{x} 2\lambda e^{-\lambda t} dt \right) \lambda e^{-\lambda x} dx \\
	&= 1 - \int_{0}^{\infty} \left( -2\lambda x e^{-\lambda x} +  \left. \left[ -2 e^{-\lambda t} \right]\right._{0}^{x} \right) \lambda e^{-\lambda x} dx \\
	&= 1 - \int_{0}^{\infty} \left( 2 - 2 e^{-\lambda x} - 2\lambda x e^{-\lambda x} \right) \lambda e^{-\lambda x} dx \\
	&= 1 - \int_{0}^{\infty} \left( 2\lambda e^{-\lambda x} - 2\lambda e^{-2\lambda x} - 2\lambda^2 x e^{-2\lambda x} \right) dx \\
	&= 1 - \left(  \left. \left[ -2 e^{-\lambda x} + e^{-2\lambda x} + \lambda x e^{-2 \lambda x} \right] \right._{0}^{\infty} - \int_{0}^{\infty} \lambda e^{-2\lambda x} dx \right) \\
	&= 1 - \left( 2 - 1 + \left. \left[ \frac{1}{2} e^{-2 \lambda x} \right] \right._{0}^{\infty}\right) \\
	&= 1 - \left( 2 - 1 - \frac{1}{2}\right) = \frac{1}{2} ~.
	\end{align*}

    The Gini index has many drawbacks. Primarily, it is a relative measure. If income over the entire population decreases, the Gini index may stay the same. Conversely, if income increases only for the top earners but remains the same for low earners, the Gini index will increase. It does not measure absolute quality of life, just the difference between top and bottom earners.
    \\In the example just given, there is an entire family of income distributions, $f(x) = \lambda e^{-\lambda x}$, whose Gini index is identical and independent of the parameter $\lambda$. It would be a useful exercise to explore the effect of net income for small, and large, values of  $\lambda$ for this family of distributions, to see the outcomes that the Gini index fails to pick up.

    \fi
\end{enumerate}

\newpage
\section{Exercises: Introduction to Proofs}
\begin{enumerate}
    \item \textbf{Topic: \textit{Propositional Logic}
    \\Allied disciplines: \textit{Computer Science}
    \\Year level: 1-2
    \\\hyperref[list:pillars]{Ethical pillars}: 5}
    
    \begin{enumerate}
        \item Consider the statement: `If you don't do it, then someone else will.' 
        
        Express this statement in symbolic notation. Find its contrapositive and its negation, giving each both in symbols and in words. 
    
        \item Your boss has given you a task. The task is well within your technical capability, but you are not sure whether it would be legal or ethical. You feel uneasy, but your boss tells you: `If you don't do it, then someone else will.'
    
        Do you think that the boss' argument is cogent? Does that depend on who you are or what the project is? How would you answer your boss?

    \end{enumerate}
    
    \ifdefined\solution
        \solution \textit{The point of this question is to get the students to show, mathematically, that the paradigm ``If you don't do it, someone else will'' is quite flawed as a argument or incentive to do something. Hopefully this can equip the students with enough realisation to push back on such an argument, should they ever encounter it in their professional life.}

       \begin{enumerate}
        \item  In symbols, the statement is
            $$ \neg P(x) \implies \exists y : (y\neq x) \cap P(y), $$
        where $P(y)$ means `$y$ will do it', and you are $x$.
        
        The contrapositive is
            $$ \neg\left(\exists y :( y\neq x) \cap P(y) \right) \implies P(x), $$
        or, simplifying the left-hand side,
            $$ \forall y: (y=x) \cup \neg P(y) \implies P(x). $$
        In words, the contrapositive says: `If nobody else does it, then you will.' Note that `if-then' clauses in English usually imply some sort of \emph{causal} relationship; this is not the case in mathematics.
        
        The negation is
            $$ \neg P(x) \cap \neg \left(\exists y : y\neq x \cap P(y) \right), $$
        which simplifies to 
            $$ \neg P(x) \cap \forall y: y=x \cup \neg P(y). $$
        In words, `You won't do it, and nobody else will'. 
        
        Note that \emph{we have not specified the domain of discourse}. In the present context, the domain might be the set of people who would be capable of doing the task: if you are the only member of that set, then $\nexists y : y\neq x$.

        As a specialist, you have power, and you can use that to answer and influence your boss. 
        \item This argument ``If you don't do it, someone else will'' is flawed because it \textit{assumes} that the activity will happen. But that it not a given. 
        \\Compare this to the scenario where there is 1 red sock and 9 blue socks in a bag. 10 people are each invited to come a pick out a sock. Then the statement ``If you don't pick the red sock, someone else will'' is true, IF all 10 people take up the offer to pick out a sock (because all scenarios are accounted for). But if not all 10 pick out a sock, then ``If you don't pick the red sock, someone else will'' is no longer true. What is true instead is that ``If you don't pick the red sock, someone else \textit{might}''. Or they might not.
        
        So the statement ``If you don't do it, someone else will'' \textit{assumes} that at least one person on earth will carry out the task. But this is not necessarily true. Perhaps you are the only person capable of doing it, in which case it is clearly false. Or perhaps the people capable of doing it (eg: your work colleagues) all refuse, in which case if you don't do it, it will not be done. 
        
        The (correct) statement that your boss can give is ``If you don't do it, someone else \textit{might}'' So now we are talking about probabilities. If the task is very easy (say, over 1 billion people are capable of doing it) and many of those have the incentive to do so, then the ``might'' is with high probability. But if very few people are capable of doing it, and/or very few have the incentive to do it, then the ``might'' is with very low probability. Now think of the sorts of tasks you may be asked to do as a mathematician, how many other people could do it, and of those how many others are incentivised to do it. The ``might'' in that statement could be very small indeed. So it could well be ``If you don't do it, someone else \textit{might} with probability 0.0000000001'', so your refusal to do it could overwhelmingly influence the chance of it being done. 
        \end{enumerate}
\fi
        
\end{enumerate}
\newpage
\section{Exercises: Analysis}

\begin{enumerate}
\item \textbf{Topic: Analysis
    \\Year level: 1
    \\\hyperref[list:pillars]{Ethical pillars}: 2,5,6}
    
    This question asks you to consider what is natural about the natural numbers $\mathbb{N}$. Briefly recall the Peano axioms which we have used in our lectures to construct $\mathbb{N}$. Let $\mathbb{N}$ be a set satisfying 
    
    \begin{enumerate}
        \item $\mathbb{N}$ contains a special element which we call $1$.
        \item There exists a injective map $\sigma: \mathbb{N}\rightarrow \mathbb{N} \setminus \{1\}$, $\sigma(n) \neq n$ for all $n \in \mathbb{N}$. 
        \item For every subset $S\subset \mathbb{N}$ such that $1\in S$ and if $n\in S$, then $\sigma(n)\in S$, it follows that $S = \mathbb{N}$.
    \end{enumerate}
    Which (if any) of $(a), (b)$ and $(c)$ feel natural to you, and why? Now consider that there are indigenous tribes that potentially count differently. Some only have the language for the first few numbers, and then use ``many'' for every larger set. For example, they might count out loud $1,2,3,4,5, \text{ many}$. In your opinion, which of $(a), (b)$ and $(c)$ would feel natural to them? Do you see why we call it ``Peano's axioms'' and not ``Peano's Theorem''? 
    \\(Hint: Do Peano's axioms work for a finite set? Do you think that someone's vocabulary fully determines their ability to count?) 
    \newline \newline You can find a brief discussion of this phenomenon in an article by Butterworth\footnote{B. Butterworth (2004). What happens when you can't count past four?. \textit{The Guardian}. \href{https://www.theguardian.com/education/2004/oct/21/research.highereducation1}{https://www.theguardian.com/education/2004/oct/21/research.highereducation1}}.

    \ifdefined\solution
    \solution
    \textit{This question is designed to make students aware that what seems natural to you, might not be natural to someone else. It uses recent research from ethnomathematics to show students that even the basics of what we perceive as pure mathematics and counting were constructed in a social context and that there is an ongoing scientific discussion about the role of language in mathematics and counting. It also teaches students that the Peano axioms don't work for a finite set of numbers. There cannot be a bijective map from a finite set into a strictly smaller set. \\
    Additional commentary detailing the politics of this and other questions can be found in M\"uller and Chiodo\footnote{D. M\"uller, M. Chiodo (2023). \textit{Mathematical Artifacts have Politics: From Examples to Embedded Ethics.} \href{https://arxiv.org/abs/2308.04871}{https://arxiv.org/abs/2308.04871}} }
    
\begin{itemize}
    \item Which of the items (a), (b) or (c) feels natural to the student can differ. From our experience, (a) feels very natural (having the number $1$ is non-controversial; if you exchange $1$ for $0$, this might already be a different story). (b) and (c) often feel natural if the student properly understands the definition of the Peano axioms, but they can also feel unnatural in the sense that it is a highly abstract way of ``adding $1$''.  
    \item Students should understand that (b) fails if you only have a finite set of numbers, since you cannot have a bijection from a finite set into a strictly smaller subset. (c) also fails in the context of \{1,2,3,4,5 many\}. In this case, one has ``many+many = many'', ``3+4=many'', etc. In this case, one could argue that only (a) would feel natural to someone who only has the set \{1,2,3,4,5, many\} available. But students should also understand that the situation is much more complicated than this. 
\end{itemize}
    
     \fi
\newpage
\item \textbf{Topic: \textit{Analysis}
\\Year level: 1-2 
\\\hyperref[list:pillars]{Ethical pillars}: 2,5} 

\begin{enumerate}
    \item Write down what might characterise a mathematical proof as ``ugly'', and as ``nice''. Is every proof either ugly or nice? Can a proof be both ugly, and nice, simultaneously?
    \item Let $f_n: [0,1]\rightarrow[0,1]$ be a continuous function for each $n\in\mathbb{N}$. Show that their pointwise maximum $h_m(x) = \max\{f_1(x), f_2(x), \dots, f_m(x)\}$ is continuous on $[0,1]$ for each $m$. Give two proofs of this result, that fit your descriptions of `ugly' and `nice' respectively.
    \item Compare your proofs to those provided by other students, and see if theirs match your descriptions of ugly and nice.
    \item Must the function $h$ defined by $h(x) = \sup\{f_m(x): m\in\mathbb{N}\}$ be continuous on $[0,1]$?
\end{enumerate}

\ifdefined\solution
\solution
\textit{The purpose of this question is to prompt the students to think about, and realise, that mathematics is viewed beyond ``correct vs incorrect''; we have our own stylistic and professional norms and beliefs on what is ``good'' mathematics, and what is ``bad'' mathematics, which go so far as to reflect the aesthetics of mathematics as ``ugly'' and ``nice''. This perception of mathematical aesthetics might influence how mathematicians model or mathematise a scenario; they may choose a model to use, or curve to fit, because it makes the equations and mathematics ``nice'', rather than for reasons of accuracy or meaningfulness.}

\begin{enumerate}
    \item Students may list all sorts of descriptors for ugly and nice. Here are a few common ones, though this list is by no means exhaustive (and as an instructor, you could present these lists to students and ask if they agree with them):
    \\\textbf{Ugly}: Long. Messy. Requires a lot of calculation. Uses a computer for part of the proof. Uses unmotivated constructions or examples/counterexamples. Breaks the problem into  the treatment of many sub-cases. Combines many disparate ideas.
    \\\textbf{Nice}: Short. Succinct. Has one key idea. Combines some disparate ideas in a novel way. Can be generalised.
    \\Notice how some of the descriptors for ugly/nice above are very similar. eg: `Combines many disparate ideas' (for ugly) and `Combines some disparate ideas in a novel way' (for nice). So some might see a proof as nice, whereas others might see it as ugly.
    \\It is unlikely that students will say that all proofs are either ugly or nice. It would be worth asking them for a proof that is neither ugly nor nice.

\item For a `nice' proof, consider the following:	
\\For $n = 1$, the result holds by definition of $f_1 (x)$. For $n = 2$, if $f(x)$ and $g(x)$ are continuous functions, then
	\begin{itemize}
		\item proposition: $\max\lbrace f(x), g(x)\rbrace$ is continuous
		\\
		proof:
		\begin{equation*}
		\max\lbrace f(x), g(x)\rbrace = \frac{1}{2} (f(x) + g(x)) + \frac{1}{2}|f(x) - g(x)|
		\end{equation*}
		(which we verify by inspection), and this is continuous since it is a sum of continuous functions $\frac{1}{2} f(x)$, $\frac{1}{2} g(x)$, $\frac{1}{2} |f(x) - g(x)|$, where the last is continuous by composition of continuous functions ($|x|$ is continuous). 
	\end{itemize}
	This proves the result for $n = 2$. For $n > 2$, we note that
	\begin{align*}
	h_n (x) &= \max \lbrace h_{n-1} (x), f_n (x)\rbrace
	\end{align*}
	and then the result follows for all $n$ by induction, using the above proposition.
    \\This proof might be considered nice because it starts with showing a straightforward but clever way to write the maximum of two functions as a continuous combination of the original functions (the verification of which takes 2 lines and a known background result that $|x|$ is continuous). It then makes a simple inductive argument to give the overall result. There is no case-bashing to deal with, and the entire proof fits in about 10 lines with only one `sub result'. The proof is possible by the clever observation of the way to write the maximum of two functions, and doesn't require much additional algebra after that.
	\\\mbox{}
    \\For an `ugly' proof, one might start with the $n=2$ step of a similar inductive argument above, but then try to define the values $x_{i}$ for which $f_{1}(x_{i})=f_{2}(x_{i})$. But then students will need to argue things like:
    \begin{enumerate}
        \item The $x_{i}$ can be partitioned into 1) a collection of closed intervals, and 2) some discrete points.
        \item There are only finitely many such intervals, and finitely many such discrete points.
        \item The function $h_2(x) = \max\{f_1(x), f_2(x)\}$ is continuous at each such discrete $x_{i}$, and continuous on each such interval.
    \end{enumerate}
    This is possible to complete, but it is indeed an ugly proof. It starts off by trying to prove the result in the way that students might initially be visualising it (looking for the points where the functions `swap' being largest in the collection). But then there are many sub-cases to deal with, and many many results to verify in order to properly make use of the sub cases. It requires a lot to `follow the proof' and keen in one's mind what has been done, what is being done, and what is yet to be done. It would probably require over a page of algebra to write this out fully, and there is no one core observation that lies at the heart of the proof.
    
	\item $h(x) = \sup\lbrace f_n (x) : n\in \mathbb{N}\rbrace$ need not be continuous on $[0,1]$: consider
	\begin{equation*}
	f_n : [0,1] \rightarrow [0,1] ~, ~~ f_n (x) = x^{\frac{1}{n}}
	\end{equation*}
	which has $f_n (0) = 0$, $f_n (1) = 1$, and $f_{n+1} (x) > f_n (x)$ for $x \in (0,1)$ for all $n \in \mathbb{N}$.
	\\ 
	We thus have $h (0) = 0$, but for all $x \in \left(0,1\right]$ we have $h (x) = \lim\limits_{n \rightarrow \infty} x^{\frac{1}{n}} = 1$.
	\\ 
	Hence this $h$ is discontinuous at $x = 0$.
\end{enumerate}
 
\fi

\newpage
\item \textbf{Topic: Analysis
\\Year level: 1-2
\\\hyperref[list:pillars]{Ethical pillars}: 6} 

You are an intern at an engineering firm. Your boss gives you two complicated convex functions $g(x)$ and $h(x)$. She then asks you if $f(x) = \max\left(g(x), h(x)\right)$ has a unique minimum. You remember that convex functions have a unique minimum (if it exists), so you set out to show that $f(x)$ is convex.
\begin{enumerate}
    \item Show that a function is convex iff its epigraph is convex. The epigraph of $f$ is defined as $epi(f) = \{ \left(x,y\right) : x \in dom(f) \& y \geq f(x) \}$.
    \item Show that the intersection of two convex sets is convex.
    \item Is the maximum of two convex functions always convex?
    \item Your boss is very busy this morning and she asks you if you can quickly explain to her why $f(x)$ has a unique minimum. How do you use $(a)-(c)$ to quickly convince her?
    \item Your boss now asks you to prepare a half-page presentation describing the work you've just done, to be used in a sales brochure targeted at CEOs and high-level procurement agents advertising the services your firm can provide. What would you prepare?

\end{enumerate}

\ifdefined\solution
\solution \textit{This question is designed to bridge the (truth) requirements and formalities of pure mathematics and the requirements and constraints of an applied context. It teaches students to make a selection about what information they need/want to talk about or present, and how they might present it, when under external uncontrollable pressures and constraints, such as time limitations, or perhaps limitations in the mathematical knowledge of the audience. It also introduces student to the challenge of explaining higher-level mathematics to an audience that is vastly under-trained to understand it fully.}
\\ \\Recall: A function is convex if its domain is convex and for all $\lambda \in [0,1]$ and all $x,x' \in dom(f)$ we have $f(\lambda x + (1-\lambda)x') \leq \lambda f(x) + (1-\lambda) f(x').$ A set $C$ is convex if for all $a,b \in C$ and all $\lambda \in [0,1]$, we have $\lambda a+ (1-\lambda)b \in C$.
\begin{enumerate}
    \item Let $f$ be convex and let $(x,y), (x',y') \in epi(f)$. Let $\lambda \in [0,1].$ Then $\lambda x + (1-\lambda)x' \in dom(f)$ because $dom(f)$ is convex. Furthermore, $f(x) \leq y$ and $f(x') \leq y'$ by definition of $epi(f)$. Using the convexity of $f$, we obtain
    $$f(\lambda x + (1-\lambda)x') \leq \lambda f(x) + (1-\lambda) f(x') \leq \lambda y + (1-\lambda) y'.$$
    In particular, $(\lambda x + (1-\lambda)x', \lambda y + (1-\lambda) y') \in epi(f).$ Thus, $epi(f)$ is convex. \\ \\ Conversely, if $epi(f)$ is convex, then by definition $dom(f)$ must already be convex. We know that $(x,f(x)), (x', f(x')) \in epi(f)$. By the convexity of $epi(f)$, we then know that $\left(\lambda x + (1-\lambda)x', \lambda f(x) + (1-\lambda)f(x' ) \right) \in epi(f)$, but by definition of $epi(f)$ this exactly means that
    $$f(\lambda x + (1-\lambda)x') \leq \lambda f(x) + (1-\lambda) f(x').$$
    Thus, $f$ is convex.
    
    \item Let $C_1, C_2, C_3, \dots$ be convex sets. Let $a,b\in C_i$ for all $i = 1, 2, 3, \dots$ and let $\lambda \in [0,1]$. Then it follows from the convexity of $C_i$ that $\lambda a + (1-\lambda)b \in C_i$ for all $i = 1, 2,3, \dots$. Thus,
    $\lambda a + (1-\lambda)b \in C_i \in \cap_{i} C_i$. Thus, the intersection of convex sets is convex.
    \item The epigraph of the maximum of two functions is the intersection of the epigraphs. The maximum of two convex functions is therefore convex. Both $g$ and $h$ are convex, hence from it immediately follows that $f$ is convex.
    \item The student should realise that the lengthy mathematical derivation can be skipped by drawing a few pictures. All individual results in $(a)-(c)$ are easy to understand based on the corresponding pictures alone. Students should learn that depending on the situation more or less mathematical precision might be necessary to convince someone. As a mathematician, they will often be required to understand the formal proof, but in applied situations explaining a result to someone else does not always require the entire proof. 
    \newline Note to instructor: (d) can be changed to ask the students to produce a single slide that explains why $f$ is convex.
    \item Students will first need to identify that the target audience (CEOs and high-level procurement agents) probably have very little mathematical knowledge, unlike their boss who probably has rather substantial mathematical knowledge. Upon realising this, the problem probably becomes quite hard for students. They may have experience explaining mathematical concepts to people who have the (approximate) requisite knowledge to understand them fully. But explaining mathematical concepts to people with insufficient knowledge to understand them fully might be an unfamiliar activity for students. They will need to identify the core results to be explained, understand what detail can be left out, reformulate the terminology for a lay audience, and then still have `enough left' to say something meaningful/useful. They will need to think very carefully about \textit{why} half a page has been dedicated to their work in the brochure.

\end{enumerate}

\fi

\newpage
\item 
\textbf{Topic: Analysis
\\Year level: 1-2
\\\hyperref[list:pillars]{Ethical pillars}: 5,6} 

    Let $f: D\to \mathbb{R}$ be a function from $D \subseteq \mathbb{R}$ to $\mathbb{R}$. 
    
    We say that $f$ is \textit{Type 1 Continuous} (Con1)  on $D$ if, $\forall a \in D$:
    \[
    \forall \epsilon>0 \ \exists \delta>0 \textnormal{ such that } \big(|x-a|<\delta\big) \Rightarrow \big (x \in D \textnormal{ and } |f(x)-f(a)|<\epsilon\big )
    \]
      We say that $f$ is \textit{Type 2 Continuous} (Con2)  on $D$ if, $\forall a \in D$,
      \[
    \big( \exists b<c \textnormal{ with } a \in [b,c]\subseteq D \big)\textnormal{ and } \big( \forall \epsilon>0 \ \exists \delta>0 \textnormal{ such that } \big(|x-a|<\delta \textnormal{ and } x \in D) \Rightarrow \big ( |f(x)-f(a)|<\epsilon\big ) \big) 
    \]
    We say that $f$ is \textit{Type 3 Continuous} (Con3) on $D$ if, $\forall a \in D$:
    \[
    \forall \epsilon>0 \ \exists \delta>0 \textnormal{ such that } \big(|x-a|<\delta \textnormal{ and } x \in D) \Rightarrow \big ( |f(x)-f(a)|<\epsilon\big )
    \]
  
     We say that $f$ is \textit{Type 4 Continuous} (Con4)  on $D$ if, $\forall a,b \in D$ with $a<b$:
    \[
     \big(f(a)<y<f(b) \textnormal{ or }f(a)>y>f(b) \big)  \Rightarrow \big(\exists c \in D \textnormal{ such that } a<c<b \textnormal{ and }f(c)=y\big)
    \]
    
    Consider the following functions:
        \begin{align*}
            f_{1}&:[0, \infty) \to \mathbb{R} \ ; \ f_{1}(x):=0
            \\f_{2}&: \mathbb{Q}\to \mathbb{R} \ ; \ f_{1}(x):=0
            \\f_{3}&: \mathbb{R}\to \mathbb{R} \ ; \ f_{3}(x):= 
             \begin{cases} 
                 0 \hspace{1.1cm} \textnormal{ if } x=0  \\
                 \sin\big(\frac{1}{x}\big) \ \  \textnormal{ if } x \neq 0  
             \end{cases} 
        \end{align*}

    \begin{enumerate}
        \item How would you describe the \textit{idea} or \textit{motivation} behind a ``continuous function''?
        \item Which of these types of continuity do $f_{1}, f_{2}, f_{3}$ satisfy on their domain of definition? Give reasons.
        \item Which of these is the \textit{right} definition of continuity, if any? Why? Who should decide what definition we use?
        \item When defining a mathematical object, be it in an abstract setting like this or a more applied setting, how do we establish the ``right'' definition, if at all? How do we choose what \text{ideas} we want to reflect, and how do we ensure that our definition captures the things we want it to capture, and avoids the things we want it to avoid?

    \end{enumerate}

    \ifdefined\solution
\solution
\textit{The purpose of this question is to help students see several things:
\\1. There is often no straightforward, obvious, or ``natural'' way to define a mathematical object, 
\\2. Many suitable contenders may exist, each with their own characteristics and consequences, and
\\3. Mathematicians don't necessarily agree on which one is ``right'', and it may take many years, or centuries, to come to such a consensus. (The definition of continuity has taken around $200$ years to achieve consensus on).
\\4. The intentions we had in our mind when making our definition my not always be fully realised in the way we had hoped or planned. (This manifests with the Alignment Problem in AI development).
\\This is not simply an artefact of ``abstract mathematics''. Students will often go and work in industrial or more practical roles, and make such ``arbitrary'' definitions and choices all the time. Mathematics may be precise, but its formulation is far from deterministic or ``provable''.
\\It is crucial that the students actively engage with the latter parts of the problem, and think about whether the \textit{right} definition is something that always, or even ever, exits, and what the intended and unintended consequences of their definitions might be. This is seen in industrial applications, such as when mathematics is used to define/measure things such as credit worthiness (with a computed credit score), or as part of the \textit{AI alignment problem} of ensuring that an AI system carries out all the actions that the developers hoped for, and none of the actions that they didn't want to occur.}

\begin{enumerate}
    \item This aims to find what mental image students have when they think of continuity. Some may say things like ``can be drawn on a page without taking your pen off the paper'', or ``a curve with no breaks in it'', and so on. All of these are perfectly reasonable concepts that mathematicians have tried to capture with their formal definition of ``continuous.''
    \item Let us look at each of Con1-Con4 in turn.
    \\\mbox{}
    \\Con1:
    \\For $f_{1}$: $0 \in [0, \infty)$, but given any $\delta>0$ we have $|0-\frac{\delta}{2}|=|\frac{\delta}{2}|<\delta$, but $-\frac{\delta}{2}\notin [0, \infty)$. So $f_{1}$ is not Con1.
    \\For $f_{2}$: $0 \in \mathbb{Q}$ but given any $\delta>0$ there is some irrational $r\in \mathbb{R}\setminus \mathbb{Q}$ with $0<r<\delta$, and so we have $|0-r|=|r|<\delta$, but $r\notin \mathbb{Q}$. So $f_{2}$ is not Con1.
    \\For $f_{3}$: $0 \in \mathbb{R}$, but given any $\delta>0$ we can find $0<x<\delta$ with $\sin(\frac{1}{x})=1$. Namely, find $k \in \mathbb{N}$ with $\frac{1}{\delta}\in[(2k-2+\frac{1}{2})\pi,(2k+\frac{1}{2})\pi]$, and then choose $x=\big((2k+2+\frac{1}{2})\pi\big)^{-1}$. Thus, if $\epsilon=\frac{1}{2}$ was chosen, there is no such $\delta$ satisfying the requirement. So $f_{3}$ is not Con1.
    \\\mbox{}
    \\Con2:
    \\For $f_{1}$: Given any $a \in [0, \infty)$, we have $a \in [a,a+1]\subseteq  [0, \infty)$. Now, given any $\epsilon >0$ simply take $\delta=1$; then taking any $x\in [0, \infty)$ with $|x-a|<1$ we have $|f(x)-f(a)|=0<\epsilon$. So $f_{1}$ is Con2.
    \\For $f_{2}$: There is no closed interval $[b,c]\subseteq \mathbb{Q}$, for any $b<c$. So $f_{2}$ is not Con2.
    \\For $f_{3}$: The same reasoning applies as for Con1, so $f_{3}$ is not Con2.
    \\\mbox{}
    \\Con3:
    \\For $f_{1}$: This is just continuity ``in the usual sense'', so $f_{1}$ is Con3.
    \\For $f_{2}$: This is just continuity ``in the usual sense'', so $f_{2}$ is Con3.
    \\For $f_{3}$: The same reasoning applies as for Con1, so $f_{3}$ is not Con3.
    \\\mbox{}
    \\Con4:
    \\For $f_{1}$: As $f_{1}$ is ``continuous in the usual sense'' on a real interval, then by the Intermediate Value Theorem it is Con4.
    \\For $f_{2}$: We cannot use the intermediate value theorem here, because $\mathbb{Q}$ contains no intervals. However, the initial test condition of the definition is never satisfied, as $f_{2}=0$ on its entire domain so we can never have $a<b\in D$ with $f(a)\neq f(b)$. Thus it holds vacuously.
    \\For $f_{3}$: A function that satisfies Cond4 is known as a \textit{Darboux function}, and $f_{3}$ (known as the \textit{Topologist's sine curve}) is one such example. So, given $a,b \in D$ with $a<b$, we split our reasoning into two cases:
    \\Case 1: $0 \notin [a,b]$. In this case we can restrict $f_{3}$ to $[a,b]$, on which it is ``continuous in the usual sense'', and thus apply the Intermediate Value Theorem.
    \\Case 2: $0 \in [a,b]$. In this case, at least one of $a,b$ is not 0; wlog say $b$ (a similar proof works if we assume $a$). Whatever the values of $a,b$, we know that $f_{3}(a), f_{3}(b)\in [-1,1]$.  So take our $y$ satisfying $f(a)<y<f(b)$ or $f(a)>y>f(b)$; this implies $y \in [-1,1]$. So now we find $k \in \mathbb{N}$ with $\frac{1}{b}\in(2k\pi,2(k+1)\pi]$, and then choose $z=(2(k+2)\pi)^{-1}$. Then, on the interval $[\frac{1}{2(k+2)\pi},\frac{1}{2(k+1)\pi}]$,  $f_{3}$ is ``continuous in the usual sense'' and takes all values in $[-1,1]$ (and thus value $y$). So we can apply the intermediate value theorem here and find $c$ satisfying $a<\frac{1}{2(k+2)\pi}\leq c\leq\frac{1}{2(k+1)\pi}< b$ with $f_{3}(c)=y$.
    \\So $f_{3}$ is Con4.

    \item Most students might argue that ``The definition we see in lectures/textbook'' is the right one. As 200 years of history of mathematics taught us, there is certainly no \textit{obvious} right definition. A definition was chosen, and now all the machinery of modern mathematics has been manufactured to fit around that definition. Books and papers have been written, that definition is fed into countless proofs, etc. It was the \textit{most widely adopted} definition, but it may not be the right one. Betamax video players  were unable to displace VHS, Dvorak keyboards were unable to displace Qwerty, and there are many (adjacent) countries that use different railway track gauges.
    \\As for ``why'' it is the right decision, they may argue ``because it captures all the functions we want to be continuous, and none of the ones that we don't''. Of course, their perspective might already be tainted with circular reasoning here; they may ``want'' functions to be continuous/ not continuous because they were trained to like that definition.
    \\As for \textit{who} should decide what definition is used, they may make reference to concepts like ``professors'', ``award winners'' etc, which are human power constructs and nowhere near objective/impartial entities. Or they may go the other way, and take the libertarian approach of ``any mathematician can make any definition they like''. This overlooks the fact that mathematics, even academic mathematics, does not operate in a vacuum. If one tries to publish theorems about definitions no-one else finds interesting, then one never publishes. More seriously, a mathematician working in an industrial setting needs to choose relevant and meaningful objects/objectives/aims to define, or their work becomes meaningless (at best), and extremely harmful at worst.
    \\For some details about the history of the definition of continuous functions, see Harper\footnote{
     J.F. Harper (2016). Defining continuity of real functions of real variables, \textit{BSHM Bulletin: Journal of the British Society for the History of Mathematics}, 31:3, 189-204. \href{https://doi.org/10.1080/17498430.2015.1116053    }{https://doi.org/10.1080/17498430.2015.1116053    }}.

    \item This is \textbf{hard}; it forms (part of) the alignment problem for mathematics and technology. The problem is currently extensively discussed for AI in the literature. The point here is not to get students to articulate ``the right way to define mathematical objects'', but rather for them to see how difficult it actually is, in an abstract setting, or applied setting. Setting objective/reward/cost functions in machine learning is making definitions. Mathematical modelling of a system (physical, social, etc) is making definitions. There is somewhat of a spectrum between ``making a mathematical definition'' and ``making a mathematical decision''. Defining continuous functions, or bell curves, are definitions. Defining ``the curve that models phenomenon X'' is half-way between definition and decision. Choosing which data to use in a machine learning process is making a decision. But all of these are the same process: mathematicians choosing what to do, and how to go about it. Everything follows logically \textit{after} the choice is made, but those choices are critical in determining what happens after.
   
\end{enumerate}

\fi

\end{enumerate}
\newpage
\section{Exercises: Linear Algebra}

\begin{enumerate}
\item

\textbf{Topic: \textit{Matrices}
 \\Allied disciplines: \textit{Economics, Engineering, Statistics, Computer Science, Natural Sciences}
\\Year level: 1
\\\hyperref[list:pillars]{Ethical pillars}: 5,7,8}

\begin{enumerate}
\item  A square matrix with entries in $\mathbb{R}$ is said to be \textit{column-stochastic} if all of its entries are nonnegative and the entries in each column sum to one. Show that every column-stochastic matrix has 1 as an eigenvalue.
	 
	 \item The Google PageRank algorithm (simplified) works as follows: Each webpage $w_{i}$ on the web is assigned a value $v_{i}$ such that, if $L_{i}$ is the set of pages that link to $w_{i}$, and $n_{i}$ is the number of outgoing links from page $w_{i}$, then
	\[
	 v_{i}=\sum_{w_{j}\in L_{i}}\frac{v_{j}}{n_{j}}
	\]
	That is, each page $w_{i}$ ``donates'' $\frac{1}{n_{i}}$ of its value $v_{i}$, uniformly, to each page that it links to, and then the pages are ``ranked'' by simply ordering the $v_{i}$'s. Show that, if every page on the web links to at least one other page, then there is at least one way of assigning values to each webpage that satisfies the above relation.

 \item What are some effects of choosing this, or indeed any, ranking algorithm on a system as widespread as Google Search?

\end{enumerate}

\ifdefined\solution
    \solution\textit{This question serves two purposes: \\First, it shows that a relatively straightforward piece of first-year mathematics can be used to have a huge impact on the way the world functions; the Google pagerank algorithm is (or more accurately, was, as it has now been updated) one of the most influential algorithms in the world. It was often referred to as the \$100,000,000,000 algorithm (the market value of Google at the time). Not every high-impact use of mathematics requires PhD-level content.
    \\Second, there may well be a whole suite of algorithms that could ``work'' in this scenario; how does Google decide which one to use? And how does this decision introduce bias into the search results; just because a computer does it, doesn't mean it is impartial!}

    \begin{enumerate}
        \item  For a square matrix  $A$, the eigenvalues of $A$ are the same as the eigenvalues of the transpose $A^{T}$. Now, the transpose of a column-stochastic matrix is a row-stochastic matrix (rows sum to 1), and such a matrix has as an eigenvector the column vector with $1$ in every entry, with corresponding eigenvalue $1$.

        \item Form a square matrix $A$ where in column $i$ we put the value $1/n_{i}$ in each entry $(i,j)$ where $w_{i}$ links to $w_{j}$, and 0 in all other entries; there should be $n_{i}$ such nonzero entries in each column, summing to 1. We now need to find a column vector $x$ of ``values'' such that $Ax=x$. That is, if there are a total of $m$ webpages on the internet, then we are seeking $\{v_{1}, \ldots, v_{m}\}$ such that
        
        \begin{align}
    A \begin{bmatrix}
           v_{1} \\
           v_{2} \\
           \vdots \\
           v_{m}
         \end{bmatrix}
         =
         \begin{bmatrix}
           v_{1} \\
           v_{2} \\
           \vdots \\
           v_{m}
         \end{bmatrix}
  \end{align}
        
        So we are looking for an eigenvector of $A$ with eigenvalue 1. Now apply the previous result to show that at least one such eigenvector exists.

\item Using a human-designed/chosen algorithm on a search engine as large as Google poses a serious (ethical) question here: which algorithm should be ``chosen''? How do we know what the benefits and drawbacks are of such an algorithm when deployed on a search engine that billions of people use every day? This may seem like a purely mathematical choice; it is an algorithm that \textit{somewhat} gives reasonable search results (as perceived by the developers and managers in Google), and is implementable (reduces compute time/resources to levels that are affordable and, more pressingly, lead to a profitable business model). But Google search is used by billions of people every day, so this (hidden) mathematical choice will affect the way a vast proportion of the global population finds (or does not find) webpages. This ``choice of  algorithm'' then becomes a massive, unchecked, proprietary, global censorship object on the most important information tool of the modern age: the internet. There may be no malicious intent here, but that is beside the point. Harm may come to society if the way Google ranks webpages overly-suppresses, or overly-promotes, various content.

    \end{enumerate}

    \fi

\end{enumerate}

\newpage
\section{Exercises: Abstract Algebra}

\begin{enumerate}

\item \textbf{Topic: \textit{Group Theory}
 \\Allied disciplines: \textit{Computer Science}
\\Year level: 1-2
\\\hyperref[list:pillars]{Ethical pillars}: 1,9}

	The \textit{15-puzzle} consists of $15$ small square tiles, numbered $1$ to $15$, which are mounted in a $4 \times 4$ frame in such a way that each tile can slide vertically or horizontally into an adjacent square (if it is not already occupied by another tile), but the tiles cannot be lifted out of the tray. In the early 20th century, a cash prize was offered for a solution to manoeuvre the tiles from the first to the second of the configurations shown below. 
 
	\begin{center}
		\begin{tabular}{| c | c | c | c |}
			\hline
			$1$ & $2$ & $3$ & $4$ \\
			\hline
			$5$ & $6$ & $7$ & $8$ \\
			\hline
			$9$ & $10$ & $11$ & $12$ \\
			\hline
			$13$ & $14$ & $15$ & \\
			\hline
		\end{tabular}
	$~~~~~$
		\begin{tabular}{| c | c | c | c |}
			\hline
			$15$ & $14$ & $13$ & $12$ \\
			\hline
			$11$ & $10$ & $9$ & $8$ \\
			\hline
			$7$ & $6$ & $5$ & $4$ \\
			\hline
			$3$ & $2$ & $1$ & \\
			\hline
		\end{tabular}
	\end{center}

\begin{enumerate}
	    \item Give such a solution, or show that none exist.
     \item Outline the merits and drawbacks of offering a prize for such a puzzle. Are there any modern-day instances of scenarios where people use their mathematical knowledge to gain (financial, or other) advantage over others?
	\end{enumerate}

	\ifdefined\solution
    \solution 
	\textit{The idea of this question is to prompt the students into asking the question ``Should I be using my understanding of mathematics to trick people?'' This is a fairly harmless example, but as soon as we consider such puzzles as being sold for profit, then it becomes an issue of exploiting knowledge asymmetry. Of course, someone with sufficient mathematical training would realise that the puzzle is impossible, but not everyone has a mathematics degree; students may not be actively aware of this. A version of this was actually created and cunningly used by the American riddler Sam Loyd (1841-1911) with a prize of \$1000 offered to anyone who solved it, and of course no one could\footnote{You can read more about Sam Loyd on Wikipedia (\href{https://en.wikipedia.org/wiki/Sam_Loyd}{https://en.wikipedia.org/wiki/Sam\_Loyd}). If you like mathematical puzzles, we recommend the many books about Loyd's work and ideas, such as M. Gardiner's (1959) collection: \textit{Mathematical Puzzles of Sam Loyd} (Dover Publications).}. Mathematics is used extensively in the modern world to gain advantage over others, such as in stock market trading (in particular, high-frequency trading).}

 \begin{enumerate}
    \item There is no such solution. To get from the first configuration to the second requires 7 (i.e., and odd number of) transpositions: $(1\ 15)(2\ 14) \cdots (7\ 9)$. So it is an odd permutation. However, to keep the empty square in the bottom right corner requires an even number of vertical transpositions, and an even number of horizontal transposition; thus we must have an even permutation. Since permutations with differing signs are necessarily different, it means that we cannot get from the first configuration to the second.
    
    \item Some merits of offering a prize for this puzzle: the prize giver will never lose (and probably profit from sales of the puzzle), it will generate publicity, it will give people a toy to play with, it will make people think more about mathematics.    
    \\Some drawbacks of offering a prize for this puzzle: people will waste a lot of time trying to solve an unsolvable puzzle, the manufacturer will unfairly profit from the sale of such puzzles, such ``losses'' (of time and/or money) might disproportionately affect those worst off in society who are least likely to have received a high-level education and understand that such a puzzle is impossible.
    \\Mathematics is used extensively in the modern world to gain advantage over others, such as in stock market trading (in particular, high-frequency trading), in gambling (lotteries, scratch-and-win cards, sports betting, etc), with the (mis)use of statistics (in politics, in court cases, etc), and so on. This is a useful question to get students thinking about when their mathematical work \textit{adds value} to the world overall, and when their work simply \textit{extracts value} for their own personal gain.

\end{enumerate}

\fi

\newpage
     \item \textbf{Topic: \textit{Group Theory}
     \\Year level: 1-2
     \\\hyperref[list:pillars]{Ethical pillars}: 6,9}
     
     We say that a finite group $G$ is interesting if it satisfies $g^{2}=e$ for all $g \in G$.
     \begin{enumerate}
         \item  Show that if a finite  group is interesting, then it is also abelian.
        \item   Show that the group $C_{2}^{n}$ (=$C_{2} \times \ldots \times C_{2}$, $n$ times) is interesting for all  integers $n \geq 0$.
       \item   Show that if a finite group $G$ is interesting, then $G \cong C_{2}^{n}$ for some integer $n\geq 0$.
        \item   Is ``interesting'' a suitable name for such groups? Who should decide the name of a mathematical object, and how should it be done?
 \end{enumerate}

 \ifdefined\solution
    \solution 
	\textit{The purpose of this question is to illustrate that names can be misleading, or skew your interpretation of the purpose, functionality, or relevance of a mathematical object. They can change the way we view, value, and use mathematics, even before we have seen what it is and understand what it does (and doesn't) do.} 
 \begin{enumerate}
  \item Take $g,h\in G$. Then $g=g^{-1}$ and $h=h^{-1}$. Thus $[g,h]=ghg^{-1}h^{-1}=ghgh=(gh)^2=e$, so $G$ is abelian. \\(Note that $G$ need not be finite for this proof to work, but in the context of the whole question it seemed easier just to work with finite groups).
  \item If $n=0$ the $G=\{e\}$ so it is obviously interesting.
  \\If $n>0$ then consider an element $g=(a_{1},\ldots,  a_{n})\in C_{2}^{n}$, where each $a_{i}\in C_{2}$. Then
  \[
    g^{2}=(a_{1},\ldots,  a_{n})^{2}=(a_{1}^{2},\ldots,  a_{n}^{2})=(0, \ldots, 0)=e
  \]
  So $C_{2}^{n}$ is interesting.
  \item If $G$ is finite, consider  a minimal generating set for $G$ (that is, a generating set of smallest possible size); $\{g_{1}, \ldots, g_{n}\}$. Then any word in the $g_{i}$'s can be rearranged to look like $g_{1}^{\epsilon_{1}}\cdots g_{n}^{\epsilon_{n}}$, where $\epsilon_{i}\in \{0,1\}$. So there are at most $2^{n}$ such words.
  \\We claim these words are all distinct. So suppose they are not all distinct; take two such elements $g_{1}^{\epsilon_{1}}\cdots g_{n}^{\epsilon_{n}}=g_{1}^{\delta_{1}}\cdots g_{n}^{\delta_{n}}$. Let $k$ be the first place where $\epsilon_{k}\neq \delta_{k}$; wlog assume $\epsilon_{k}=1$ and $\delta_{k}=0$. Then we can express $g_{k}$ in terms of the remaining $g_{i}$'s, and so our set is not a minimal generating set for $G$; contradiction.
  \\So each $g \in G$ can be uniquely represented by some  $g_{1}^{\epsilon_{1}}\cdots g_{n}^{\epsilon_{n}}$, where $\epsilon_{i}\in \{0,1\}$. 
  \\(Note that this question stops here if we can assume the classification of finite abelian groups. If not, proceed with the rest of the proof).  
  \\Now construct a homomorphism
  \begin{align*}
      \phi: G &\to C_{2}^{n}
      \\ \phi(g_{1}^{\epsilon_{1}}\cdots g_{n}^{\epsilon_{n}})&:=(\epsilon, \ldots, \epsilon_{n})
  \end{align*}
Well-defined: for any $g \in G$, its representative $g=g_{1}^{\epsilon_{1}}\cdots g_{n}^{\epsilon_{n}}$ is unique, so $\phi(g)=(\epsilon, \ldots, \epsilon_{n})$ is well-defined.
\\Homomorphism: take $g=g_{1}^{\epsilon_{1}}\cdots g_{n}^{\epsilon_{n}}$, $h=g_{1}^{\delta_{1}}\cdots g_{n}^{\delta_{n}}$. Then
\[
\phi(gh)=\phi(g_{1}^{\epsilon_{1}}\cdots g_{n}^{\epsilon_{n}}g_{1}^{\delta_{1}}\cdots g_{n}^{\delta_{n}})=\phi(g_{1}^{\epsilon_{1}+\delta_{1}}\cdots g_{n}^{\epsilon_{n}+\delta_{n}})=(\epsilon_{1}+\delta_{1}, \ldots, \epsilon_{n}+\delta_{n})=\phi(g)\phi(h)
\]
Surjective: take $(\epsilon, \ldots, \epsilon_{n})$. Then $\phi(g_{1}^{\epsilon_{1}}\cdots g_{n}^{\epsilon_{n}})=(\epsilon, \ldots, \epsilon_{n})$.
\\Injective: this comes from the fact that $\phi$ is a surjective map between two groups of the same size.
    \item The name ``interesting'' suggests some sort of mathematical depth or intrigue to the object. However, here we have $C_{2}^{n}$, which can simply be viewed as a binary register of length $n$. There is no additional structure, and many would argue that such groups are not at all ``interesting''. But this then contrasts the often-made statement that ``Mathematics is simply abstract objects and concepts, so surely we give an object whatever name we like and it makes no difference to the mathematics''. Students may struggle to reconcile these two observations which in isolation seem reasonable (1. These groups are not interesting, and 2. The name of an object means nothing).
    
    ``Who should decide the name of a mathematical object, and how should it be done?'' opens up an interesting discussion point, and one that can actually be raised in a larger setting (eg: large group discussion). Who has the \textit{right} to name a mathematical object? What are their responsibilities in creating such a name? What obligations are there on others to follow that naming convention? This relates to power structures within the mathematical community; something students may not be very aware of as their training to date might have led them to believe that the mathematics community is self-administrated purely by logical argument.
    
 \end{enumerate}

\fi
    
\end{enumerate}

\newpage
\section{Exercises: Cryptography and Number Theory}

\begin{enumerate}
    \item 
    \textbf{Topic: \textit{Number theory / cryptography}
     \\Allied disciplines: \textit{Computer Science}
    \\Year level: 1-2
    \\\hyperref[list:pillars]{Ethical pillars}: 1,2,5,9}

\begin{enumerate}
    \item  
    Alice and Bob are using RSA public keys $(N, e_{1})$ and $(N, e_{2})$ with $e_{1}, e_{2}$ different and coprime. Their corresponding private keys are $(N, d_{1})$ and $(N, d_{2})$. They ask their colleagues to always send each of them the same message, encrypted to their respective keys. An eavesdropper Eve is monitoring their communications. Can Eve decipher the encrypted messages? Who would you tell, and what would you tell them?
    \item It turns out that Alice and Bob have an arrest warrant out for them, and Ana the police analyst is monitoring their communications. Now who would you tell, and what would you tell them?
    \item It turns that the arrest warrant is from a country with a long history of human rights violations. Now who would you tell, and what would you tell them?

    \item Amy and Bruce were friends. Bruce had given Amy an RSA public/private keypair $(M, f_{1})$ / $(M, h_{1})$, based on Bruce's public/private keypair $(M, f_{2})$ / $(M, h_{2})$, but with $f_{1}, f_{2}$ different; Bruce thought this was safe. Amy and Bruce had a bit of a falling out because Amy felt Bruce's line of work was unethical, but Amy is still able to monitor messages sent to Bruce. Can Amy decipher Bruce's incoming messages? Who would you tell, and what would you tell them?
  
    \item They fell out because Bruce works for an investigative journalism organisation wanting to protect its confidential sources. Now who would you tell, and what would you tell them?
\end{enumerate}
   
\ifdefined\solution
\solution \textit{This question illustrates how the context and framing of a mathematics question influences ``what problem should we be solving, and how?''. Here, the same (mathematical) setup is posed with various different contexts, and at the end of each sub-question the student is asked ``ok, now that you have done some mathematics, what will you actually DO with it?''. The answer to this gets more and more ambiguous and ethically challenging with each sub-question, and students should be warned off thinking that there is a ``right'', or ``obvious'' conclusion in each case.}

\begin{enumerate}
    \item Let $m$ be a message that is encrypted to both RSA keys, transmitted, and intercepted by Eve. Then we have 
    \begin{align}
        m^{e_1} &\equiv  c_1 \ \  (mod\ N)
        \\m^{e_2} &\equiv  c_2 \ \ (mod\ N)
    \end{align}
    where $c_{1}, c_{2}$ are the encrypted messages transmitted to Alice and Bob respectively.

    As $e_{1},e_{2}$ are coprime and public, Eve can use Euclid's algorithm to compute integers $k,k$ such that $he_{1}+ke_{2}=1$.
    
    As $N$ is public, and Eve has intercepted $c_{1}, c_{2}$, Eve can now compute $c_{1}^{h}$ mod $N$, and $c_{2}^{k}$ mod $N$, by simply raising powers. Note that if either $h$ or $k$ is negative (say $h$), Eve can first compute a multiplicative inverse of $c_{1}$ mod $N$ using Eulid's algoritm again, as
    \[
    c_{1}x \equiv 1\ (mod\ N) \iff c_{1}x +tN=1 
    \]
    so as it is extremely likely that $c_{1},N$ are coprime ($N$ has only two prime factors, which are both very large), Eve can use Euclid to find $x,t$ above, and $x$ is the desired multiplicative inverse. [Note: if $c_{1},N$ share a factor, the Euclid's algorithm will yield such a factor, and thus enable Eve to factorise $N$, giving total information about each RSA key and enabling Eve to decipher everything anyway.]
    Now Eve can compute $c_{1}^{h}c_{2}^{k}$ mod $N$ and use this to recover the original message $m$, as
    \[
    m=m^{1}=m^{he_{1}+ke_{2}}\equiv c_{1}^{h}c_{2}^{k} \ (mod \ N)
    \]
    (observing that $m<N$, and so $m=c_{1}^{h}c_{2}^{k}$ mod $N$).

    The semi-obvious thing students might suggest to do at this point is to tell Alice and Bob that Eve can potentially decipher messages sent to them, and that at least one of them should change their RSA keys to use a different modulus. However, there is an implicit assumption that ``all evesdroppers are bad'', and an astute student may ask more questions about the situation before deciding who to tell and what to tell them. This leads to the next question.

    \item In this more detailed scenario, the mathematics of the system is exactly the same. Ana would be able to decipher the messages sent to Alice and Bob. The semi-obvious thing students might suggest is to tell Ana \textit{how} to decipher messages sent to Alice and Bob. However, there is an implicit assumption that ``all people who have a warrant for their arrest are bad'' and ``all police analysts are good''. An astute student may ask more questions about the situation before deciding who to tell and what to tell them. This leads to the next question.

    \item In this even more detailed scenario, the mathematics of the system is again exactly the same, and moreover the actors are the same. But now we have more information about some of the actors. This may lead some students to pause and reconsider their choice of who to tell, and what to tell them. It may be that the country does have a bad human rights record, but that Alice and Bob are charged with murder with a very strong case against each of them. It may be that Alice and Bob have been accused of politically-motivated offences and have not even been charged with any crime. The point of this question is to get students to start and consider as many scenarios as possible, and how each scenario changes the optic of the approach.

    \item We first show that, by using $M, f_{1},h_{1}$, Amy is able to factorise $M=pq$. Note that, for now, Amy does not know $p$ or $q$. 

    Amy knows $f_{1},h_{1}$, so can write $f_{1}h_{1}-1=2^{s}k$ for some odd $k$. We define $ord_{p}(x)$ to be the order of an element $x$ in the multiplicative group $\mathbb{Z}_{p}^{\times}$. Now we define the set
    \[
    H:=\{x \in \mathbb{Z}_{M}^{\times} \ | \ ord_{p}(x^{r}) \neq ord_{q}(x^{r})\}
    \]

    If Amy happens to have some $x \in H$, then it becomes easy to factorise $M$. This is because, if $x \in H$, then we can set $y:=x^{r}$ and then observe that, in $\mathbb{Z}_{M}^{\times}$, we have:
    \[
    y^{2^{s}}=x^{2^{s}r}=x^{de-1}=1\ (\textnormal{all in } \mathbb{Z}_{M}^{\times})
    \]
    as $de\equiv 1\ (mod\ \varphi(M))$ and so $x^{de}\equiv x^{1}\ (mod\ M)$.
    
    So $ord_{p}(y)$ and $ord_{q}(y)$ must both be
powers of 2; without loss of generality we can suppose that $ord_{p}(y) = 2^{t} < ord_{q}(y)$.

    Thus  $y^{2^{t}} \equiv  1\ (mod\ p)$, but $y^{2^{t}} \not\equiv 1\ (mod\ q)$, and hence
    \[
    gcd(y^{2^{t}}-1,M) = p
    \]
    Hence, with $y^{2^{t}}-1$, Amy can easily use Euclid's algorithm to find $p$ and so factorise $M$. Even without knowing $t$ exactly, Amy can just try $t=1,2, \ldots, s$ in order until Euclid's algorithm yields $p$; this is not a long sequence to test because $s \leq \log_{2}(x)$.

    However, how does Amy find an $x \in H$ to start this process? Well, it turns out that the following is true (by a standard theorem in this subject):
    \[
    |H|>\frac{|\mathbb{Z}_{M}^{\times}|}{2}
    \]
    Thus, a random choice of element from $\mathbb{Z}_{M}^{\times}$ has at least a $\frac{1}{2}$ chance of being in $H$. So repeatedly choosing elements from $\mathbb{Z}_{M}^{\times}$ will quickly yield an element in $H$, as running through $k$ choices will find $p$ with probability $1-\frac{1}{2^{k}}$.

    Having found $p$, and hence $q=\frac{M}{p}$, Amy can now compute $\varphi(M)-(p-1)(q-1)$, and hence find $h_{2}$ such that $f_{2}h_{2}\equiv 1\ (mod\ \varphi(M))$ and thus Bruce's private key $(M, h_{2})$.

    This question poses an interesting sub-problem: what are the types of people (like Bruce) who possess and RSA keys, but do NOT understand the technicalities of how that are used? Not everyone who uses cryptography knows all the intricacies of how it works. Hopefully this leads on well to the next question.

   \item This is hard. In this more detailed scenario, the mathematics of the system is  exactly the same, and moreover the actors are the same. And now we have more information about some of the actors. But the system now presents a distinct ethical challenge, and one where society at large has no broad and well-established stance.  There is much debate about whether, and under what circumstances, journalists should protect their sources. Moreover, this debate changes as you go to different parts of the world, or different times in history. Investigative journalism, and the informants and whistleblowers that they often deal with, present distinct ethical challenges. And yet these are exactly some of the people who will make frequent use of cryptographic tools.

   This question is here to prompt the student to consider, possibly for the first time ever, such a scenario, and to see if they can acknowledge the complexities of the situation. There is certainly no right or wrong answer here, and it is an area that is fiercely debated in society, but yet one that students may find themselves working in later in their career.

\end{enumerate}

\fi

\newpage
\item \textbf{Topic: \textit{Number theory / cryptography}
 \\Allied disciplines: \textit{Computer Science}
    \\Year level: 1-2
    \\\hyperref[list:pillars]{Ethical pillars}: 1,6,7,9,10}
    
    \begin{enumerate}
        \item What are some of the consequences - mathematical, and social - of the existence of a fast (ie: polynomial time) factorisation algorithm?
        \item What would you do if you found a fast factorisation algorithm?
    \end{enumerate}

\ifdefined\solution
\solution \textit{This question looks at the notion of responsible disclosure, which is well-understood in computer science (when bugs or exploits are found in code), but not as frequently discussed in mathematics. There are certain mathematical tools, algorithms, truths, etc which, if discovered or developed, would have a monumental impact on society; fast factorisation is one of the most significant of these. \\ Students should try and understand the impact on society, and the impact on them both professionally and personally, of making such a revelation, especially if it is done irresponsibly. The notion of ``I can publish whatever mathematical truth I want, whenever I want, and however I want'' doesn't quite work here. Just because something is proven true, doesn't necessarily mean it should be released. \\ Students should also explore the options available for releasing such works. It is not a binary choice between ``publish'' and ``don't publish''; there are other intermediate options.}

\begin{enumerate}
    \item There are many mathematical consequences of having a fast factorisation algorithm. One of the more significant ones is that it would lead to the breaking of RSA encryption, as that whole system relies on the (believed) difficulty of factorising integers.
    \\The social consequences of breaking RSA cannot be understated. Many internet security systems use RSA, including HTTPS (webpage encryption), PGP (encrypted emails), and digital signatures for verifying software updates. In short: most internet security would fall apart overnight if a fast factorisation algorithm were to become public knowledge. As such, anything done over the internet that involved ``trust'' (eg: commerce, banking, secure messaging, software updates, etc) would collapse, and the internet would be reduced to one giant public-readable notice board. There would be significant disruption to the buying and selling of goods, to banking, and to the global supply chain. In short: given our current reliance on the internet for global supply systems, it would quickly become very difficult to do things like access money, buy food, etc.

\item What to do is a \textbf{very hard} question. Responsible disclosure is an active field of research, and even the research communities in adjacent fields where such things arise more frequently (computer exploit discovery, infectious disease discovery, etc) do not have an exact, agreed methodology for what to do. It depends on the situation.

Before deciding \textit{what} to do, students should first try and write down a list of what the \textit{can} do, and then explore the pros and cons of each option. There are more options than just ``publish'' and ``don't publish''. The following is a (non-exhaustive) list of other considerations:
\begin{enumerate}
    \item Can a 0-knowledge proof be released (eg: solve a few factorisation challenges). This verifies to the community that the algorithm exists, without releasing it, and motivates the IT industry to quickly migrate away from RSA and to other cryptosystems.
    \item Can the algorithm be placed with a ``custodian'' (eg: government)?
    \item What are the risks and threats to the researcher if the whole world knows they have a (secret) algorithm for fast factorisation? Will their safety, or the safety of their family and friends, be at risk from nefarious actors (criminals, hostile governments, their own government, etc)?
\item How effectively can the researcher control this secret? Is it being kept in a standard office with minimal security (standard locks, a ground-floor window, etc)?
\item Who else already knows about this (were there any collaborators), and do they pose a risk that the researcher needs to take into account (will their collaborator try and sell the algorithm to criminals? Does this change what the researcher does, and when they do it?)
    
\end{enumerate}

\end{enumerate}

\fi

\end{enumerate}
\newpage
\section{Exercises: Differential Equations}\label{diff-Qs}\begin{enumerate}
\item \textbf{Topic: \textit{Differential Equations}
 \\Allied disciplines: \textit{Economics, Engineering, Statistics, Natural Sciences}
\\Year level: 1
    \\\hyperref[list:pillars]{Ethical pillars}: 1,3,6,7}
    
A detective arrives at the scene of a crime at 5:00pm, finds a warm cup of tea, and measures its temperature at 40°C. By 5:30pm the tea's temperature has dropped to 30°C. 
\begin{enumerate}
    \item The police approach you with this data and ask you when the tea was likely made. Briefly discuss any questions that you still need to ask the police officers and their potential ethical relevance. What are potential barriers of communication?
\item The police are unable to provide you with more information, but ask you to give an estimate based on idealised conditions and a constant room temperature of 20°C. Giving all mathematical details and assumptions, use Newton's law of cooling to estimate when the tea was likely made. What error margins might your computation  have, and what might the consequences of these be?
\end{enumerate}

\ifdefined\solution
    \solution
\textit{The purpose of this question is to teach the students to think about the assumptions that go into typical physical models. It also exposes the student to the fact that the information about these assumptions (values for boundary conditions, etc) do not come out of a vacuum, but can present together with some real communicative challenges. In this question, students also learn that making ``obvious'' mathematical assumptions can have ethical consequences, and it teaches them to be humble when making any assumptions.
\\ Additional commentary detailing the politics of this and other questions can be found in M\"uller and Chiodo\footnote{ D. M\"uller, M. Chiodo (2023). \textit{Mathematical Artifacts have Politics: From Examples to Embedded Ethics.} \href{https://arxiv.org/abs/2308.04871}{https://arxiv.org/abs/2308.04871}}.}

\begin{enumerate}
\item The student should realise that the surrounding temperature of the room is unknown and that a varying room temperature will lead to different times for when the tea was made. They should discuss the ethical relevance of this in the context of pursuing a wrong subject. They should also consider if the cup of tea was moved outside (i.e. the surrounding temperature was not constant over time), heated up manually or drank from (i.e. the volume of tea was not constant over time), and that this would lead to a variation in temperature and different brewing times. The student should understand that the police officer likely has very little mathematical training. Hence, they should be able to break down their communication and explain the relevance of these questions in everyday language. 

The student should understand that in certain situations they cannot simply make assumptions without further input. 

\item Let $T_e(t)$ be the unknown temperature of the surrounding room. 
    We assume $T_e$ is the typical ideal room temperature and constant at 20°C. Let $T(t)$ be the temperature of the tea, where we set $5pm \leftrightarrow  0$ and $5.30pm \leftrightarrow 1$. We assume that the cup was not heated up nor moved since it was first put on the table. According to Newton's law of cooling
$$\frac{dT}{dt} = k\left(T_e - T \right).$$
Then $T$ satisfies the conditions $T(0) = 40$ and $T(1) = 30$, leading us to solve the differential equation
$$\frac{dT}{dt} = k\left(20 - T \right)$$
$$\Rightarrow \frac{dT}{dt} + kT= 20k$$
By setting $u = e^{kt}$, we obtain
$$ u\frac{dT}{dt} + T \frac{du}{dt} = 20uk$$
$$ \Rightarrow \frac{d}{dt}\left( uT \right) = 20uk$$
$$  \Rightarrow e^{kt} T = e^{kt} 20 + c$$
$$ \Rightarrow T(t) = 20 + c e^{-kt}$$
We now use the given temperature constraints to obtain the constants $c$ and $k$:
$$T(0) = 20 + c = 40 \Rightarrow c = 20$$
$$T(1) = 20 + 20 e^{-k} = 30 \Rightarrow k = -\ln\left(\frac{1}{2}\right).$$
Thus, the temperature of the tea at time $t$ is given by
$$T(t) = 20 + 20 e^{\ln\left(\frac{1}{2}\right)t}.$$
To find the time of brew, we solve
$$100 = 20 + 20 e^{\ln\left(\frac{1}{2}\right)t}$$
$$\Rightarrow t = \frac{\ln(4)}{\ln(1/2)} = - 2 $$
Thus the tea was made at about 4pm.

However, there are significant error margins. The ambient temperature could well have been below 15°C (if it were winter and the crime scene had not been occupied for a while, or alternatively was outdoors), or it could have been over 30°C if it were a hot day in a building with no airflow, or a room that was excessively heated. The student could, at the very least, re-run the computations with these ``extremal'' ambient temperatures. Also, at this point the students  might try and think about any correlation between 1) the room/area not having idealised conditions and 2) the room/area being the scene of a crime.

The tea may have been made with water that was below 100°C ; some people are impatient and take the kettle before it finishes boiling, or the scene may have been at altitude where the boiling point of water is lower. Again, the student could re-run the computation with lower starting temperatures (95°C, 90°C).

The tea may have also had something cold put into it, such as milk from a fridge, a large metal spoon, etc, which would have decreased the initial/early temperature even further, by perhaps another 5°C. And, of course, some tea could have been drunk early on, changing the rate at which the temperature drops (a smaller cup of tea will experience a faster temperature drop).

The student might see some, or all, of these implicit assumptions/factors, or may well come up with more of their own. They should hopefully try and do at least one ``error'' calculation, whereby the stack all the errors in one direction (e.g., all the small errors that might lead to a underestimate of the time the tea was made), and see how far off their original computation they were.

This computation might be used as evidence in a court, which could have a significant impact on the life and liberty of the accused. So the consequences on that person/people would be significant. However, in a court setting, the evidence may need to be presented in ``worst case scenario'' format, giving extremal upper and lower bounds for when the tea was made with all reasonable assumptions stress-tested (while it might be unreasonable to factor in aliens coming and cooling down the tea, it would be perfectly reasonable to assume that the kettle had been switched off before the water had fully boiled).

In addition (and perhaps something the student may not have picked up on) is that fact that the computation relies on several assumptions which were imposed on them by the police: the assumption of ``idealised conditions and a constant room temperature of 20°C''. The police might be seeking a prosecution and conviction, for which such ``irrefutable mathematical evidence'' may be of great assistance. The student, acting as a mathematician, may inadvertently be part of a miscarriage of justice if they were to not make it clear to any court using their evidence that such assumptions were requested by the police investigating the crime.

\end{enumerate}

\fi

\newpage
\item \textbf{Topic: \textit{Differential Equations}
 \\Allied disciplines: \textit{Economics, Engineering, Statistics, Natural Sciences}
\\Year level: 0-1
    \\\hyperref[list:pillars]{Ethical pillars}: 1,9}
    
    Consider two armies $R$ (red) and $B$ (blue). Let $m_R(t)$ and $m_B(t)$ be the number of soldiers of army $R$ and $B$ respectively. 
\begin{enumerate}
    \item Assume that the losses of each army are proportional to the strength of the other army with proportionality constants $a_R$ and $a_B$. Set up a system of differential equations that describe the strength of each army at time $t$.
    \item Derive Lanchester's square law
$$a_B\left(m_B^2(t) - m_B^2(0) \right) = a_R\left(m_R^2(t)-m_R^2(0)\right).$$
\item How would you describe the proportionality constants $a_R$ and $a_B$ in everyday language? Would you be comfortable working on such problems in your career?
\end{enumerate}

\ifdefined\solution
\solution
\textit{This question is designed to show that abstract mathematics can be very ``real'', and can have very substantial life-and-death consequences, once you put it into everyday language. Here we have mathematics for military engagement; this is very real, and very serious.}
\begin{enumerate}
    \item This question asks you to set up Lanchester's model $$ \frac{dm_B}{dt} = - a_R m_R$$
    $$ \frac{dm_R}{dt} = -a_B m_B $$
    with initial conditions $m_B(0) = m_{B_0}$ and $m_R(0) = m_{R_0}.$
    \item Consider the quotient and separate variables
    $$ \frac{dm_B}{dm_R} = \frac{a_R m_R}{a_B m_B}$$
    $$\Rightarrow a_B m_B dm_B  = a_R m_R dm_R $$
    $$ \Rightarrow a_B \int\limits m_B dm_B = a_R \int\limits m_R dm_R$$
    
   Integrating gives
   \[
   a_B \frac{m_B^2(t)}{2}+c= a_R \frac{m_R^2(t)}{2}+d
   \]
   Subtracting the substituted term $t=0$ from both sides,  gives the desired result:
$$a_B\left(m_B^2(t) - m_B^2(0) \right) = a_R\left(m_R^2(t)-m_R^2(0)\right)$$

    \item The coefficients $a_B$ and $a_R$ can be interpreted as the killing rate of the armies. 
    \\This piece of mathematics might influence whether a general sends their army into battle or not, by analysing whether ``More of their troops will die than our troops'' ($a_B$ and $a_R$ are ``death rates'' when viewed by the other side). Students may have never considered the fact that even a 5-line piece of calculus can have such a significant impact on warfare (and thus on human life). There is no right or wrong answer to the last part of this question; it only serves to \textit{initiate} these sorts of considerations among students, and should be discussed with open-minded, thoughtful, and non-judgemental discourse.
\end{enumerate}
\fi

\newpage
\item \textbf{Topic: \textit{Differential Equations}
 \\Allied disciplines: \textit{Economics, Engineering, Statistics, Natural Sciences}
\\Year level: 1-2
    \\\hyperref[list:pillars]{Ethical pillars}: 3,6}

The evolution of an infectious disease in a population can be modelled by
    \begin{align*}
        \dod{U}{t} &= U \left(1 - (U+I)\right) - \beta UI \\
        \dod{I}{t} &= I \left(1 - (U+I)\right) + \beta UI - \delta I,
    \end{align*}
    where $U$ is the uninfected population and $I$ is the infected population. 
    
    \begin{enumerate}
    \item Explain this model to a biologist. What are $\beta$ and $\delta$? Suppose we model herpes simplex or influenza using these equations: do you think $\beta$ and $\delta$ are small or large for these diseases?

    \item For $\beta=3/4$, determine the location and stability of the critical points of the above system in the three cases (i) $\delta=1/5$, (ii) $\delta=2/5$, (iii) $\delta=3/5$. Thus determine the long-term outcome for the population in each case.
    
    \item Consider a disease with $\beta=3/4$ and $\delta=3/5$. A drug is discovered that can reduce $\delta$ to $2/5$. Would you distribute the drug to everybody? If you had contracted the disease, and were offered the drug, would you take it?
    \end{enumerate}
    \ifdefined\solution

    \solution\textit{The purpose of this question is to identify what the term} worse \textit{means, and in particular, worse} for whom? \textit{A government or hospital may prefer $\delta=\frac{3}{5}$ as it leads to an overall lower death toll; an individual being infected by the disease would almost certainly prefer $\delta=\frac{2}{5}$. Utilitarianism is great, unless you're on the losing end of it! Mathematicians are good at solving optimisation problems, but are seldom made to question} what \textit{to optimise for.}
    
    \begin{enumerate}
    \item In this simplistic model, the parameters $\beta$ and $\delta$ represent the rates of infection and mortality, respectively.
    
    A disease with a very high mortality rate kills anybody who catches it very quickly, but those people would not be able to infect very many people. A disease with a low mortality rate, such as herpes simplex, can spread to a large population\footnote{See also: Wikipedia, \textit{Epidemiology of herpes simplex}: \href{https://en.wikipedia.org/wiki/Epidemiology_of_herpes_simplex}{https://en.wikipedia.org/wiki/Epidemiology\_of\_herpes\_simplex}.}.

    \item Substituting $\beta = \frac{3}{4}$, we have
	\begin{align*}
	\dot{U} &= U(1 - U - \tfrac{7}{4} I) \\
	\dot{I} &= I((1 - \delta) - \tfrac{1}{4} U - I)
	\end{align*}
	Thus we find the critical points
	\begin{align*}
	1) ~~ (U,I) &= (0,0) ~, \\
	2) ~~ (U,I) &= (1,0) ~, \\
	3) ~~ (U,I) &= (0,1-\delta) ~, \\
	4) ~~ (U,I) &= (\tfrac{4}{9} (7\delta - 3), \tfrac{4}{9} (3 - 4\delta)) ~,
	\end{align*}
	and the Jacobian of the system is
	\begin{equation*}
	J(U,I) = \begin{pmatrix}
	1 - 2U - \frac{7}{4} I & -\frac{7}{4} U \\
	-\frac{1}{4} I & 1-\delta - \frac{1}{4} U - 2 I
	\end{pmatrix} ~.
	\end{equation*}
	At the first two critical points we find
	\begin{equation*}
	J(0,0) = \begin{pmatrix}
	1 & 0 \\
	0 & 1-\delta
	\end{pmatrix} ~, ~~ \text{and} ~
	J(1,0) = \begin{pmatrix}
	-1 & -\frac{7}{4} \\
	0 & \frac{3}{4} -\delta
	\end{pmatrix} ~.
	\end{equation*}
	Thus by inspection, in each of the three given cases, $(U,I) = (0,0)$ is an unstable node and $(U,I) = (1,0)$ is a saddle (the Jacobian here has negative determinant).
	
	Now considering each case individually:
	\begin{enumerate}
	    \item $\delta = \frac{1}{5}$: Point $3)$ is at $(U,I) = (0,\frac{4}{5})$ with Jacobian
	\begin{equation*}
	J(0,\tfrac{4}{5}) = \begin{pmatrix}
	- \frac{2}{5} & 0 \\
	-\frac{1}{5} & -\frac{4}{5}
	\end{pmatrix} ~,
	\end{equation*}
	This has positive determinant and negative trace, hence it is stable.
	
	Point $4)$ is not in $U \geqslant 0, I \geqslant 0$.
	
	\item $\delta = \frac{2}{5}$: Point $3)$ is at $(U,I) = (0,\frac{3}{5})$ with Jacobian
	\begin{equation*}
    	J(0,\tfrac{3}{5}) = \begin{pmatrix}
    	- \frac{1}{20} & 0 \\
    	-\frac{1}{4} & - \frac{3}{5}
    	\end{pmatrix} ~.
	\end{equation*}
	Again this has positive determinant and negative trace, hence it is stable.
	
	Point $4)$ is still not in $U \geqslant 0, I \geqslant 0$.
	
	\item $\delta = \frac{3}{5}$: Point $3)$ is at $(U,I) = (0,\frac{2}{5})$ with Jacobian
	\begin{equation*}
	J(0,\tfrac{2}{5}) = \begin{pmatrix}
	\frac{6}{20} & 0 \\
	-\frac{2}{20} & -\frac{2}{5}
	\end{pmatrix} ~.
	\end{equation*}
	This is now a saddle (it has negative determinant).
	
	Point $4)$ is at $(U,I) = (\tfrac{8}{15},\tfrac{4}{15})$ with Jacobian
	\begin{equation*}
    	J(\tfrac{8}{15},\tfrac{4}{15}) = \begin{pmatrix}
    	- \tfrac{8}{15} & - \tfrac{14}{15} \\
    	- \tfrac{1}{15} & - \frac{4}{15}
    	\end{pmatrix} ~,
	\end{equation*}
	This has positive determinant and negative trace, hence it is stable.
    \end{enumerate}

	Thus our long-term outcomes for the population are
	\begin{enumerate}
		\item $\delta = \frac{1}{5}$: $(U, I) \rightarrow (0, \tfrac{4}{5}) ~,$ total population $= \tfrac{4}{5}$
		\item $\delta = \frac{2}{5}$: $(U, I) \rightarrow (0, \tfrac{3}{5}) ~,$ total population $= \tfrac{3}{5}$
		\item $\delta = \frac{3}{5}$: $(U, I) \rightarrow (\tfrac{8}{15}, \tfrac{4}{15}) ~,$ total population $= \tfrac{4}{5}$
	\end{enumerate}
	The intermediate value of $\delta$ is the most harmful for the total long-time stable population; this happens because it strikes a balance between the infectiousness and the mortality of the disease -- a low mortality allows for a large completely infected population, a high mortality inhibits the spread of the disease. This is perhaps somewhat counter-intuitive, as one might initially assume that a lower mortality rate would necessarily lead to lower overall deaths. But a lower mortality rate can lead to an overall larger number of people contracting the disease. But it is subtle; if the mortality rate is very low, we see that overall deaths are also lowered.

    \item     Reducing $\delta$ would make it more likely for an infected person to survive, but it would also increase the number of people that get infected, and increase the total number of deaths! A politician might argue against distributing the drug, ``for the greater good'' (i.e., to reduce the number of overall deaths). But what if your friend or family member contacted the disease; would you refuse them the drug? Is it ethical to refuse someone a drug that could save their life, in order to minimise overall lives lost? If you contracted the disease and were offered the drug, would you refuse the drug and increase the chance of your own death ``for the greater good''?
    \\This problem again highlights the need to be careful when interpreting a model: the model itself cannot advise whether it is more ethical to increase the mortality rate for individuals or to allow a less fatal disease to spread more widely. The decision requires an additional value judgement. 

\end{enumerate}
\fi 

\newpage
\item \textbf{Topic: \textit{Differential Equations} \\ Year level: 1 
 \\Allied disciplines: \textit{Economics, Engineering, Statistics, Natural Sciences}
\\\hyperref[list:pillars]{Ethical pillars}: 1,5,10} 

``\textit{Dual Use Research} is defined as research conducted for legitimate purposes that generates knowledge, information, technologies, and/or products that could be utilised for both benevolent and harmful purposes.'' \footnote{Source: Boston University Office of Research (n.d.). \textit{Dual Use Research of Concern}. \href{https://www.bu.edu/research/ethics-compliance/safety/biological-safety/ibc/dual-use-research-of-concern/}{https://www.bu.edu/research/ethics-compliance/safety/biological-safety/ibc/dual-use-research-of-concern/}}
\begin{enumerate}
    \item How does dual use come up in the study of differential equations?
    \item Do you know a differential equation that can be applied in a benevolent and harmful way? (Hint: consider differential equations and their applications from the lectures. Can you apply some of them somewhere else?)
\end{enumerate}
\ifdefined\solution
\solution
\textit{This question asks students to look for beneficial and harmful scenarios for their mathematics. Such a task is typically not explored in most undergraduate classes, but it is a valuable lesson on ethics: while most people are generally good, some people may actively seek to do harm. Students should be aware that ``just because it's mathematics'' doesn't necessarily mean it is ``morally good''. The question may be well suited for in-class discussions or recitation sections. \\Additional commentary detailing the politics of this and other questions can be found in M\"uller and Chiodo\footnote{D. M\"uller, M. Chiodo (2023). \textit{Mathematical Artifacts have Politics: From Examples to Embedded Ethics.} \href{https://arxiv.org/abs/2308.04871}{https://arxiv.org/abs/2308.04871}.}.}
\begin{enumerate}
    \item Students should realise that the way we teach differential equation is often through teaching about abstract classes of differential equations (e.g. ``this is how linear differential equations of the first order look like, and they can be solved like that in their general form'') that can then be applied to many situations.
    \item Useful examples depend on the course where this question is used. For first-year differential equations, it can make sense to consider the mixing problem or diffusion problems in beneficial situations (How long do I have to wait until polluted water is clean again?) and harmful situations (How much pollutant is needed to poison a water supply?). For later courses (e.g. on PDEs), the Neutron Diffusion Equation is a good example (used in the design nuclear reactors and atomic bombs).

    A very broad example of dual-use differential equations are the Maxwell Equations in electrodynamics:
    
\begin{align*}
    \nabla \cdot \mathbf{E} &= \frac{\rho}{\varepsilon_0} \\
    \nabla \cdot \mathbf{B} &= 0 \\
    \nabla \times \mathbf{E} &= -\frac{\partial \mathbf{B}}{\partial t} \\
    \nabla \times \mathbf{B} &= \mu_0 \mathbf{J} + \mu_0 \varepsilon_0 \frac{\partial \mathbf{E}}{\partial t}
\end{align*}

They formed the foundation of much of the new technology of the 20th century (and beyond). They enabled the transmission of electricity to homes and institutions around the world for good and bad. They enabled TV and radio to broadcast the sound and images of people, including the propaganda of political dictators. They increased the availability of heat and light into homes, and thus encouraged the consumption of fossil fuels in the generation of electricity.

For a more specific example, one could consider the \textit{Neutron Diffusion Equation}, which governs how neutrons are produced and diffused in a uniform material:
    \[
    \frac{\partial \Phi}{\partial t}= \frac{\sigma-1}{\tau}\Phi + \frac{\lambda^{2} }{3 \tau}\nabla^{2} \Phi
    \]
    where $\Phi(\mathbf{x},t)$ is the (free) neutron density, $\sigma$ is the average number of neutrons released in a fission event, $\tau$ is the average time between fissions, and $\lambda$ is the average distance a neutron travels before being absorbed by a nucleus.

This is worked through in full later in this section. It models nuclear fission sufficiently well to enable technologies such as nuclear power generation, but also devices such as atomic weapons.

\end{enumerate}

\newpage
\item\label{chem-spill} \textbf{Topic: \textit{Differential Equations}
 \\Allied disciplines: \textit{Economics, Engineering, Natural Sciences}
\\Year level: 1
    \\\hyperref[list:pillars]{Ethical pillars}: 5,8,10}
    
A chemical accident happened near the small village [city name]. The city's local water reservoir has volume $V$. The inflow and outflow of the reservoir is given by the flow rate $r$. Let $x(t)$ be the amount of [chemical] in the reservoir at time $t$. Assume that the reservoir was clean at the beginning, i.e. $x(0) = 0$. Let $C_e(t)$ be the concentration of the [chemical] flowing into the reservoir.
\begin{enumerate}
    \item Set up and solve a differential equation describing the concentration of the reservoir.
    \item How can you use your solution to model repeated pollution (e.g. criminals dumping [chemical] near the reservoir every weekend)?
    \item What are some relevant questions you can ask about the concentration of [chemical] in the reservoir?
    \item Suppose the polluter is caught and after some cleaning the incoming water is clean. How can you use your model to analyse when the water in the reservoir will be safe again? How sure are you of your answer, and how much does it matter?
\end{enumerate}
\ifdefined\solution
\solution
\textit{This question is designed to show students that very simple mathematics can be used to model local environmental disasters, which can often be an example of how it may be used unsustainability. It teaches students to find good questions instead of merely answering someone else's questions.}
\begin{enumerate}
    \item Consider ``rate of change = rate of chemical inflow - rate of chemical outflow''.
    $$ \frac{dx}{dt} = rC_e - r \frac{x}{V}$$
    $$ \frac{dx}{dt} + r\frac{x}{V} = r C_e $$
    Let $C = \frac{x}{V}$ be the concentration of [chemical] in the reservoir. Plugging in we obtain a differential equation purely involving concentrations.
    $$V\frac{dC}{dt} + rC = rC_e$$
    $$\Rightarrow \frac{dC}{dt} + \frac{r}{V}C = \frac{r}{V}C_e$$
    To solve this,  multiply with the integrating factor $e^{\frac{r}{V}t}$ to obtain
    $$ e^{\frac{r}{V}t} \frac{dC}{dt} + \frac{r}{V} C e^{\frac{r}{V}t} = \frac{r}{V}C_e e^{\frac{r}{V}t} $$
    $$ \Rightarrow \frac{d}{dt}\left(C e^{\frac{r}{V}t} \right) = \frac{r}{V}C_e e^{\frac{r}{V}t} $$
    $$ \Rightarrow  C e^{\frac{r}{V}t} = \int \frac{r}{V}C_e e^{\frac{r}{V}t} dt$$
    $$ \Rightarrow C  =e^{-\frac{r}{V}t}\int \frac{r}{V}C_e e^{\frac{r}{V}t} dt$$
    $$ \Rightarrow C  = C_{e}+ e^{-\frac{r}{V}t} \cdot const$$
    Plug in the initial condition $x(0) = 0 \Rightarrow C(0) = 0$ to get
    $$ \Rightarrow 0  = C_{e}+ e^{0} \cdot const$$
    $$ \Rightarrow const  = -C_{e}$$
    $$ \Rightarrow C  = C_{e}(1- e^{-\frac{r}{V}t})$$

    \item This can be modelled using a sinusoidal input concentration $C_e$. Alternatively, you can also use a linear combination of shifted/stretched indicator functions $\pi(t) = \mathbf{1}_{t\in [-1,1]}$ if you want to be able to model different concentrations for the different dumps of the chemical:
    $$\sum\limits_{i=1}^n c_i \pi\left(\frac{t-t_i}{T_i}\right).$$ 
    The differential equation is linear, so we know that the solution will itself be a linear combination of the solutions corresponding to the individual dumps of the chemical.
    \item Possible questions include
    \begin{itemize}
        \item Will the pollution of the reservoir ever reach a dangerous level in the reservoir?
        \item What is deemed a ``safe'' level of [chemical] in the reservoir?
        \item How closely does the concentration of the reservoir follow the inflow of pollutant chemicals?
        \item Will the reservoir reach an equilibrium concentration of [chemical]?
    \end{itemize}
    \item Use the differential equation from (a) to model the behaviour until the time when the input stream is clear. This gives the concentration of [chemical] in the reservoir at that time. Then set $C_e = 0$ and use this as the new initial condition to solve the equations again, but now seeking to find $t$ such that the level of [chemical] in the reservoir is deemed safe (a number that the student should have mentioned in the previous part of the question). Students need to understand that not all relevant numbers will be given to them at the start of a problem, and that often they will need to specifically ask for such numbers when carrying out their mathematical work.

    For ``How sure are you?'', students will need to explore what the known unknowns are:  errors in the measurement apparatus, non-uniform mixing, samples taken in a very clean/dirty part of the stream or reservoir, etc. They will also need to discuss any unknown unknowns: other sources of pollutant (other streams leading into the reservoir), samples being tampered with accidentally or deliberately, etc.

    For ``How much does it matter?'', students should identify that we are dealing with poison in drinking water, so it matters immensely! They should understand that this is an \textit{estimate}, which is there to \textit{forecast} when the water \textit{might} be safe to drink. The only way to actually know if it is safe is to thoroughly test it before including it in the water supply; the forecast is there to help government etc plan things such as how much water from alternate sources they need to obtain (another reservoir, bottled water, etc), when to start testing the reservoir (no point wasting resources testing it too early), etc. The mathematics here is \textit{part of} the solution; it is not \textit{the} solution, and should not be relied upon as a definitive answer to a question as serious as the safety of drinking water.

\end{enumerate}
\fi

\newpage
\item \textbf{Topic: \textit{(Partial) Differential Equations}
 \\Allied disciplines: \textit{Engineering, Natural Sciences}
\\Year level: 2-3
    \\\hyperref[list:pillars]{Ethical pillars}: 1,7,9,10} 
    
    The following is the \textit{Neutron Diffusion Equation}, which governs how neutrons are produced and diffused in a uniform material: 
    \[
    \frac{\partial \Phi}{\partial t}= \frac{\sigma-1}{\tau}\Phi + \frac{\lambda^{2} }{3 \tau}\nabla^{2} \Phi
    \]
    where $\Phi(\mathbf{x},t)$ is the (free) neutron density, $\sigma$ is the average number of neutrons released in a fission event, $\tau$ is the average time between fissions, and $\lambda$ is the average distance a neutron travels before being absorbed by a nucleus.

    \begin{enumerate}
    \item (Using separation of variables, and observing the spherical symmetry of the system) Find the most general solution to Neutron Diffusion Equation for a uniform material arranged in a ball of radius $R$ (Hint: Using the substitution $\Phi = V/r$).
    \item Using the Dirichlet boundary condition $\Phi=0$ on the surface, find (in terms of $\sigma$ and $\lambda$) the \textit{critical radius} of the ball, for which the neutron density starts to grow exponentially in time.
    \item For Uranium-235, it can be experimentally determined that $\sigma\sim 2.3$ and $\lambda \sim  6$cm. Compute the critical radius of a ball of pure Uranium-235.
   
    \item What have you just done?
    \end{enumerate}
    
    \ifdefined\solution
    \solution\textit{This question leads to the radius of a critical mass of Uranium-$235$. ie; an atomic bomb. This is not a number that can be determined experimentally (at least, not in an efficient or safe way), but rather is determined from OTHER values that CAN be determined experimentally. This is the key piece of mathematics that gives some physics (of measuring some properties of Uranium-$235$) the predictive power of building an atomic bomb. The earlier parts of the question slightly obfuscate what is being done, but it is by no means ``hidden'', and one hopes that the students will read through to part (d) before starting the computation, and will at least reflect on the impact of such work before carrying it out (and if not at the start, then at least at the end when they reach (d) ).}
    
    \textit{This is exactly the mathematics that was taught to incoming scientists/researchers at the Manhattan Project during the Second World War, as it demonstrated the feasibility of building an atomic bomb, as well as a rough idea of how to build it. There was even a set of lecture notes given to new arrivals, called the Los Alamos Primer\footnote{See \href{https://en.wikipedia.org/wiki/Los_Alamos_Primer}{https://en.wikipedia.org/wiki/Los\_Alamos\_Primer}}.}

    \begin{enumerate}
    \item We are working in a sphere, so writing $\Phi(r,\theta, \phi,t)$ in spherical coordinates we use spherical symmetry to observe that there is no $\theta$ or $\phi$ dependence. ie: $\Phi=\Phi(r,t)$.  This is then a standard separation of variables question.

    The Laplacian in spherical coordinates is
    \[
    \nabla^{2} = \frac{1}{{r^2}} \frac{{\partial}}{{\partial r}} \left( r^2 \frac{{\partial}}{{\partial r}} \right) + \frac{1}{{r^2 \sin \theta}} \frac{{\partial}}{{\partial \theta}} \left( \sin \theta \frac{{\partial}}{{\partial \theta}} \right) + \frac{1}{{r^2 \sin^2 \theta}} \frac{{\partial^2}}{{\partial \phi^2}}
    \]
    As there is no $\theta$ or $\phi$ dependence, we can rewrite this as
    
    \[
    \nabla^{2} = \frac{1}{{r^2}} \frac{{\partial}}{{\partial r}} \left( r^2 \frac{{\partial}}{{\partial r}} \right) = \frac{\partial^{2}}{\partial r^{2}}+\frac{2}{r}\frac{\partial}{\partial r}
    \]
    So our PDE is now:
    \[
   \frac{\partial \Phi(r,t)}{\partial t}= \frac{\sigma-1}{\tau}\Phi(r,t) + \frac{\lambda^{2} }{3 \tau}\bigg(\frac{\partial^{2}\Phi(r,t)}{\partial r^{2}}+\frac{2}{r}\frac{\partial \Phi(r,t)}{\partial r} \bigg)
    \]

    Using our substitution trick $\Phi(r,t)=V(r,t)/r$  and the chain rule simplifies the Laplacian part of the equation to:
    \[
    \frac{\lambda^{2} }{3 \tau}\bigg(\frac{\partial^{2}\Phi}{\partial r^{2}}+\frac{2}{r}\frac{\partial \Phi}{\partial r} \bigg)=\frac{\lambda^{2} }{3 \tau}\bigg(\frac{\partial^{2}}{\partial r^{2}}\bigg(\frac{V}{r}\bigg)+\frac{2}{r}\frac{\partial }{\partial r}\bigg(\frac{V}{r}\bigg) \bigg)=\frac{\lambda^{2} }{3 \tau r}\frac{\partial^{2}V}{\partial r^{2}}
    \]
    So our full PDE for $V(r,t)$ is now:
    \[
   \frac{1}{r}\frac{\partial V(r,t)}{\partial t}= \frac{\sigma-1}{\tau r}V(r,t) +\frac{\lambda^{2} }{3 \tau r}\frac{\partial^{2}V}{\partial r^{2}}
    \]
    which, after cancelling $\frac{1}{r}$ becomes
     \[
   \frac{\partial V(r,t)}{\partial t}= \frac{\sigma-1}{\tau }V(r,t) +\frac{\lambda^{2} }{3 \tau }\frac{\partial^{2}V}{\partial r^{2}}
    \]
    We now use separation of variables to rewrite
    \[
    V(r,t)=X(r)Y(t)
    \]
    and so our PDE becomes the following ODE:
    \[
    X(r)Y'(t)=\frac{\sigma-1}{\tau }X(r)Y(t)+\frac{\lambda^{2} }{3 \tau }Y(t)X''(r)
    \]
    Rearranging gives
    \[
    \frac{Y'(t)}{Y(t)}=\frac{\sigma-1}{\tau }+\frac{\lambda^{2} }{3 \tau }\frac{X''(r)}{X(t)}=\frac{a}{\tau}
    \]
    where $a$ is some separation constant; as both sides are time-dependent functions that match for all $t$, they must be equal to a constant. Note that we chose $\frac{a}{\tau}$ for computational convenience later; most students will probably just choose a constant $c$, and it makes no difference at the end.

    Our two ODEs then become
    \begin{align}
    0 &=Y'(t)-\frac{a}{\tau}Y(t) 
    \\ 0 &= X''(r)+\frac{3(\sigma -a -1)}{\lambda^{2}}X(r)
    \end{align}

    For the first ODE (1), we assume a solution of the form $Y(t)=Ae^{kt}$. Then (1) becomes
    \begin{align*}
        0&= \frac{d}{dt}\bigg(Ae^{kt}\bigg)-\frac{a}{\tau}Ae^{kt}
        \\&=kAe^{kt}-\frac{a}{\tau}Ae^{kt}
        \\&=(k-\frac{a}{\tau})Ae^{kt}
    \end{align*}
    and so $k=\frac{a}{\tau}$, which then gives
    \[
    Y(t)=Ae^{\frac{a}{\tau}t}
    \]

    For the second ODE (2), we assume a solution of the form $X(r)=Be^{hr}$. Then (2) becomes
     \begin{align*}
         0&= \frac{d^{2}}{dr^{2}}\bigg(Be^{hr}\bigg)+\frac{3(\sigma -a -1)}{\lambda^{2}}Be^{hr}
         \\ &=h^{2}Be^{hr}+\frac{3(\sigma -a -1)}{\lambda^{2}}Be^{hr}
         \\&=(h^{2}+\frac{3(\sigma -a -1)}{\lambda^{2}})Be^{hr}
     \end{align*}
     which gives
     \[
     h^{2}+\frac{3(\sigma -a -1)}{\lambda^{2}}=0
     \]
     and thus
     \[
     h = \pm i \bigg( \frac{3(\sigma -a -1)}{\lambda^{2}}\bigg)^{\frac{1}{2}}
     \]
     For temporary convenience, let us set
     \[
     \delta :=\bigg( \frac{3(\sigma -a -1)}{\lambda^{2}}\bigg)^{\frac{1}{2}}
     \]
    So the general solution to (2) is now
    \[
    X(r)=De^{i\delta r}+Ee^{-i\delta r}
    \]
    where $F,G$ are constants.
    
    Expanding out $e^{i\delta r}=\cos(\delta r)+ i \sin(\delta r)$ and rearranging gives:
    \begin{align*}
        X(r)&=D(\cos(\delta r)+ i \sin(\delta r)) +E(\cos(\delta r)- i \sin(\delta r))
        \\&= F\cos(\delta r) +G \sin(\delta r)
    \end{align*}
    where $F,G$ are new (complex) constants.

    Combining all this, and using our original substitution $\Phi=\frac{V}{r}$, we get
    \begin{align*}
        \Phi(r,t)&=\frac{V(r,t)}{r}
        \\& = \frac{1}{r}Y(t)X(r)
        \\&= Ae^{\frac{a}{\tau}t}\bigg(F \frac{\cos(\delta r)}{r} + G \frac{\sin(\delta r)}{r}\bigg)
    \end{align*}
    Note that this must be defined on the entire domain, and in particular in the middle of the ball (at $r=0$). But in the limit as $r \to 0$, we have that $\frac{\cos(\delta r)}{r}\to \infty$. So the $F \frac{\cos(\delta r)}{r}$ term in our general solution must be removed, ie, $F=0$. So our general solution now becomes (with $\delta$ expanded out):
    \[
     \Phi(r,t)=\frac{H}{r}e^{\frac{at}{\tau}}\sin\bigg(  \frac{r}{\lambda}(3\sigma -3a-3)^{\frac{1}{2}}\bigg)
    \]
    where $H$ is some complex constant; $H, a$ are yet to be determined.

    \item The critical radius occurs when $a=0$ (if $a<0$ the density decays exponentially, if $a>0$ it grows exponentially). So we need to compute the radius $R$ for when this happens. Substituting in the boundary condition $\Phi(R,t)=0$ (with $a=0$) gives
    \[
    0=\frac{H}{R}\sin\bigg(  \frac{R}{\lambda}(3\sigma -0 -3)^{\frac{1}{2}}\bigg)
    \]
    And thus
    \[
    \frac{R}{\lambda}(3\sigma -3)^{\frac{1}{2}} = \pi
    \]
    Rearranging for $R$ gives
    \[
    R=\frac{\pi\lambda}{(3\sigma -3)^{\frac{1}{2}}}
    \]

    \item Direct substitution gives $R\sim 9.5$cm.

    \item By solving this equation and finding the critical radius you have enabled the construction of an atomic bomb, using principles of physics and some experimentally-determined values. You could go and investigate the difficulty of obtaining/refining that much Uranium, and whether it would be \textit{feasible} to obtain the amounts required (maybe it would take 100 years to obtain that much Uranium; maybe it would take 3 days. The difference here is the difference between being able to actually build such a bomb, and not being able to).
    
    \textit{Note to Instructor: at this point, it might be good to talk to your students about how they feel about having solved this problem. There is no obvious right or wrong answer here: some will be happy that they can understand something as powerful using simple mathematics, while others might feel sad when they see how easy (and thus dangerous) it is to solve such a task. This is a deep social, political, and philosophical issue, that the whole of humanity has been grappling with for the past 80 years; it is indeed quite complicated. Be open, be communicative, and understand that some students just had their own small Oppenheimer moment.}
    
    For further reading, Dan Jackson's articles Cantor's Paradise provide a good discussion involving both mathematics, physics and historical background.\footnote{
Dan Jackson (2023).  Derivation of the neutron-diffusion-equation. \textit{Cantor's Paradise}.
\\\href{https://www.cantorsparadise.com/derivation-of-the-neutron-diffusion-equation-64087afbdd20}{https://www.cantorsparadise.com/derivation-of-the-neutron-diffusion-equation-64087afbdd20}
\newline 
Dan Jackson (2023). Solving the neutron diffusion equation and calculating the critical mass. \textit{Cantor's Paradise}. \href{https://www.cantorsparadise.com/solving-the-neutron-diffusion-equation-and-calculating-the-critical-mass-of-uranium-cf90cc3bd16f}{https://www.cantorsparadise.com/solving-the-neutron-diffusion-equation-and-calculating-the-critical-mass-of-uranium-cf90cc3bd16f} \newline 
Dan Jackson (2023). How one miscalculation stopped Nazi Germany from developing the atomic bomb. \textit{Cantor's Paradise}. \href{https://www.cantorsparadise.com/how-one-miscalculation-stopped-nazi-germany-from-developing-the-atomic-bomb-3305b442cabc}{https://www.cantorsparadise.com/how-one-miscalculation-stopped-nazi-germany-from-developing-the-atomic-bomb-3305b442cabc}}

    \end{enumerate}

    \fi

\end{enumerate}
\newpage
\section{Exercises: Vector Calculus}
\begin{enumerate}
    \item \textbf{Topic: \textit{Vector Calculus}
     \\Allied disciplines: \textit{Engineering, Natural Sciences}
\\Year level: 1
    \\\hyperref[list:pillars]{Ethical pillars}: 2,7} 
    
Blood flowing in an artery is modelled as an axisymmetric flow in a cylinder of length $l$ and radius $R$ with axis parallel to $\hat{\mathbf{z}}$. The velocity is $\mathbf{v}(r,\theta,z) = v(r)\hat{\mathbf{z}}$, and satisfies
    $$
        \nabla^2 v = -\frac{p}{\eta l}, \quad \text{with boundary conditions} \quad v'(0) = v(R) = 0,
    $$
    where $p$ is the pressure difference between the ends of the artery, and $\eta$ is the viscosity of blood.

    \begin{enumerate}
    \item \textit{(Requires some knowledge of fluid dynamics)} What assumptions are we making in this model, and how important are they? 
    
\item  Find $v(r)$, and calculate the volumetric flow rate
    $$ Q = \iint_A \mathbf{v}\cdot\hat{\mathbf{z}} \dif A $$
    through a cross-section $A$ of constant $z$. If a build-up of cholesterol reduces the radius $R$ by a factor of $1/5$, by what factor does the pressure increase? 
    
    \item If a person's blood pressure is $50\%$ higher than the `normal' level then they are said to be in a \emph{hypertensive crisis} and require immediate medical attention. Compared to the radius of the artery, how thick would a layer of cholesterol need to be to trigger this? How sure are you of your answer?
    \end{enumerate}

    \ifdefined\solution
        \solution
    \textit{This is a modelling question, where the mathematics of model itself is quite straightforward, but the implications for the thing being modeled (in this instance, the health of a human) are incredibly serious. This is a real model, which describes the real health of real people. Sometimes it is easy for mathematicians to become detached from the consequences of their work, through the abstract nature of mathematics. This question is designed to help students realise that their work can have real consequences and effects, in a way they can all understand and relate to. It prompts them to question the assumptions they have used, implicitly or explicitly, and to evaluate these against the scale of the potential harms or outcomes.}

 \begin{enumerate}
 \item This model will assume things such as incompressibility of the fluid (Is blood \textit{actually} incompressible? How much does it change the results if it isn't, and how significant might those tiny changes be to a patient in a medical setting?). It will assume that the artery is a perfectly straight cylinder, with uniform radius (How many biological constructs are so precise? And again, what might the consequences be in terms of accuracy, and impact, if they are not?) It assumes blood has a uniform, constant viscosity (Will \textit{one} value of viscosity be used? Do people from different backgrounds, of different ages or genders or races, have the same blood viscosity?). It is important to have students articulate their assumptions, and be aware of the limitations in these.
 
\item	We have $\frac{dv}{d\theta} = \frac{dv}{dz} = 0$, hence
	\begin{align*}
	\nabla^2 v = \frac{1}{r} \frac{d}{dr}\left( r \frac{dv}{dr} \right) &= -\frac{p}{\eta l} \\
	r \frac{dv}{dr} &= c_1 -\frac{p}{2\eta l} r^2 ~~ \text{for some constant $c_1$,}
	\end{align*}
	and $\left.\frac{dv}{dr}\right\rvert_{r=0} = 0$ gives $c_1 = 0$:
	\begin{align*}
	\frac{dv}{dr} &= -\frac{p}{2\eta l} r \\
	v(r) &= \frac{p}{4 \eta l} (c_2 - r^2) ~,
	\end{align*}
	and $v(R) = 0$ gives $c_2 = R^2$:
	\begin{align*}
	v(r) = \frac{p}{4 \eta l} (R^2 - r^2) ~.
	\end{align*}
	This gives the volumetric flow rate through $A$ as
	\begin{align*}
	Q &= \frac{p}{4\eta l} \int_{\theta = 0}^{2\pi} \int_{r = 0}^{R} (R^2 - r^2) r dr d\theta \\
	&= \frac{p \pi}{2\eta l} \int_{0}^{R} (R^2 r - r^3) dr \\
	&= \frac{p \pi}{2\eta l} \left( \frac{1}{2} R^4 - \frac{1}{4} R^4 \right) \\
	&= \frac{p \pi}{8\eta l} R^4 ~.
	\end{align*}
	This gives
	\begin{equation*}
	p = \frac{8\eta l}{Q \pi R^4} ~,
	\end{equation*}
	thus if $R$ is reduced to $R_1 = \frac{4}{5} R_0$, we find
	\begin{align*}
	\frac{p_1}{p_0} &= \left(\frac{R_0}{R_1}\right)^4 \\
	& = \left(\frac{5}{4}\right)^4 \\
	& = \frac{625}{256} ~~(\approx 2.44) ~.
	\end{align*}
	That is, if $R$ decreases by a factor of $\frac{1}{5}$ the pressure more than doubles.
	
	\item If, on the other hand, we have $p_{1}=1.5p_{0}$, then 
	\[
	 1.5= \frac{p_1}{p_0} = \left(\frac{R_0}{R_1}\right)^4
	\]
	And thus $R_{1}=\frac{R_{0}}{(1.5)^{1/4}} \approx 0.9R_{0}$. So a layer of cholesterol as thick as 10\% of the radius of the artery is enough to trigger a serious medical condition.
 \\However, this assumes that the model is correct. We made many many assumptions, and even if those assumptions are by-and-large correct, we are only talking of small numbers here (10\% of the radius). Do the small errors in our assumptions compound to change this number significantly, to, say 5\% or 15\% of the radius? What might be the consequences for a patient receiving a diagnosis from a doctor? A slight underestimate and a patient with a serious medical condition might be given a clean bill of health. A slight overestimate and a completely healthy patient might be subjected to invasive and risky medical intervention unnecessarily.
\end{enumerate}
 
\end{enumerate}

\newpage 
\section{Exercises: Physics}

\begin{enumerate}
\item \textbf{Topic: \textit{Physics}
 \\Allied disciplines: \textit{Engineering, Natural Sciences}
    \\Year level: 1
    \\\hyperref[list:pillars]{Ethical pillars}: 1,5,9}
    
     Write a MATLAB function that plots
$$ y = y_0 + \frac{g}{k}t - \frac{g}{k^2}\left(1-e^{-kt}\right) $$
$$ x = \frac{v_0}{k}\left(1 - e^{-kt} \right)$$
in 2-dimensional space, the variables $k$, $v_0$ and $y_0$ should adjustable by the user. Plot this for $v_0 = 800$ and $y_0 = 4000$ and varying values of $k$.
Can you guess where these equations are coming from? Based on the variable names and your plots, try to reverse-engineer the source of the problem. What are the variables $v_0$, $y_0$ and $k$?

\ifdefined\solution
\solution
\textit{This question is designed to show students that they might be working in an information-constrained environment where their employer asks them to solve mathematical problems without giving them proper contextual knowledge. It teaches them to closely look at all variable names and produce plots of equations to reverse-engineer the source of the problem. The reverse-engineered problem is that of aerial bombing. Most students will only see it after having completed the task (i.e. coded the function). The variable $y_0$ is the height of the bomber, $v_0$ its velocity when dropping the bomb and $k = \frac{D\rho A}{2m}$ is coefficient of the drag force.} 

We created this question by using the work of Pepper and Wilson\footnote{Pepper, C. and Wilson, C. (2009). Aerial Bombing Techniques. Available online at \href{https://mse.redwoods.edu/darnold/math55/DEproj/sp09/WilsonPepper/CPandCWfinal.pdf}{https://mse.redwoods.edu/darnold/math55/DEproj/sp09/WilsonPepper/CPandCWfinal.pdf}}.

Given the complexity of including source code and plots, we omit these in this solution.

\newpage
\item \textbf{Topic: \textit{Dynamics} 
 \\Allied disciplines: \textit{Engineering, Natural Sciences}
\\ Year level: 1
\\\hyperref[list:pillars]{Ethical pillars}: 4,5}

A projectile of mass $m$ is fired from the (flat) ground with velocity $\mathbf{v}_{0}$, while a wind blows with constant velocity $\mathbf{u}$. The gravitational acceleration is $\mathbf{g}$ and the air exerts a drag force on the projectile equal to $mk$ times the velocity of the projectile \emph{relative to the wind}.
    \begin{enumerate}
        \item Write down and solve the equation model to find the trajectory of the projectile (assuming the projectile starts at the origin $\mathbf{x_0 = 0})$
        \item Suppose that $\mathbf{v}_{0}$ is vertical and $\mathbf{u}$ is horizontal. Show that the travel time $T$ of when the projectile hits the ground satisfies the implicit equation
        $$
        1 - e^{-kT} = \frac{kT}{1+\lambda},
        $$
        where $\lambda = kv_{0}/g$, for $g=|\mathbf{g}|$, $v_{0}=|\mathbf{v}_{0}|$. Verify that $\lambda$ is nondimensional, and give a physical interpretation of $\lambda$. 
       
        \item Although no analytic form is available for $T$, find an expression for the range $R$ of the projectile (that is, the distance it travels before it lands) in terms of $T$ and the other parameters.
        \item Real measurements can never be perfectly accurate or precise. Given that $u$ is uncertain, estimate the uncertainty in $T$ and $R$. \textit{Hint:} Linearise the equations about the true values by writing $u = u_0 + \delta u$, where $u_0$ is the true value and $\delta u$ is the uncertainty. Find an equation relating $\delta T$ and $\delta R$ to $\delta u$.

        \textit{Extension:} Now, estimate the error if $v_0$ is uncertain i.e. $v_0 = v_{00} + \delta v_{0}$ (you will need to use Taylor expansions to get the error estimate). 

        Why is it important to carry out uncertainty estimates in these calculations?

        \item This question has asked you to investigate \textit{projectile motion}. Can you list any uses of projectile motion that do not relate to military applications?

    \end{enumerate}
    
    \ifdefined\solution
    \solution\textit{Experimental scientists are well aware of the uncertainties that their experiments can have. On the other hand, mathematicians are somewhat suspicious of uncertainty, but they should not fall into the trap of thinking that `exact' answers are ever exact! This is especially important when the application, in this case the aiming of a projectile, could have serious effects on safety. The final part of the question shows how some areas of mathematics are immensely influenced by industries such as the military.}

\begin{enumerate}
    \item \textit{Note to instructor: This part can be given as a derivation exercise by stating the expression for $x(t)$ within the question, or as an exercise asking students to work out the expression fully for themselves.}
    \\The motion will be governed by Newton's Second Law $\mathbf{F}=m\mathbf{a}=m\mathbf{\dot v}$. Taking the origin to be the point where the projectile is fired from, we have:

$$m\mathbf{\dot v} = m\mathbf{g} -mk\left(\mathbf{v} - \mathbf{u}\right)$$ 
$$\Rightarrow m\mathbf{\dot v} =  -mk\mathbf{v} + mk\mathbf{u} + m\mathbf{g} $$
$$\Rightarrow \mathbf{\dot v} =  -k\mathbf{v} + k\mathbf{u} + \mathbf{g} $$
$$\Rightarrow \mathbf{\dot v} + k\mathbf{v} = k\mathbf{u} + \mathbf{g} $$
Using the integrating factor $e^{k t}$ we obtain
$$ \frac{d}{dt}\left( \mathbf{v} e^{kt} \right) = e^{kt}\left(k\mathbf{u} + \mathbf{g}\right)$$
$$\Rightarrow \mathbf{v} e^{kt} =  \frac{1}{k}e^{kt}\left(k\mathbf{u} + \mathbf{g}\right) + \mathbf{c}$$
$$\Rightarrow \mathbf{v} = \left(\mathbf{u} + \frac{\mathbf{g}}{k}\right) + e^{-kt}\mathbf{c} $$
But  $\mathbf{v}(0) = \mathbf{v}_{0}$, hence
$$\mathbf{c} =\textbf{v}_{0}-\mathbf{u} - \frac{\mathbf{g}}{k}$$
$$\Rightarrow \mathbf{v} = \left(\mathbf{u} + \frac{\mathbf{g}}{k}\right) + e^{-kt}\left(\textbf{v}_{0}-\mathbf{u} - \frac{\mathbf{g}}{k}\right) $$

$$\Rightarrow \textbf{x}=\left(\mathbf{u} + \frac{\mathbf{g}}{k}\right)t -\frac{1}{k}e^{-kt}\left(\textbf{v}_{0}-\mathbf{u} - \frac{\mathbf{g}}{k}\right)  + \mathbf{d} $$

But $\mathbf{x}(0) = \mathbf{0}$, hence
$$ \mathbf{d} = \frac{1}{k}\left(\textbf{v}_{0}-\mathbf{u} - \frac{\mathbf{g}}{k}\right) $$

So the position of the projectile at time $t$ is given by
$$\Rightarrow \textbf{x}=\left(\mathbf{u} + \frac{\mathbf{g}}{k}\right)t +\frac{1}{k}(1-e^{-kt})\left(\textbf{v}_{0}-\mathbf{u} - \frac{\mathbf{g}}{k}\right)  $$

\item Assume that the $z$-axis points up and the $x$-axis points into the direction of the wind. So $\mathbf{u}=(u,0,0)$,  $\mathbf{v}_{0}=(0,0,v_{0})$, and $\mathbf{g}=(0,0,-g)$. Then the equation transforms into
$$\begin{bmatrix}
x \\ y\\ z
\end{bmatrix} = \begin{bmatrix}ut \\ 0 \\ -gt/k \end{bmatrix} + \frac{1-e^{-kt}}{k}
\begin{bmatrix} - u \\ 0 \\ v_{0} + g/k \end{bmatrix}$$
The projectile returns to the ground when $z = 0$, i.e. when
$$\frac{gT}{k} = \frac{1-e^{-kT}}{k}\left(v_{0}+\frac{g}{k}\right)$$
$$\Leftrightarrow gT = \left(1-e^{-kT}\right) \left(v_{0}+\frac{g}{k}\right)$$
$$\Leftrightarrow \frac{gT}{v_{0}+\frac{g}{k}} = \left(1-e^{-kT}\right)  $$
$$\Leftrightarrow \frac{T}{\frac{v_{0}}{g}+\frac{1}{k}} = \left(1-e^{-kT}\right)  $$
$$\Leftrightarrow \frac{kT}{\frac{v_{0}k}{g}+1} = \left(1-e^{-kT}\right)  $$
Setting $\lambda = kv_{0}/g$ gives the desired result.
\\$\lambda$ is clearly non-dimensional, as $k$ is in units of $\frac{1}{time}$, $v_{0}$ is in units of $\frac{distance}{time}$, and $g$ is in units of $\frac{1}{time^{2}}$.
\\As a physical quantity, $\lambda$ is measuring the ratio of the drag and initial velocity, to the gravitational force.

\item The range $R$ of the projectile is the value of $x(T)$, where $T$ satisfies 
   $$
        1 - e^{-kT} = \frac{kT}{1+\lambda},
        $$
From the previous part, we have in general that
$$x(t)=  u\left(t-\frac{1 - e^{-kt}}{k}\right)$$
At $t=T$, we have $1 - e^{-kT} = \frac{kT}{1+\lambda}$, and so
$$R=x(T)=  u\left(T-\frac{kT}{k(1+\lambda}\right)=uT\left(1-\frac{1}{(1+\lambda)}\right)=\frac{uT\lambda}{(1+\lambda)}$$

\item

We start by treating each parameter and its uncertainty separately to estimate their effect on $T$ given by the implicit equation
   $$
        1 - e^{-kT} = \frac{kT}{1+\lambda},
        $$
and the explicit equation $$R=\frac{uT\lambda}{(1+\lambda)}$$. 

We start by considering $u = u_0 + \delta u$, which suggests that $T$ is unaltered and is still given by    
$$
        1 - e^{-kT} = \frac{kT}{1+\lambda},
        $$
as there is no $u$ term present in the equation for $T$ i.e. $T= T_0$ with $\delta T = 0$. \\

The equation for $R$ however, becomes $$R=\frac{(u_0 + \delta u) T\lambda}{(1+\lambda)},$$ which is of the form $R = R_0 + \delta R$, where $R_0=\frac{u_0 T\lambda}{(1+\lambda)}$ and $\delta R = \frac{\delta u T\lambda}{(1+\lambda)}$. 

For the extension part, we now consider $v_0 = v_{00} + \delta v_0$, and find that the equation for $T$ is altered as obtained from the new implicit equation and so can write its new value $T = T_0 + \delta T$ 
$$
        1 - e^{-k(T_0+\delta T)} = \frac{k(T_0+\delta T)}{1+\frac{k (v_{00} + \delta v_0)}{g}}=\frac{g(T_0+\delta T)}{\frac{g}{k}+v_{00} + \delta v_0} 
        $$
after replacing $\lambda = kv_{00} /g$.

Looking at the left hand side, if we Taylor expand $e^{-k\delta T}$ (expanding the function $e^{x}$ about the point $x=0$, to two terms) in the left hand side of this expression, we get
\[
1 - e^{-k(T_0+\delta T)} = 1-e^{-kT_0}e^{-k\delta T} \approx 1-e^{-kT_0}(1-k\delta T)
\]

If we Taylor expand $\frac{1}{\frac{g}{k}+v_{00} + \delta v_0} $ (expanding the function $\frac{1}{x}$ about the point $x=\frac{g}{k}+v_{00}$, to two terms) in the right hand side of this expression, we first observe that the Taylor expansion about $x=a$ would be
\[
\frac{1}{x} \approx \frac{1}{a}-\frac{1}{a^{2}}(x-a)
\]
So setting $a=\frac{g}{k}+v_{0}$ gives
\[
\frac{1}{\frac{g}{k}+v_{00} + \delta v_0} \approx \frac{1}{\frac{g}{k}+v_{00} } - \frac{1}{(\frac{g}{k}+v_{00})^{2}}\Big(\frac{g}{k}+v_{00} + \delta v_0 -(\frac{g}{k}+v_{00})\Big) = \frac{1}{\frac{g}{k}+v_{00} } - \frac{\delta v_{0}}{(\frac{g}{k}+v_{00})^{2}}
\]

Now combining these two approximations gives
\[
1-e^{-kT_0}+ke^{-kT_0}\delta T \approx \frac{gT_0}{\frac{g}{k}+v_{00} } - \frac{gT_0\delta v_{0}}{(\frac{g}{k}+v_{00})^{2}}+ \frac{g\delta T}{\frac{g}{k}+v_{00}} - \frac{g\delta T \delta v_{0}}{(\frac{g}{k}+v_{00})^{2}}
\]

Cancelling off the original expression for $T_0$ (as known from the original implicit equation), and neglecting the second-order error term $\delta T \delta v_{0}$ (as both errors are assumed small), we end up with

\[
ke^{-kT_{0}}\delta T \approx  \frac{g\delta T}{\frac{g}{k}+v_{00}} - \frac{gT_0\delta v_{0}}{(\frac{g}{k}+v_{00})^{2}}
\]

We can now solve this for $\delta T$ in terms of $T_0, v_{00}, \delta v_{0}, g, k$, which after significant algebraic manipulation yields

\[
\delta T \approx \frac{g T_0 \delta v_0}{\Big(\frac{g}{k} + v_{00}\Big) \Big(g - k e^{-k T_0} (g + k v _{00})\Big)}
\]

The equation for $R$ then becomes $$R\approx \frac{ u (T_0 + \delta T) (v_{00} + \delta v_0)}{\frac{g}{k} + v_{00} + \delta v_0},$$ with $\delta T$ given by its new value above. \\

We can use the same Taylor expansion for the denominator term from the calculation earlier to subsequently get 

$$R\approx \frac{ u} {\Big(\frac{g}{k} + v_{00}\Big)^2} \Big[(T_0 + \delta T) (v_{00} + \delta v_0)\Big(\frac{g}{k} + v_{00} - \delta v_0\Big)\Big],$$

which is of the form $R = R_0 + \delta R$, where $R_0= \frac{ u} {(\frac{g}{k} + v_{00})^2} T_0 v_{00} \Big(\frac{g}{k} + v_{00}\Big)$ and \\
$\delta R \approx \frac{ u} {(\frac{g}{k} + v_{00})^2} \Big((T_0 \delta v_0 + \delta T v_{00} )(\frac{g}{k} + v_{00}) - T_0 v_{00} \delta v_0\Big)$, after neglecting the second-order and third-order error terms. Substituting in $\delta T$ from above then gives
\begin{align*}
 \delta R&\approx \frac{ u} {(\frac{g}{k} + v_{00})^2} \Big((T_0 \delta v_0 + \frac{g T_0 \delta v_0}{\Big(\frac{g}{k} + v_{00}\Big) \Big(g - k e^{-k T_0} (g + k v _{00})\Big)} v_{00} )\Big(\frac{g}{k} + v_{00}\Big) - T_0 v_{00} \delta v_0\Big)
 \\&=\frac{ uT_{0}\delta v_{0}} {(\frac{g}{k} + v_{00})^2} \Big(\Big[1 + \frac{gv_{00} }{\Big(\frac{g}{k} + v_{00}\Big) \Big(g - k e^{-k T_0} (g + k v _{00})\Big)}  \Big]\Big(\frac{g}{k} + v_{00}\Big) - v_{00} \Big)
 \\&=\frac{ uT_{0}\delta v_{0}} {(\frac{g}{k} + v_{00})^2} \Big(\frac{g}{k} + v_{00} + \frac{gv_{00}(\frac{g}{k} + v_{00}) }{\Big(\frac{g}{k} + v_{00}\Big) \Big(g - k e^{-k T_0} (g + k v _{00})\Big)}  - v_{00} \Big)
 \\&=\frac{ ugT_{0}\delta v_{0}} {(\frac{g}{k} + v_{00})^2} \Big(\frac{1}{k} + \frac{v_{00} }{ \Big(g - k e^{-k T_0} (g + k v _{00})\Big)}   \Big)
\end{align*}

These error estimates are indeed inelegant; there is no reason to assume or expect that such estimates would be ``nice'' calculations.
\\

It is important to carry out uncertainty estimates in these calculations because, as mentioned above, no such calculation is ever exact. The mathematics might be exact, and might even have a closed-form ``exact'' solution. But there will always be errors in the initial measurements of parameters, initial conditions, etc. If such calculations are related to projectiles (possibly of explosives such as weapons, fireworks, etc), then error margins need to be understood \textit{before} firing the projectile. \\

Moreover, the model in use might not properly reflect the physical system, or might over-simplify and ignore important terms or other parameters. For additional discussion, you can discuss with students what might happen for long-distance projectiles which must also account for the rotation of the earth; ask them to estimate the lengthscale over which rotation becomes important. \\

\item The purpose of this part of the question is to help students see that, up to now, most of the uses of projectile motion (an interesting, abstract, problem that is tractable at the early undergraduate level) come directly from military applications; in particular, firing guns or dropping bombs from an aeroplane. There are other uses, such as:
\begin{itemize}
    \item Calculating the path of emergency supplies being dropped from an aeroplane into a natural disaster recovery zone.
    \item Computing the final resting place of an object or ship that has sunk in the ocean. This is projectile motion where the fluid medium is water rather than air.
    \item Computing the return path for spacecraft re-entry to Earth.
    \item Computing the path, and hence range, of a water jet coming out of a firehose. Here the projectile is water, and the medium is air, but water is still ``a thing moving in air'', so similar principles apply.

\end{itemize}

\end{enumerate}

\fi  
    
\end{enumerate}

\newpage
\section{Exercises: Optimisation}\label{optimisation}
\begin{enumerate}
    \item \textbf{Topic: \textit{Optimisation}
     \\Allied disciplines: \textit{Economics, Statistics, Computer Science}
    \\Year level: 1
    \\\hyperref[list:pillars]{Ethical pillars}: 5,6} 
    
    Let $S$ be an n-dimensional random vector of asset returns. Let $U: \mathbb{R}\rightarrow \mathbb{R}$ be the concave utility function of the investor and let $r:\mathbb{R}^n \rightarrow \mathbb{R}_{\geq 0}$ be the convex risk function associated with an investment $x \in \mathbb{R}^n$. 
\begin{enumerate}
    \item Set up an optimisation problem that maximises the investor's expected utility given an upper bound $b \geq 0$ for the risk. You may assume that the investor has no other monetary constraints.
    \item Set up an optimisation problem that minimises the risk given a lower bound $u \geq 0$ on the expected utility. You may assume that the investor has no other monetary constraints.
    \item You can assume that the optimal solutions of (a) and (b) all lie on the same frontier in space. By choosing the appropriate values for $b$ and $u$ the optimisation problems lead to the same values. Are the optimisation problems the same?
    Put into words how we use and understand risk in each case. Are we using risk in the same way?
\end{enumerate}
    \ifdefined\solution
    \solution \textit{This question is designed for students to explore different optimisation models. It shows them that they can solve one problem using two different approaches, they can mathematically get the same result. But the approaches differ at a psychological and communicative level and thus elicit different responses from managers, customers, etc.}
    
\begin{enumerate}
    \item 
        \begin{equation*}
        \begin{aligned}
        & \underset{x}{\text{maximize}}
        & & \mathbb{E}\left[U\left(S^Tx\right)\right] \\
        & \text{subject to}
        & & r\left(x\right) \leq b \\
        \end{aligned}
        \end{equation*}
    \item 
        \begin{equation*}
        \begin{aligned}
        & \underset{x}{\text{minimize}} 
        & &  r\left(x\right) \\
        & \text{subject to}
        & & \mathbb{E}\left[U\left(S^Tx\right)\right] \geq u \\
        \end{aligned}
        \end{equation*}
    \item 
     In the first case we see risk as a constraint and in the second case minimizing risk is the objective. The student should understand that equivalent mathematical models can elicit different responses in everyday situations when the underlying variables, terms and ideas are used in different ways. Model (a) is commonly understood as more ``optimistic'', while model (b) is typically seen as more ``pessimistic''. Many people intuitively prefer to maximise their utility function, instead of minimising risk. To many, it feels as if they get more out of it. 
\end{enumerate}

\newpage
    \item\label{oil-pipe} \textbf{Topic: \textit{Optimisation}
     \\Allied disciplines: \textit{Economics, Computer Science, Engineering, Natural Sciences}
    \\Year level: 1
    \\\hyperref[list:pillars]{Ethical pillars}: 1,5,9}

    An oil company wants to build a pipeline connecting an offshore oil platform to a refinery (on land). The coastline is straight. The oil platform is a distance $D_1$ from the coast. The refinery is on the coastline, distance $D_2$ from the point on the coast closest to the platform. 
    Building the pipeline has a cost per unit length $c_1$ at sea, and $c_2$ on land. How should the pipeline be built? What are the factors that need to be considered when providing a response to this question?
    \\\textit{Here is an alternative formulation of the same problem, inspired by Shulman:}
    \\An offshore oil well platform is located at point A, which is a distance $D_1$ km from the nearest point Q on a straight shoreline. The oil is to be piped from A to a refinery at a point T on the shoreline by piping it straight under water to a point P on the shoreline between Q and T and then to T by a pipe along the shoreline. Suppose that the distance from Q to T is $D_2$ km and that it costs $c_1$ per km to lay an underwater pipe and $c_2$ per km to lay the pipe along the shoreline. What should the distance $x$ of the pipeline from P to T be? What are all the factors that need to be considered when building this pipeline?

    \ifdefined\solution
    \solution \textit{The cost-minimising path is given by Snell's law, and is an exercise in trigonometry. But who said we were optimising over cost? This is an assumption that is trained in to mathematicians while they are students, but it need not always be the right way to optimise.
    Economic actions almost always have externalities, such as possible damage to the environment (the pipe may go through a coral reef or protected habitat) or to existing infrastructure (it may go through a school, or a site of archaeological significance). Policymakers have to take a holistic view of these effects, and the mathematicians that advise them need to be honest about the simplifications and the assumptions that go into their models.}

    \textit{The idea for this question comes from Shulman\footnote{B. Shulman (2002). Is There Enough Poison Gas to Kill the City?: The Teaching of Ethics in Mathematics Classes. \textit{The College Mathematics Journal} 33(2), pp. 118-125. \href{https://www.tandfonline.com/doi/abs/10.1080/07468342.2002          .11921929           }{https://www.tandfonline.com/doi/abs/10.1080/07468342.2002.11921929                             }}.}

The cost-minimising path (ASSUMING there are no other associated costs) will be a diagonal straight line under water from the platform at A to the point P on the shoreline, then a straight line down the shoreline to the refinery from P to T.
So let us take ANY such pipe layout (diagonal under water, then straight along the shoreline). Let $x$ be the length of the pipe along the shoreline from P to T. Then the length of the (diagonal) pipe under water, $L$ from A to P is given by Pythagoras:
\[
L=\sqrt{D_{1}^{2}+(D_{2}-x)^{2}}
\]
Thus the total cost of the pipe in terms of $x$, $C(x)$, is given by
\[
C(x)= c_{2}x + c_{1}\sqrt{D_{1}^{2}+(D_{2}-x)^{2}} 
\]
Differentiating with respect to $x$, and looking for a minima, we get
\[
0=\frac{dC}{dx}=c_{2} - \frac{c_{1}(D_{2}-x)}{\sqrt{D_{1}^{2}+(D_{2}-x)^{2}}}
\]
And thus
\[
\frac{c_{2}}{c_{1}}=\frac{(D_{2}-x)}{\sqrt{D_{1}^{2}+(D_{2}-x)^{2}}}
\]
Let $\theta$ be the angle the underwater pipe makes with the shoreline. Observe that we have 
\[
\cos(\theta)=\frac{(D_{2}-x)}{\sqrt{D_{1}^{2}+(D_{2}-x)^{2}}} = \frac{c_{2}}{c_{1}}
\]
Thus
\[
\frac{D_{1}}{D_{2}-x}=\tan(\theta)=\tan(\arccos(\frac{c_{2}}{c_{1}}))
\]

Rearranging gives
\[
x=D_{2}-\frac{D_{1}}{\tan(\arccos(\frac{c_{2}}{c_{1}}))}
\]
so the minimum cost pipe has 2 segments:
\\Segment 1: a straight pipe on the shore, going up distance $x$ from the refinery.
\\Segment 2: a straight pipe under water, from the end of the segment 1 pipe to the platform.

But what we have computed is ``the minimum cost pipe, assuming there are no other costs involved, and not factoring in any other impacts or concerns'', which is \textit{not} what the question asks. The question asks ``How should the pipeline be built?''; this is not a mathematics problem, it is a problem about building a pipe. Mathematics is a useful tool in providing an answer, but there are many other things to consider. Figuring out where and how to build the pipe involves many considerations, including: What is under the sea (coral reef, etc)? What is on the shoreline (protected habitat, important buildings, areas of cultural or religious significance, etc)? What are the risks of having a pipeline in that area (oil spills, etc)? What are the ongoing maintenance costs of underwater pipes, compared to overland (rust, difficulty in fixing leaks, etc)? How susceptible to damage or sabotage is an underwater pipe, compared to overland (damage from protest groups, accidental impact from vehicles, etc)? What are the political consequences of building an oil refinery (protests, reputational damage, etc)? Is there a better thing to be building instead (an offshore wind farm)?

There are numerous factors to consider here, more than can be listed. The point is to get students to realise that there are at least SOME factors, and that optimising over (naive) cost, with no other considerations, is not at all a good way to answer the question.

\fi

\newpage
     \item \textbf{Topic: \textit{Optimisation}
    \\Allied disciplines: \textit{Economics, Statistics, Computer Science}
    \\Year level: 1
    \\\hyperref[list:pillars]{Ethical pillars}: 5,6}
    
\begin{enumerate}
    \item You have just been hired as an intern by an investment firm and are given $n$ possible investment opportunities, each having a deterministic profit $p_i$ but random investment costs $a_i$. You are told that you have a fixed budget constraint $b$. Someone else has already modelled this for you as a Knapsack problem using the mean $\mu_i$ of $a_i$.
    \begin{equation*}
\begin{aligned}
& \underset{x}{\text{maximize}}
& & \sum\limits_{i=1}^n p_ix_i \\
& \text{subject to}
& & \sum\limits_{i=1}^n \mu_i x_i \leq b \\
& 
& & x_i \in \{0,1\}, \; i = 1, \ldots, n.
\end{aligned}
\end{equation*}
Briefly discuss why this might not be a good idea.
\item A simple way to tackle the randomness of the problem is by considering chance constraints. In this case we accept that side conditions are violated with a probability $\epsilon.$
\begin{equation*}
\begin{aligned}
& \underset{x}{\text{maximize}}
& & \sum\limits_{i=1}^n p_ix_i \\
& \text{subject to}
& & \mathbb{P}\left(\sum\limits_{i=1}^n a_i x_i \leq b\right) \geq 1- \epsilon \\
& 
& & x_i \in \{0,1\}, \; i = 1, \ldots, n.
\end{aligned}
\end{equation*}
Show that it can be rewritten as 
\begin{equation*}
\begin{aligned}
& \underset{x}{\text{maximize}}
& & \sum\limits_{i=1}^n p_ix_i \\
& \text{subject to}
& & \Phi^{-1}\left(1-\epsilon\right) \sqrt{\sum\limits_{i=1}^n \sigma_i^2 x_i^2} + \sum\limits_{i=1}^n \mu_i x_i \leq b \\
& 
& & x_i \in \{0,1\}, \; i = 1, \ldots, n.
\end{aligned}
\end{equation*}
You can assume that the $a_i$ are independent and follow a normal distribution $a_i \sim N\left(\mu_i, \sigma_i^2\right)$. How would you describe the special form of the first constraint in the second version in simple language?
\end{enumerate}

\ifdefined\solution
    \solution
    \textit{The question is designed to show students that ignoring randomness and assuming idealised deterministic conditions can lead to extreme costs in optimisation conditions. Students should learn that this can often be avoided through very simple concepts. The student should understand that it may be necessary to question the validity of an optimisation model they are given by a 3rd party, and not blindly rush to solve it.}
\begin{enumerate}
    \item The classical formulation of the knapsack problem assumes that the $a_i$, $p_i$ and $b_i$ are known and fixed an known in advance. The student should intuitively understand that if randomness is involved, naively solving the classical problem might no longer be a good idea. Random variations of parameters can lead to large violations of the constraints.
    
    Note to Instructor: Ben-Tal et al.\footnote{Ben-Tal, A., and Nemirovski, A. (2000) Robust solutions of linear programming problems
contaminated with uncertain data. \textit{Mathematical Programming} 88(3),pp. 411–424.} found in their study of the NETLIB library of optimization problems, that a deviation in the coefficients of as little as $0.1\%$ can lead to the violation of the inequality constraints by $450\%$. In $18\%$ of the cases the violations were larger than $150\%$. The classical problem formulation is often unusable in practice when the parameters are uncertain.
    \item We can rewrite the constraint by considering
\begin{align*}
    \mathbb{P}\left(\sum\limits_{i=1}^n a_i x_i \leq b\right) &= \mathbb{P}\left(\frac{\sum\limits_{i=1}^n a_ix_i - \sum\limits_{i=1}^n \mu_ix_i}{\sqrt{\sum\limits_{i=1}^n \sigma_i^2x_i^2}} \leq \frac{b - \sum\limits_{i=1}^n \mu_ix_i}{\sqrt{\sum\limits_{i=1}^n \sigma_i^2x_i^2}} \right) \\ &= \mathbb{P}\left(Z \leq \frac{b - \sum\limits_{i=1}^n \mu_ix_i}{\sqrt{\sum\limits_{i=1}^n \sigma_i^2x_i^2}} \right) \geq 1 - \epsilon
\end{align*}
where $Z \sim N\left(0,1\right)$. After using the cumulative distribution function $\Phi$ of the standard normal distribution and some rearranging, we obtain the solution:
$$\Phi\left(\frac{b - \sum\limits_{i=1}^n \mu_ix_i}{\sqrt{\sum\limits_{i=1}^n \sigma_i^2x_i^2}} \right) \geq 1- \epsilon$$
$$\Rightarrow \frac{b - \sum\limits_{i=1}^n \mu_ix_i}{\sqrt{\sum\limits_{i=1}^n \sigma_i^2x_i^2}} \geq \Phi^{-1}\left(1-\epsilon\right).$$
$$\Rightarrow  b - \sum\limits_{i=1}^n \mu_ix_i \geq \left(1-\epsilon\right)  \sqrt{\sum\limits_{i=1}^n \sigma_i^2x_i^2} \Phi^{-1}$$
$$\Rightarrow \Phi^{-1}\left(1-\epsilon\right) \sqrt{\sum\limits_{i=1}^n \sigma_i^2 x_i^2} + \sum\limits_{i=1}^n \mu_i x_i \leq b$$
Taking a closer look at the constraint reveals its special form (in simple language):
$$ \text{penalty from randomness } + \text{ classical cost when using the mean $\mu_i$ of $a_i$} \leq b.$$
\end{enumerate}

\end{enumerate}
\newpage
\section{Exercises: Probability}
\begin{enumerate}
\item \textbf{Topic: \textit{Probability}
    \\Allied disciplines: \textit{Economics, Statistics, Computer Science}
    \\Year level: 1
    \\\hyperref[list:pillars]{Ethical pillars}: 2,5,6}
    
You are on holiday with your friend Kelly, and the two of you go hiking. Unfortunately, Kelly falls on a rock and breaks her right foot. She is transported to a small remote hospital and needs to stay at the hospital for another night. The power of the hospital fails each hour independently with probability $p$.
\begin{enumerate}
\item Derive a formula for the expected number of hours until the hospital loses power. 
\item How small does $p$ have to be, such that the expected number of power failures in a given year is less than 1?
\item Kelly has to stay in the hospital for another $x$ hours. If the power does not fail, she will be discharged in time and the two of you still catch your return flight. What is the probability that the hospital loses power in the first $x-1$ hours?
\item Kelly asks you what value of $p$ would make you feel comfortable. How do you decide? 
\end{enumerate}

\ifdefined\solution
    \solution
    \textit{This question is designed to show students that mathematics is not always black and white, but that there are many subjective components in setting up models and choosing values.}
\begin{enumerate}
\item Let $N$ be the number of hours until the hospital first loses power. Then we need to calculate $E\left[ N\right] = \sum\limits_{i=0}^\infty P\left( N > i\right).$ The probability that the power does not fail in the first $i$ hours is given by $\left( 1-p\right)^i$. Hence,
$$E\left[N\right] = \sum\limits_{i=0}^\infty \left(1-p\right)^i = \frac{1}{1-(1-p)} = \frac{1}{p}$$
\item The expected number of power failures in a year is given by $$365\cdot 24 \cdot p.$$ Solving $365 \cdot 24 \cdot p < 1$ gives $p < \frac{1}{365 \cdot 24}$.
\item The probability that the power fails in the $i$-th hour for the first time is given by $\left(1-p\right)^{i-1}p$. Hence, 
$$ P\left(N < x\right) = \sum\limits_{i=1}^{x-1} P\left(N = i\right) = \sum\limits_{i=1}^{x-1} \left(1-p\right)^{i-1}p.$$
\item The final part of the question is designed to test if the student gives a solution only considering their friend, the chance of catching the return flight or if they also think about the other patients. As mathematicians, we often need to choose values for variables that make us feel comfortable or seem right. The question is designed to show students that outside of pure mathematics, ``finding a value that feels right'' can be tricky. Some will decide to optimise $p$ according to either (a), (b) or (c), others might compare different solutions or tell you that this ``question is not well defined''. When mathematics is applied in a social context, questions are rarely well-defined and leave plenty of room for different interpretations. 
\end{enumerate}
\fi

\newpage
    \item \textbf{Topic: \textit{Probability}
    \\Allied disciplines: \textit{Economics, Statistics, Computer Science}
    \\Year level: 1
    \\\hyperref[list:pillars]{Ethical pillars}: 5}
    
    Two players A and B flip a biased coin alternately with A going first, and the first player to flip a heads wins the game. Each flip is independent of all previous flips, and has a probability $p$ of coming up heads. Find an expression in terms of $p$ for $\Pr(A)$, the probability that A wins, and sketch $\Pr(A)$ against $p$. What happens as $p\rightarrow1$? What about at $p=1$? What happens as $p\rightarrow0$? What about at $p=0$?
    
    \ifdefined\solution
    \solution \textit{The purpose of this question is to help mathematicians see that their mathematics might still yield a solution outside their domain of definition of the problem. In this instance, there is no winner when $p=0$, even though in the limit as $p\rightarrow0$ the chance of player A winning approached $1/2$. It is easy to forget what the original problem was, and its constraints, after it has been `mathematised', and this can lead to a mathematical model being misused and giving nonsensical outputs that are harmful if relied upon in the physical world.}

    This can be solved by seeing that the chance of B winning is the chance of A flipping tails x the chance of A winning (as, after flipping a tail, A passes the coin to B, and B now has the same chance of winning as A originally did). So this can be solved using two simultaneous equations:
    $$\Pr(A)+\Pr(B)=1$$
    $$\Pr(B)=(1-p)\Pr(A)$$
    \\Which gives $$\Pr(A)=\frac{1}{2-p}.$$  (This can also be solved by summing an infinite series; it gives the same answer). But this is only valid if there is an actual chance of heads; that is, $0<p\leq 1$. If $p=0$, then the equation $\Pr(A)+\Pr(B)=1$ is not valid, as no-one can win, so in that case $P(A)=0$.
    
    In the limiting case $p \to 1$ we have $\lim_{p\to 1}\frac{1}{2-p} = 1$. That is, as $p$ gets closer to $1$, A stands an overwhelming chance of winning. And moreover it is clear that when $p=1$ we have $P(A)=1$ (as soon as A flips the coin, they will get heads and win).

    In the limiting case $p \to 0$ we have $\lim_{p\to 0}\frac{1}{2-p} = \frac{1}{2}$. That is, as $p$ gets closer to $0$, the game becomes closer and closer to being perfectly fair (the advantage of going first is greatly reduced). But, when $p=0$ we have $P(A)=0$, as it is not possible for A to ever flip heads. The limiting behaviour of the model is is no way indicative of the outcome of the system at that limit point.

    This is a very important lesson: after mathematising, the initial assumptions, domain of definition, etc can be easily lost or forgotten, and the model misused as a result. Just because the expression $\frac{1}{2-p}$ has a well-defined limit as $p\to 0$, and is even defined at $p=0$, does not mean it is indicative of the physical behaviour of the system at $p=0$.

\fi

\newpage
\item \textbf{Topic: \textit{Probability}
    \\Allied disciplines: \textit{Economics, Statistics, Computer Science}
    \\Year level: 1
    \\\hyperref[list:pillars]{Ethical pillars}: 4,5,6}
    
\begin{enumerate}
   
 \item  Consider the following methods of randomly choosing a chord in a circle:

	\begin{enumerate}
		\item Fix a point on the circumference as one end of the chord, and choose an angle $\Theta \sim \text{Unif}(0, \pi)$ to be the angle subtended by the chord at the centre.

		\item Fix a radius to intersect the chord at right-angles, and choose a distance from the centre $H \sim \text{Unif}(0, 1)$ to be the midpoint of the chord.

		\item Choose a uniformly distributed point in the circle $(X,Y)$ to be the midpoint of the chord (to do this, choose uniform variables $D \sim \text{Unif}(0,1)$, $\Theta \sim \text{Unif}(0,2\pi)$, then set $X = \sqrt{D}\cos\Theta$, $Y = \sqrt{D}\sin\Theta$).

	\end{enumerate}

	For each of these, find the distribution function $F_L (l)$, and the corresponding probability density function $f_L (l)$, for the length $L$ of the chord. In each case evaluate the probability that $L$ is greater than the side length of an inscribed equilateral triangle.

\item How would you answer the following question: ``Given a circle of unit radius, what is the probability that a random chord of the circle is longer than the sides of an inscribed equilateral triangle?''
\end{enumerate}

    \ifdefined\solution
    \solution \textit{This is a crucial observation for students studying probability; how we define ``random'' determines/changes what we can deduce about random processes. To say ``A randomly chosen item from a particular set has a $50\%$ probability of possessing property X'' is a meaningless statement} unless \textit{it is fully qualified with} how \textit{that item was ``chosen randomly''. This is one of the ways mathematics can, and is, used to dupe the public, especially in when probability theory and statistics are involved.}

\textit{The final question is (deliberately) not well-formulated, because there are many ways in which a chord could be chosen `randomly'. When we informally say `pick a number from 1 to 10 randomly', we usually implicitly mean `under a uniform distribution'; that is, $Pr(n)=0.1$ for all $n \in \{1, \ldots, 10\}$. This interpretation does not work in this problem, where a chord has multiple degrees of freedom. }

    \begin{enumerate}

    	\item Computing the probability functions and densities: 
	\begin{enumerate}
		\item We have $L = 2 \sin(\Theta /2)$, thus
		\begin{align*}
		F_L (l) &= \mathbb{P}\lbrace L \leqslant l \rbrace = \mathbb{P}\lbrace 2 \sin(\Theta /2) \leqslant l \rbrace \\
		&= \mathbb{P}\lbrace \Theta \leqslant 2 \arcsin(l/2) \rbrace = F_\Theta (2 \arcsin(l/2)) \\
		&= \frac{2}{\pi} \arcsin(l/2) ~~~ (0 < l < 2) ~,
		\end{align*}
		giving
		\begin{equation*}
		f_L (l) = \frac{d}{dl} F_L (l) = \frac{1}{\pi} \frac{1}{\sqrt{1-l^2/4}} ~~~ (0 < l < 2) ~.
		\end{equation*}
		\item We have $L = 2\sqrt{1 - H^2}$, giving
		\begin{align*}
		F_L (l) &= \mathbb{P}\lbrace L \leqslant l \rbrace = \mathbb{P}\lbrace 2\sqrt{1-H^2} \leqslant l \rbrace \\
		&= \mathbb{P}\lbrace H \geqslant \sqrt{1 - l^2/4} \rbrace = 1 - F_H (\sqrt{1 - l^2/4}) \\
		&= 1 - \sqrt{1-l^2/4} ~~~ (0 < l < 2) ~,
		\end{align*}
		and
		\begin{equation*}
		f_L (l) = \frac{d}{dl} F_L (l) = \frac{l}{2} \frac{1}{\sqrt{1-l^2/4}} ~~~ (0 < l < 2) ~.
		\end{equation*}
		\item We have $L = 2\sqrt{1-(X^2 + Y^2} = 2\sqrt{1-D}$, giving
		\begin{align*}
		F_L (l) &= \mathbb{P}\lbrace L \leqslant l \rbrace = \mathbb{P}\lbrace 2\sqrt{1-D} \leqslant l \rbrace \\
		&= \mathbb{P}\lbrace D \geqslant 1 - l^2/4 \rbrace = 1 - F_D (1 - l^2/4) \\
		&= l^2/4 ~~~ (0 < l < 2) ~,
		\end{align*}
		and
		\begin{equation*}
		f_L (l) = \frac{d}{dl} F_L (l) = \frac{l}{2} ~~~ (0 < l < 2) ~.
		\end{equation*}
	\end{enumerate}
	The probability that $L$ is greater than the side length of an inscribed equilateral triangle is $1 - F_L(\sqrt{3})$, giving the results:
	\begin{enumerate}
		\item $1 - F_L(\sqrt{3}) = 1 - \frac{2}{\pi} \arcsin\left(\frac{\sqrt{3}}{2}\right) = 1 - \frac{2}{3} = \frac{1}{3} ~.$
		\item $1 - F_L(\sqrt{3}) = \sqrt{1 - \frac{3}{4}} = \frac{1}{2} ~.$
		\item $1 - F_L(\sqrt{3}) = 1 - \frac{3}{4} = \frac{1}{4} ~.$
	\end{enumerate}

 \item The difficulty with this question is that even when we interpret the term `random' -- here as some form of `uniformly distributed' -- it is still not enough to determine \textit{how} exactly to sample a chord. The \textit{mathematical choice} of which definition of random to use will heavily influence the outcome of the computation, as we saw above.

One might see this more clearly when sampling between sets. For ease of measurement, a set S of people might be split into two subsets A,B. Then, to sample height from the whole set S, we might uniformly choose between A and B (probability 0.5 in each case), and then uniformly choose one member from that set and measure their height. If set B had fewer people than set A, but had a lot of tall people, then this would give a different (taller) average than just sampling uniformly from the whole set S.

    \end{enumerate}
    
    \fi

   \newpage 
    \item \textbf{Topic: \textit{Probability}
    \\Allied disciplines: \textit{Statistics, Computer Science, Natural Sciences}
    \\Year level: 1
    \\\hyperref[list:pillars]{Ethical pillars}: 2,5}
    
    Colourblindness is usually caused by a defective gene in the X chromosome. Let the defective allele, which is recessive, be written as `x', and the healthy allele, which is dominant, be written as `X'. Most women have two X chromosomes, while most men have one X and one Y chromosome. Epidemiologists have determined that, globally, $p \approx 8\%$ of men are colourblind.
    
    \begin{enumerate}
        \item If a man is not colourblind, what is his genotype? If a woman is colourblind, what is her genotype?
        
        \item Approximately what proportion of women are colourblind? Approximately what proportion women are carriers of the colourblindness gene, but not colourblind themselves?

        \item Neither of Sam's parents are colourblind. What is the probability that Sam's child is colourblind?

        \item What other issues do you need to consider to provide meaningful and helpful probabilities here, and why?
    \end{enumerate}

    \ifdefined\solution
    \solution\textit{This question is all about implicit assumptions, how they can be erroneous, and the complications that can arise from them. Implicit assumptions (a particular name implies a particular gender) can seem trivial to solve or address, from a mathematical perspective, but only IF they are actually identified. If they remain undetected, then it doesn't matter how easy they would have been to address; they remain unaddressed, and the problem the create thus persists.}
    
    \begin{enumerate}
        \item A non-colourblind man must have a non-defective X chromosome, and his genotype must be XY. A colourblind woman must have two defective chromosomes, so her genotype must be xx.
    
        \item Assuming that birth and death rates are independent of colourblindness (so that the value of $p$ does not change between generations), about $p$ of all X chromosomes carry the defective allele. Therefore, a woman has probability $p^2\approx0.6\%$ of having genotype xx, and probability $2p(1-p) \approx 15\%$ of being a carrier, with genotype Xx.

        \item Did the student assume that Sam is male? (It is important to be aware of one's implicit biases, such as `male by default', especially in a community that is majority-male.) Therefore both cases must be considered to account for the fact that Sam may either be male or female.
        
        Sam's father must have genotype XY, and their mother's genotype is distributed with probabilities
        \begin{center}
        \begin{tabular}{ccc}
             XX (healthy) & Xx (carrier) & xx (colourblind) \\ \hline\hline
             $ \displaystyle \frac{1-p}{1+p}$ & 
             $ \displaystyle \frac{2p}{1+p} $ & $0$
        \end{tabular}
        \end{center}
        when we exclude the possibility that she is colourblind. 
        
        (1) If Sam is male, then it is the X chromosome that he receives from his mother that determines whether Sam is colourblind. Sam can only get the defective chromosome if his mother is a carrier, so his genotype distribution is
        \begin{center}
            \begin{tabular}{cc}
            XY (healthy) & xY (colourblind) \\\hline\hline 
            $\displaystyle \frac{1}{1+p}$ 
            & 
            $\displaystyle \frac{p}{1+p}$ 
            \end{tabular}
        \end{center}
        (So, he has a better chance than the population of not being colourblind.) We don't know anything about Sam's partner, so we assume her genotype follows the population distribution for females
        \begin{center}
            \begin{tabular}{ccc}
            XX (healthy) & Xx (carrier) & xx (colourblind) \\ \hline\hline
            $(1-p)^2$ & $2p(1-p)$ & $p^2$
            \end{tabular}
        \end{center}
        Now consider Sam's children. For a male child, the genotype depends purely on the mother, and does not depend on Sam (or Sam's parents), so Sam's son has genotype distribution
        \begin{center}
            \begin{tabular}{cc}
            XY (healthy) & xY (colourblind) \\\hline\hline
            $1-p$ & $p$
            \end{tabular}
        \end{center}
        \textit{i.e.} no different from a population male. For a female child, we begin by calculating the probability that she is colourblind: this can occur only if Sam is colourblind, and his partner has at least one copy of the defective allele:
        $$ \Pr(\text{Sam's daughter is xx}) = \frac{p}{1+p} \left[p^2 + p(1-p) \right] = \frac{p^2}{1+p}. $$
        Similarly,
        $$
        \Pr(\text{Sam's daughter is XX}) = \frac{1}{1+p} \left[(1-p)^2 + p(1-p)\right] = \frac{1-p}{1+p}.
        $$
        Hence the genotype distribution for Sam's daughter is
        \begin{center}
            \begin{tabular}{ccc}
            XX (healthy) & Xx (carrier) & xx (colourblind) \\\hline\hline
            $\displaystyle\frac{1-p}{1+p}$ & 
            $\displaystyle\frac{2p-p^2}{1+p}$ 
            &             $\displaystyle\frac{p^2}{1+p}$
            \end{tabular}
        \end{center}
    
        (2) Now suppose that Sam is female. She inherits at least one healthy X chromosome from her father, and so Sam's genotype distribution is 
        \begin{center}
            \begin{tabular}{ccc}
            XX (healthy) & Xx (carrier) & xx (colourblind) \\\hline\hline
              $\displaystyle\frac{1}{1+p}$ 
              & $\displaystyle\frac{p}{1+p}$ 
              & $0$
            \end{tabular}
        \end{center}
        (which is even better than her mother's). Sam's partner can be assumed to have the population distribution for males,
        \begin{center}
            \begin{tabular}{cc}
            XY (healthy) & xY (colourblind) \\\hline\hline
            $1-p$ & $p$
            \end{tabular}
        \end{center}

        Now consider Sam's children. For a male child, the genotype depends purely on Sam, and does not depend on Sam's partner. So 
        \[ \Pr(\text{Sam's son is xY}) = 0.5\frac{p}{1+p}  =  \frac{p}{2(1+p)} \]
        
        So Sam's son has genotype distribution
        \begin{center}
            \begin{tabular}{cc}
            XY (healthy) & xY (colourblind) \\\hline\hline
            $\displaystyle\frac{2+p}{2(1+p)}$ & $\displaystyle\frac{p}{2(1+p)}$
            \end{tabular}
        \end{center}

        For a female child, we begin by calculating the probability that she is colourblind: this can occur only if Sam's partner is colourblind, and Sam provides one copy of the defective allele:
        \[\Pr(\text{Sam's daughter is xx}) = 0.5\frac{p}{1+p}p=\frac{p^{2}}{2(1+p)}\]
        Similarly,
         \[\Pr(\text{Sam's daughter is XX}) = (1-p)(\frac{1}{1+p})+0.5\frac{p}{1+p}=\frac{(1-p)(2+p)}{2(1+p)}\]
         Hence the genotype distribution for Sam's daughter is
        \begin{center}
            \begin{tabular}{ccc}
            XX (healthy) & Xx (carrier) & xx (colourblind) \\\hline\hline
            $\displaystyle\frac{(1-p)(2+p)}{2(1+p)}$ & 
            $\displaystyle\frac{3p}{2(1+p)}$ 
            &             $\displaystyle\frac{p^{2}}{2(1+p)}$
            \end{tabular}
        \end{center}

    \end{enumerate}

    \item There are several issues that need to be considered.
    For example:
    \begin{itemize}
        \item There are other (rarer) causes of colourblindness (we have only said that colourblindness is \textit{usually} caused by a defective X chromosome, but there may well be other genetic, or environmental, causes, which will affect the probabilities calculated. These need to be factored in at the start of the problem.
        \item There may be variability by ethnicity or region. This is important, because people do not choose random partners from the entire population to have children with; they will normally choose someone who is geographically quite near to them, and they may also have a preference for choosing someone of the same ethnicity. In some regions, and/or within some ethnicities, $p$ may differ substantially, so a more targeted value of $p$ (rather than the global average) would be needed to compute meaningful probabilities. 
        \item Are `healthy' and `defective' good terms? Could they stigmatise people? This is common usage in this case, but could be more controversial in some other cases. There are various genetic conditions (e.g., autism, dwarfism) which many people deem as their own unique character traits, rather than a ``defect''.
        \item Gender and chromosomal sex are not equivalent. For example, perhaps as many as $1.7\%$ of people are intersex.\footnote{Intersex Human Rights Australia (2013). \textit{Intersex population figures}. \url{https://ihra.org.au/16601/intersex-numbers/}} This is not necessarily a small number compared to $p$. 
    \end{itemize}

\fi

\newpage
\item \textbf{Topic: \textit{Probability - Littlewood's Law} 
\\Allied disciplines: \textit{Statistics, Natural Sciences}
\\ Year level: 1
\\ \hyperref[list:pillars]{Ethical pillars}: 5}

Suppose we define a `miracle' as an exceptional event that occurs with frequency $10^{-6}$ (`one in a million'). Suppose also that a person observes $1$ `event' every second that they are alert, and that they are alert for $9$ hours of each day. Approximately how many miracles would they expect to observe in one month? How sensible or meaningful is such a computation?

\ifdefined\solution
    \solution 
	\textit{This problems demonstrates why many quantitative statements are not in a sense `mathematical' -- here for example, an `event' is not very well defined. When we teach students to model or mathematise something, we often use physical examples or examples where we know that it makes sense in a positivist or materialist sense.
  \\However, you can mathematise many things where the underlying ``object'' or ``event'' is not very well defined or ambiguous, and the mathematisation or model then often hides that fact. In this question, we ask the students do calculate the rate of ``miracles'', but you can very well use other social, psychological or theological phenomena instead. Even ambiguous physical or biological events work very well if you highlight to students the ambiguity that can be hidden behind or through mathematical formulations.
  \\Hence, problems can arise when mathematical rigor is used to give (undeserved) authority that is built on imprecise terms or unjustified assumptions.}

    The (approximate) expected rate of `miracles' will be
	\begin{align*}
	\text{rate} &= \frac{1 ~\text{miracle}}{10^{6} ~\text{seconds}} \cdot \frac{60 ~\text{seconds}}{1 ~\text{minute}} \cdot \frac{60 ~\text{minutes}}{1 ~\text{(alert) hour}} \cdot \frac{9 ~\text{alert hours}}{1 ~\text{day}} \cdot \frac{30 ~\text{days}}{1 ~\text{month}} \\
	&= \frac{60 \cdot 60 \cdot 9 \cdot 30}{10^6} \frac{\text{miracles}}{\text{month}} \\
	&= \frac{3 \cdot 3 \cdot 9 \cdot 3}{250} ~\text{miracles per month} \\
	&= \frac{243}{250}  ~\text{miracles per month} \\
	&\approx 1 ~\text{miracle per month} ~.
	\end{align*}

 For the reasons discussed above, this is not a terribly meaningful statement. One could instead describe or quantify ``events'' in a different way (eg: A person sees one event per minute, or per hour. Or an event could be something entirely different, such as ``An item presented on the 6 O'Clock news''.) Under these different definitions, the rate of miracles would differ dramatically.

\fi

\newpage
\item \textbf{Topic: \textit{Probability - Simpson's Paradox} 
\\Allied disciplines: \textit{Economics, Statistics, Engineering, Natural Sciences}
\\ Year level: 1 
\\ \hyperref[list:pillars]{Ethical pillars}: 2,4,6,9} 

	In a particular admissions cycle, a mathematics department observes a higher success rate for male applicants than for female applicants. To investigate whether this is  the same across the two sub-departments of Pure Mathematics and Applied Mathematics, the following year the department asks each applicant to give their preference for pure or applied (they are not allowed to be ambivalent), and records the resulting statistics as follows:
 
	\begin{center}
	\begin{tabular}{| c | c | c |}
			\multicolumn{3}{c}{Total:} \\
			\hline
			& Applications & Successful \Tstrut\Bstrut \\
			\hline
			Female & $300$ & $30$ \Tstrut \\
			Male & $1000$ & $210$ \Bstrut \\
			\hline
	\end{tabular}
	\end{center}
	\begin{center}
		\begin{tabular}{| c | c | c |}
			\multicolumn{3}{c}{Prefer applied:} \\
			\hline
			& Applications & Successful \Tstrut\Bstrut \\
			\hline
			Female & $270$ & $18$ \Tstrut \\
			Male & $350$ & $15$ \Bstrut \\
			\hline
		\end{tabular}
		\begin{tabular}{| c | c | c |}
			\multicolumn{3}{c}{Prefer pure:} \\
			\hline
			& Applications & Successful \Tstrut\Bstrut \\
			\hline
			Female & $30$ & $12$ \Tstrut \\
			Male & $650$ & $195$ \Bstrut \\
			\hline
		\end{tabular}
	\end{center}
	
 \begin{enumerate}
     \item  Compare the success rates for male and female applicants that prefer applied; the success rates for male and female applicants that prefer pure; and the success rates for male and female applicants over all.
 
	\item What do you notice? Why is this possible? This is known as \textit{Simpson's Paradox}.

    \item If possible, find the admission statistics by gender and mathematics preference (pure/applied) in your mathematics department, and see if the same phenomenon occurs.

 \end{enumerate}

\ifdefined\solution
    \solution 
	\textit{The purpose of this problem is to demonstrate Simpson's paradox, in which a trend appears in several different groups of data but disappears or reverses when these groups are combined. It also attempts to highlight the immense gender disparity in many mathematics departments around the world.}

 \begin{enumerate}
     \item We calculate the success rates:
	\begin{center}
	\begin{tabular}{| c | c | c | c |}
			\hline
			& Prefer applied & Prefer pure & Total \Tstrut\Bstrut\\
			\hline
			Female & $\frac{18}{270} = \frac{14}{210}$ & $\frac{12}{30} = \frac{4}{10}$ & $\frac{30}{300} = \frac{10}{100}$ \Tstrut\Bstrut\\
			Male & $\frac{15}{350} = \frac{9}{210}$ & $\frac{195}{650} = \frac{3}{10}$ & $\frac{210}{1000} = \frac{21}{100}$ \Tstrut\Bstrut\\
			\hline
	\end{tabular}
	\end{center}

 \item We note that females with a given preference have a higher success rate than males with the same preference, but lower overall. This is Simpson's Paradox. 
 \\The heuristic reason for why this is possible is that the largest male cohort (prefer pure) has a much higher acceptance rate (0.3) than the largest female cohort (prefer applied) which is about 0.067. So what dominates the overall acceptance of men is for those who prefer pure (0.3), while what dominates overall acceptance for women is those who prefer applied (0.067).
 \\This is a great lesson in why it is usually a terrible idea to ``take averages of averages''.

\item Students may have access to this data. The main point of the question is not so much to redo the calculation (it is not a given that Simpson's paradox will always arise here), but rather to illustrate the immense gender disparity in many mathematics departments around the world, which is likely (but not certain) to include the department the students are currently in. Observations like this resonate more when they are closer to home.

 \end{enumerate}

\fi    

\newpage
    \item \textbf{Topic: \textit{Probability} 
    \\Allied disciplines: \textit{Economics, Statistics, Computer Science}
    \\ Year level: 1
    \\\hyperref[list:pillars]{Ethical pillars}: 6} 
    
    Discuss the following statements, and how they use different notions of the word ``chance'':
    \begin{enumerate}
        \item `There is a $1/2$ chance that a flip of a fair coin gives heads.'
        \item `Given that I am playing as White, there is a $1/3$ chance that I will beat you at chess.'
        \item `There is a $1/2$ chance that the foetus is male.'
        \item `If you sample a random voter from a referendum where $52\%$ of voters voted `yes', then there is a  $52\%$ chance that they voted `yes'. '
        \item `If you sample a random voter before an upcoming election, there is a $30\%$ chance that they will vote for the incumbent.'
        \item `There is currently a $40\%$ chance that AI will wipe out humanity in the next 5 years.'
    \end{enumerate}
    Which of these events would you bet on (either for or against), if any? 
    
    \ifdefined\solution
    \solution \textit{These statements are meant to highlight the different ways in which a probability (ie: chance) can be interpreted by mathematicians, or by the population at large. }
    
    \begin{enumerate}
        \item This is an example of \emph{aleatoric uncertainty}. The probability can be determined by conducting many trials and observing that $1/2$ of them produce successes. 
        \item Chess has no chance elements by itself (apart from the process to determine who plays White, hence the conditioning). The probability expresses uncertainty about the relative skills of the two players, as well as uncertainty about the factors (\textit{e.g.} psychological) that might affect the quality of their play. 
        \item This is an example of \emph{epistemic uncertainty}: the foetus' sex has already been determined. It is presently unknown, but could be measured. 
        \item Each voter's choice in the referendum has already been fixed, but the randomisation comes from the sampling process. 
        \item In contrast, voters have not yet fixed their individual choices for the next election yet. 
        \item This is a statement about an uncertain event that cannot be repeated and has no precedents. The $40\%$ probability is a subjective assessment.
    \end{enumerate}
    As these examples demonstrate, probabilities should sometimes be interpreted as subjective quantities (`degrees of rational belief'), which are to be updated as further information becomes available (using Bayes' formula).  
    
    The interpretation of randomness sometimes has legal consequences. In some jurisdictions, gambling laws apply to wagers on games of chance but not on games of skill, and the division can be somewhat arbitrary: for example, poker is considered to be a game of skill in New York, but a game of chance in Germany. 
    \fi

\end{enumerate}

\newpage
\section{Exercises: Statistics}

\begin{enumerate}
\item \textbf{Topic: \textit{Statistics}
\\Allied disciplines: \textit{Economics, Probability, Computer Science}
\\Year level: 1
    \\\hyperref[list:pillars]{Ethical pillars}: 3,4,5,6}
    
    Your city council aims to redesign the city centre to make it more bike-friendly if at least 75\% of the citizens agree. The city administration claims that this is the case. A competing political initiative believes that the actual percentage of supporters is less than 75\%. To test the claims, a hypothesis test is conducted at the 5\% significance level with a sample size of $n=100$.
\begin{enumerate}
\item How does the city council test? For which results is it confirmed?
\item How does the competing political initiative test? For which results is it confirmed?
\item For which results can neither party reject the null hypothesis? 
\item Assume that the null hypothesis could not be rejected. Without further comment, one of the parties approaches you with another 100 data points and asks you to perform the test again by combining the old and new samples to reach a sample size of $n=200$. What do you do?
\end{enumerate}
\ifdefined\solution
    \solution
\textit{The purpose of this question is to teach students some basics about good data sampling. It also teaches students to explore competing ethical interests in data collection and analysis, and the pressures that a mathematician can face in scenarios like these.}

\begin{itemize}
    \item \textit{The student should understand the difference between a right- and left-sided hypothesis test and be able to explain why the parties chose different tests in questions (a) and (b).}
    \item \textit{The student should see the problem of p-hacking and lay out some of its properties in this context. They should be able to discuss the dangers of using an existing sample multiple for multiple decisions: to perform the hypothesis test, to use the result to decide that a larger sample is needed, and then to use the data points again in a later test on combined data.
    The student should also be able to question the source of the new sample and if it was properly sampled. They should understand that competing ethical interests could mean that the standards of good sampling/interviewing may have been violated.}
\end{itemize}

Due to time constraints, we were unable to include a full solution to this exercise.

\fi

\newpage
\item \textbf{Topic: \textit{Statistics}
\\Allied disciplines: \textit{Probability, Natural Sciences}
\\Year level: 1
    \\\hyperref[list:pillars]{Ethical pillars}: 6,9,10} 
    
In  2020 the UK government tried to explain their new COVID alert levels by the following equation:
$$\text{COVID Alert Level = R (rate of infection) + number of infections}.$$
where the alert level is a discrete scale from $1-5$, and $R$ is the average number of people each infected person infects.

Why does this not make much sense? Why could it still be good to communicate this way? Why could it be dangerous to communicate this way? How would you communicate it?

\ifdefined\solution
\solution 
\textit{This question teaches students that it is often necessary to break down difficult mathematics into everyday language. It teaches them about what information is lost along the way and to be wary of the fact that others will have a lot less mathematical training (and hence may need a simpler form of communication).
You can find a full discussion in the article by Woolley \footnote{T. Woolley (2020). Coronavirus: Why the maths behind ‘COVID alert levels’ makes no sense. \textit{The Conversation}. \newline  \href{https://theconversation.com/coronavirus-why-the-maths-behind-covid-alert-levels-makes-no-sense-138634}{https://theconversation.com/coronavirus-why-the-maths-behind-covid-alert-levels-makes-no-sense-138634}}.}
\begin{enumerate}
    \item The units don't match: The alert level is unitless and ranges from 1-5, R is usually also unitless while the rate of infection is not, the number of infections is unitless. All three parameters are also defined on different ranges (R is usually quite small, while the number of infections can be in the thousands).
    \\Students might observe that actually there might have been some predetermined function $f: \mathbb{R}\times \mathbb{R}\to \mathbb{N}$ giving
    \[
    \text{ COVID Alert Level = f(R,number of infections)}
    \]
    ant that this is actually a sensible statement and probably close to what was being done at the time. They should then consider how much the general population would understand this mathematically formalised statement.

    \item The benefit of this style of communication is that people quickly learned to worry about both the ``rate of infection'' and ``the number of infections''. It breaks down difficult concepts into everyday language.
    \item People can be confused as it doesn't properly weigh the parameters. They want to tell you that both the rate and number of infections are important, but are they equally important? Is one or the other more important? Is the importance constant over time or does it vary?
\item This is hard, and that is the first thing students should acknowledge. Not everyone has a degree in mathematics, so one is constrained in the language and terminology that can be used.
\\They might lean towards variants of ``The COVID alert level \textit{depends} on the infection rate \textit{as well as} the number of infections''. But how does it depend on these? How are they weighed? And how can one communicate this to public in a quick and effective manner, that fits on one static PowerPoint slide? What is the actual point of this communication; what is it trying to achieve? (it is not there just as an education segment for the population; it is there to effect some change in society - what change?
\\Ultimately, who is harmed if this communication is not effective, and how badly are they harmed? And who will take the blame, or be blamed, if the communication is ineffective and there is resulting harm?

\end{enumerate}

\end{enumerate}
\newpage
\section{Exercises: Numerical Methods}
\begin{enumerate}
    \item \textbf{Topic: \textit{Numerical Methods} 
    \\Allied disciplines: \textit{Economics, Statistics, Computer Science, Engineering, Natural Sciences}
    \\ Year level: 1 
    \\\hyperref[list:pillars]{Ethical pillars}: 6} 

    \begin{enumerate}
        \item     Let $F:[a,b]\rightarrow\mathbb{R}$ be a smooth function with $F(a) < 0 < F(b)$, and $F'(x) > 0$. Then the equation $F(x) = 0$ has a unique solution $x=c$. Numerical approximations to $c$ can be found using Newton--Raphson iteration, which converges quadratically.

State the formula for Newton--Raphson iteration and explain it graphically. Convey the above information to someone who has no post high school mathematical training. 

\item State the formulae for the Finite Difference methods (forward, backward and central) for the approximation of the first derivative of a function $F(x)$ at point $a$ denoted by $F'(a)$ using a step size $h$, and explain it graphically. Convey this information to someone who has no post high school mathematical training.

\end{enumerate}

\ifdefined\solution
    \solution \textit{The question is designed for students to break down the mathematics into non-technical language and to communicate in a way that requires minimal mathematical knowledge. It teaches them that pictures and graphic representations are often a useful tool to explain aspects to non-technical people. For both parts of the question, it would be useful to discuss with students which words/terminology are \textbf{not} appropriate to use in explanations to a lay audience. For example, is it reasonable to use words such as tangent, derivative, differentiable, continuous, domain, root, convergence? The audience may not have studied mathematics beyond the age of 16, and even if they did to mathematics through to the end of high school, that may have been decades ago for them.}

\begin{enumerate}
    \item The iterative formula for Newton--Raphson iteration, with initial guess $a<x_{0}<b$ of the root of the function $F$, is given by:
    \[
    x_{n+1}:= x_{n}-\frac{F(x_{n})}{F'(x_{n})}
    \]
    To see this graphically, draw a hypothetical graph of $F$, and take a point $x_{n}$. Then draw the tangent to $F$ at $x_{n}$, and see where this tangent intersects the x-axis (it will definitely have some non-zero slope, as $F'(x)>0$ on $[a,b]$, and call this new point $x_{n+1}$.
    \\The iteration continues in this way; for each new point computed, draw the tangent to $F$ at that point, see where the tangent intersects the $x$-axis, and take that as the next new point.

    To explain all this to someone with no post high school mathematics training, the following terms might be used:
    \begin{itemize}
        \item \textbf{$\mathbf{F:[a,b]\rightarrow\mathbb{R}}$:} $F$ is a graph in the $xy$-plane that starts at $x=a$ and ends at $x=b$.
        \item \textbf{Smooth function:} the curve/graph of $F$ has no breaks in it, and has no sharp points, spikes, corners or kinks. 
        \item \textbf{$\mathbf{F(a) < 0 < F(b)}$:} the starting point (at $a$) of the curve of $F$ is less than 0, and the endpoint is greater than 0.
        \item \textbf{$\mathbf{F'(x) > 0}$:} the slope of the curve of $F$ is always ``uphill'', when travelling along it from left to right. It is never flat or downhill.
        \item \textbf{The equation $\mathbf{F(x) = 0}$ has a unique solution $\mathbf{x=c}$:}
        if all the conditions on $F$ listed above hold, then $F$ crosses/intersects the $x$-axis exactly once (between $a$ and $b$), at a value which we will call $x=c$.
        \item \textbf{Numerical approximations to $\mathbf{c}$ can be found using Newton-Raphson iteration:} without having an exact way to figure out where $F$ crosses the $x$-axis, we can nonetheless approximate $c$ by taking successively better and better estimates. We start with a ``guess" for $c$, and then create a better estimate using some straightforward calculations. Each time we make an estimate, we use that to then make an even better estimate (i.e.: even closer to $c$).
        \item \textbf{Converges quadratically:} The more of these estimates we calculate, the  closer and closer they get to c. And this happens rather quickly; we don't need to make too many of these estimates to get very close to the true value of $c$.
        \item \textbf{State the formula for Newton--Raphson iteration and explain it graphically} This ``estimation" method works as follow. First, we pick a random point somewhere between $a$ and $b$; the midpoint will do. We call this the `First Estimate'. Then we look at the curve $F$ at this value `First Estimate', and on the curve at that point we draw a straight line which is tangential to the curve; that is, it just touches $F$ at that point, like the edge of a long pencil touching the side of a basketball. We know that the curve $F$ is always "uphill", so this line will have some slope on it (it won't be horizontal), and will thus cross the $x$-axis somewhere. We look at where this tangential line crosses the $x$-axis, and this becomes our next estimate, which we will call `Second Estimate'. We then run this whole process again, starting with the value `Second Estimate', to generate an even better estimate which we call `Third Estimate'. And so on. This can all be done on a computer, and with perhaps only a few hundred successive estimates (which would take less than 1 second on a normal computer) we would have an estimate that is extremely close to the true value of $c$; perhaps within a few thousandths of a percent. For most industrial applications, an approximation like this is perfectly acceptable, and overall this is a much better approach than doing a lot of complicated mathematics to figure out the exact value of $c$.    
    \end{itemize}
The stress here is that this needs to also be \textit{presented} well, i.e., not in bullet-point form as above, but in regular prose (paragraphs).

\item The formulae for each of the three finite difference methods, to approximate the derivative $F'(a)$ for a function $F(x)$ at point $a$, using a step size $h$, are given by:
    \[
    \text{Forward difference:   } F'(a) \approx \frac{F(a+h) - F(a)}{h}
    \]
    \[
    \text{Backward difference:   } F'(a) \approx \frac{F(a) - F(a-h)}{h}
    \]
     \[
    \text{Central difference:   } F'(a) \approx \frac{F(a+h) - F(a-h)}{2h}
    \] \\
    To explain this graphically, draw a hypothetical graph of $F$, and take a point $a$. Then draw the tangent to $F$ at $a$ (a line just touching $F$ at $a$ but not crossing or touching anywhere else on $F$), which represents the slope of the graph $F$ at point $a$, which we can denote as $F'(a)$. \\

    \textbf{The forward finite difference method essentially tries to approximate $F'(a)$ as follows:} consider another point $a+h$ on the graph of $F$, a short distance $h$ to the right of point $a$ (assuming $h$ is positive) and draw a line joining this point $a+h$ on $F$ to the previous point $a$. This yields a straight line that gives us a closer enough approximation to the tangent line drawn earlier on $F$ at point $a$ (as the two do not coincide/are not parallel and therefore do not have the same slope). The value of this approximation to the slope $F'(a)$ is given by the the value of the slope of the line connecting the two points $a$ and $a+h$ on the graph of $F$, which is indeed given by the difference between the function values $F(a+h)$ and $F(a)$ divided by the small distance $h$ between them, where $h$ is called the `step size', which is exactly what the mathematical expression of the formula states. \\
    
    \textbf{The backward finite difference method essentially tries to approximate $F'(a)$ as follows:} consider another point $a-h$ on the graph of $F$, a short distance $h$ to the left of point $a$ (assuming $h$ is positive) and draw a line joining this point $a-h$ on $F$ to the previous point $a$. This yields a straight line that gives us a closer enough approximation to the tangent line drawn earlier on $F$ at point $a$ (as the two do not coincide/are not parallel and therefore do not have the same slope). The value of this approximation to the slope $F'(a)$ is given by the the value of the slope of the line connecting the two points $a$ and $a-h$ on the graph of $F$, which is indeed given by the difference between the function values $F(a)$ and $F(a-h)$ divided by the small distance $h$ between them, where $h$ is called the `step size', which is exactly what the mathematical expression of the formula states. \\

    \textbf{The central finite difference method essentially tries to approximate $F'(a)$ as follows:} consider two additional points $a-h$ and $a+h$ on the graph of $F$, a short distance $h$ to the left and right of point $a$ respectively (assuming $h$ is positive) and draw a line joining these two points $a-h$ and $a+h$ on $F$. This yields a straight line that gives us a closer enough approximation to the tangent line drawn earlier on $F$ at point $a$ (as the two do not coincide/are not parallel and therefore do not have the same slope). The value of this approximation to the slope $F'(a)$ is given by the the value of the slope of the line connecting the two points $a+h$ and $a-h$ on the graph of $F$, which is indeed given by the difference between the function values $F(a+h)$ and $F(a-h)$ divided by the small distance $2h$ between them, where $h$ is called the `step size', which is exactly what the mathematical expression of the formula states. \\

    A couple of points worth noting when explaining: 
    \begin{itemize}
        \item The central difference method approximation is more accurate than both the forward and backward difference methods as it is essentially equivalent to the average value of each of their individual approximations (this can be checked and verified mathematically and is a good exercise to try - it can also be confirmed graphically as it yields a straight line that more closely mirrors the slope of $F$ at point $a$ given by the tangent line). \\

        \item The smaller the step size value $h$, the greater the accuracy of the approximation to the slope of $F$ at $a$, as the additional points $a+h$ or $a-h$, will become much closer to $a$ and therefore yield straight lines that more closely correspond to the actual slope desired (this can be verified graphically using a graphical plotting software like Desmos by constructing a simple curve $F$, the tangent line at $a$ and the straight lines mentioned above in each case and then varying the parameter $h$ to see that when $h$ approaches $0$, the straight lines exactly coincide with the tangent, thereby making the approximation exact). \\

        \item In this case, we use the word `slope' to commonly refer to the first derivative of the function $F$ to make it easier to explain the concept to someone with no post high-school background in mathematics. We are concerned here only with first order derivative for the appromximations, but it is worth noting that a central finite difference approximation can also be derived and does exist for second order derivatives as well (this is also an exercise worth trying). 
    
    \end{itemize}

\end{enumerate}

\end{enumerate}
\newpage
\section{Exercises:  Pre-Calculus and Functions (catch-up / introductory)}
\begin{enumerate}
    \item \textbf{Topic: Pre-Calculus 
    \\Allied disciplines: \textit{Economics, Statistics, Computer Science, Engineering, Natural Sciences}
    \\ Year level: 0-1 
    \\\hyperref[list:pillars]{Ethical pillars}: 5,6}
    \\As you have seen during school, text problems are often highly dependent on their historical and social context. But this context naturally changes throughout time or when you go somewhere else. This is one of the reasons why textbooks are regularly updated to not just reflect the newest research in pedagogy, but also to ensure that their ethics goes with the time. 

    Read carefully, but don't solve, the following question from a German textbook during World War II. You can find this and other examples in Shulman\footnote{B. Shulman (2002). Is There Enough Poison Gas to Kill the City?: The Teaching of Ethics in Mathematics Classes. \textit{The College Mathematics Journal} 33(2), pp. 118-125. \href{https://www.tandfonline.com/doi/abs/10.1080/07468342.2002          .11921929           }{https://www.tandfonline.com/doi/abs/10.1080/07468342.2002.11921929                             }}:
    
    \setlength{\leftskip}{1.5cm}
    \setlength{\rightskip}{1.5cm}
    ``According to statements of the Draeger Works in Luebeck, when using pesticides on a field, only 50\% of the chemicals are effective. The atmosphere must be covered up to a height of 2 metres in a concentration of 45 mg per cubic metre to cover all the crops. How much chemical is needed to cover crops of four square kilometres?''

    \setlength{\leftskip}{0cm}
    \setlength{\rightskip}{0cm}
    
    Now solve the following version of this question:

    \setlength{\leftskip}{1.5cm}
    \setlength{\rightskip}{1.5cm}
    ``According to statements of the Draeger Works in Luebeck, when gassing a city, only 50\% of the evaporated poison gas is effective. The atmosphere must be poisoned up to a height of 20 metres in a concentration of 45 mg per cubic metre. How much phosgene is needed to poison a city of 50,000 inhabitants who live in an area of four square kilometres?''

    \setlength{\leftskip}{0cm}
    \setlength{\rightskip}{0cm}
    
    The mathematics you've used is almost identical. We call this the dual-use problem of mathematics: many mathematical equations and results can be simultaneously used for beneficial and harmful applications. Can you find yet another example for the same mathematics?
    
    \ifdefined\solution
    \solution
    \textit{A word of warning is necessary for this question. The dual use of mathematics is usually not explored at the high school level. When using this question, be aware that students have very little background and have almost never seen a mathematics question that actively asks them to do something bad. This is why we urge you to not have students solve the original problem directly, but rather to explore the mathematics and ethics of this question indirectly by solving a similar problem and finding other applications for the same mathematics.
    \\Before setting this question, you should also be aware of who your students are, their background and how they'd react to seeing mathematics in this context. Care and empathy are needed when going into the very dark areas of mathematical applications and education. To not leave students alone with their emotions, we suggest that this question is best suited for classroom discussions rather than as a homework problem.}
    \fi
\end{enumerate}

Due to time constraints, we were unable to include a full solution to this exercise.

\newpage
\section{Projects: Communicating your Mathematics}
\begin{enumerate}
    \item \textbf{Communicating and Documenting mathematics in a team - The Ethics Memo }

Additional question: Communicate a problem at hand to a relevant party (e.g. your boss) in an \textit{Ethics Memo}, as illustrated in the outline below.

This question can be added on to an existing question or exercise from previous sections. It teaches students a way to communicate and document their ethical concerns.

\textit{\textbf{Notes to Instructor:} Ethics memos are an experimental idea which we've internally tried and discussed. It is best suited for problems that have multiple feasible solutions (e.g. optimisation problems that allow for multiple ways to model it or when a non-mathematical solution is a good option) or face competing interests (e.g. as in the case of the oil pipeline: environmental constraints vs business interests). The style of this memo is adapted from policy memos which are used to analyse political problems and make recommendations.}

\textit{To assist students in seeing what aspect of ethical issues might be worth signposting in such a memo, it might be useful to suggest that they consult the Manifesto for the Responsible Development of Mathematical Works\footnote{M. Chiodo \& D. M\"uller (2023). Manifesto for the Responsible Development of Mathematical Works -- A Tool for Practitioners and Management. \href{https://arxiv.org/abs/2306.09131}{https://arxiv.org/abs/2306.09131}}}. 
\end{enumerate}
\newpage

\begin{mdframed}
\textbf{The Ethics Memo}
\newline \mbox{}
\newline \textbf{To:}
\newline \textbf{From:}
\newline \textbf{Re:}
\newline \textbf{Date:}

\mbox{}
\newline \textbf{Issue:} Summarise the ethical challenge that you have encountered during your mathematical work. This includes translating professional jargon into readable prose. The level of detail should be adopted to how familiar the recipient of this memo is with the problem and the mathematics. Generally speaking: The shorter the better, often a few lines is enough. 
\newline \newline \textbf{Interests:} In this section you should address why the recipient ought to care about the issue. It should address the most important ethical and other interests. You should rank all interests according to a pre-defined standard (e.g. necessary legal obligations, vital ethics/norms/interests, and other important interests that are not vital). 
\newline \emph{Legal:}
\newline \emph{Vital:}
\newline \emph{Important:}
\newline \newline \textbf{Analysis:} This section analyses the general shape of the ethical challenge. You should briefly explain its dynamics and point out key causal factors. Try to explain some of the necessary background in a way that the recipient of this memo can understand and wants to read.
\newline \newline \textbf{Objectives:} Consider what resources and capabilities can be mobilized, and formulate achievable goals: what can you feasibly do to address the ethical issue? Do not aim for too much. If the problem is too involved, make this clear. Be humble and aware of how much you can actually do.
\newline \newline \textbf{Options:} Identify and briefly describe distinct (non-)mathematical solutions that address the issue. This section must include an assessment of the pros and cons of each solution. Only consider realistic and feasible solutions. Avoid writing sandwich memos where you present two unreasonable options and the only reasonable recommendation is the one you prefer. Whenever possible, you should present strong alternatives. Understand that \textit{not} doing mathematics, or doing \textit{less} mathematics, is often an option. Not every problem requires a mathematical solution.
\newline \newline  \textbf{Recommendation:} This section is a concise summary of the (unique) option you recommend and why you recommend it. Try to get your point across as succinctly as possible. Be aware of any preconceptions that your recipient might have about your recommendation.
\newline \newline \textbf{Implementation:} Summarise the initial steps that need to be implemented. Whenever possible, avoid professional jargon while accurately describing each step. These steps can (and often do) include the consultation of other specialists and domain experts. Briefly describe how success could be measured.
\end{mdframed}

\newpage

\noindent Here is an example of a `well written' ethics memo, inspired by the Lake Contamination Problem (see Question \ref{chem-spill} in Section \ref{diff-Qs}). Note that the problem and this ethics memo are entirely fictional and do not constitute legal/moral advice for a real situation. 

\begin{mdframed}
\textbf{Ethics Memo. Project \#561 - Lake [Name] contamination.}
\newline \mbox{}
\newline \textbf{To:} Mrs Jane Anderson, senior manager of safety and compliance division. \url{j.anderson@[domain_name]}
\newline \textbf{From:} Mr Richard Smith. Engineering division 2. \url{r.smith@[domain_name]}
\newline \textbf{Re:} Calculating when  the Lake [Name] water will be suitable for drinking.
\newline \textbf{Date:} 12 March 2024.
\newline \newline \textbf{Issue:} My team was tasked with calculating when the water in Lake [Name] might be safe to release into the drinking supply again, following a chemical contamination event in late 2023. I am concerned that there might be an over-reliance on these calculations in determining when to re-introduce water from Lake [Name] into the drinking supply, and that if taken as the only factor in decision-making, could lead to unsafe contamination of the drinking supply. 
\newline \newline \textbf{Interests:} Our organisation is responsible for making the final decision on whether, and if so, when, to reintroduce this water back into the water supply. 

\textit{Legal interests:} My understanding is that our organisation is fully liable for any avoidable damages, harms, or injuries that come from such a decision, as well as costs to rectify these. The directors would also be criminally liable for any harms to persons.

\textit{Vital interests:} Lake [Name] is part of the supply to the lower-east region, affecting 7 million households, and 2 million businesses and state organisations (hospitals, schools, etc). Contamination at this scale would pose significant disruption, and could lead to business closures, general panic, mass hospitalisations, and loss of life.

\textit{Important interests:} If contaminated water is introduced into the drinking supply, even temporarily, this would take vast time and money to clean, further delaying a return to normal supply. The reputational, and financial, damage to our organisation would  be immense.
\newline \newline \textbf{Analysis:} The contaminant in question, \textit{benzene}, is known to cause serious health conditions if ingested in sufficient quantities. These include drowsiness, delerium, loss of consciousness, and in some cases, death\footnote{\textit{Benzene: toxicological overview} https://www.gov.uk/government/publications/benzene-general-information-incident-management-and-toxicology/benzene-toxicological-overview}. At the concentrations currently present in Lake [Name], most adults would feel negative effects and illness after consuming half a glass of water, and require hospitalisation after  5 glasses. 

Contamination was caused by a leaking pipe in a nearby petrochemical factory which began on 27 November 2024. This was secured on 16 December 2024. Over time, the concentration of benzene in the lake will decrease as clean rainwater continues to flow in, and contaminated water flows out into the ocean. My team was tasked with calculating when the lake water will be safe enough for drinking again.

We have calculated a `safe date' of 3 April 2024. However, our mathematical prediction makes many assumptions and estimates, several of which are difficult to confirm precisely. For example, a decrease in rainfall might decrease the inflow of clean water, whereas heavy rainfall would increase the inflow of clean water, but might also dislodge stagnant benzene between the chemical plant and the lake. Our calculation gives an ``average case'' prediction as to when the water will be safe to drink again, not an exact prediction.
\newline \newline \textbf{Objectives:} The uncertainty around our prediction must be taken into account when using it for decision-making. This can still be very useful for arranging interim measures, such as the contracting of other water sources for public supply (other lakes, bottled water, etc). But for the core, high-risk decision of reintroducing the lake water into the main supply, we should investigate obtaining additional updated data regarding the state of the lake and its surrounds.
\newline \newline \textbf{Options:} There are several solutions available at this point. Independently of these options, we need to clean the immediate area surrounding the leak and its pathway to the lake, so as to prevent recontamination in the future.

\textit{Option 1.} We can proceed with using this average-case prediction, 3 April 2024, to decide when to re-introduce the lake water, and how long to implement interim measures. While this is likely to be the most financially efficient approach in the short term, we have no way of definitely knowing if the water will be safe on 3 April 2024. If our estimate is longer than necessary, additional interim costs will have been incurred unnecessarily, as well as unnecessary disruption to the public. If our estimate is shorter than necessary, there is a chance of significant harm to the public, significant additional delays, and significant additional costs to us.

\textit{Option 2.} We can use the calculation of 3 April 2024 as an earliest-case assumption, and proceed with planning interim measures at least up until that date. Then, shortly before our predicted date, we can re-test the lake water, and use that data to update our prediction, extend interim measures if needed, and then repeat the process. Only when the lake passes full testing should the water be re-introduced into the supply, fulfilling our legal duties to ensure safe supply. This reduces the chance of significant harm to the public and significant costs to us, but may still delay reintroduction unnecessarily, inconveniencing the public, harming businesses, and increasing our costs.

\textit{Option 3.} We can repeatedly re-test the lake water, at daily intervals. As soon as the water is determined to be clean enough, we can re-introduce it into the supply. This would minimise the disruption to the public, but would increase our costs in terms of testing, and also in terms of having to repeatedly re-negotiate interim measures, such as alternate water sources, at short notice. However, our reputation would be best-preserved under this option.
\newline \newline  \textbf{Recommendation:} I recommend implementing option 2, alongside a  cleanup of the area surrounding the lake. It minimises significant risks to the public and to our organisation, and likely also meets our legal obligations. The additional costs of delaying the reintroduction unnecessarily are manageable compared to reintroducing the water too early.
\newline \newline \textbf{Implementation:} 
\\- Engage water testing professionals to test the water.
\\- Engage a clean-up team to carry out surrounding ground decontamination.
\\- Communicate test/measurement results to my team for additional calculations.
\\- Communicate timelines to the team arranging interim water sources.
\\- Communicate minimal, and probable, projected timelines to the public.
\\Our success here can be measured by seeing that we safely reintegrate the lake into the water supply in a timely manner.
\end{mdframed}

\newpage

\noindent Here is an example of a `poorly written' ethics memo, inspired by the Oil Pipe Optimisation Problem (see Question \ref{oil-pipe} in Section \ref{optimisation}. Note that the problem and this ethics memo are entirely fictional and do not constitute legal/moral advice for a real situation.
\begin{mdframed}

\textbf{Memo} \\
\newline \textbf{To:} Steve
\newline \textbf{From:} Alex
\newline \textbf{Re:} Caclulations
\newline \textbf{Date:} 2/9 \\
\newline \textbf{Issue:} The main thing here is to figure out how to connect the oil platform to the refinery for the cheapest possible cost. Using the equation $C(x) = c_2 x \, + \, c_1 \sqrt{D_1^2 + (D_2-x)^2}$, we know there's a way to get this done efficiently. Basically, we just need to solve for $x$ where $\frac{dC}{dx} = 0$. Not much else to say - it's basic maths. \\
\newline \textbf{Interests:} 
\newline \emph{Legal:} Uh, make sure we don't break laws, I guess.
\newline \emph{Vital:} Cheap pipeline construction, obviously. That's all.
\newline \emph{Important:} Finish the project quickly and on time. \\
\newline \textbf{Analysis:} This is trivial. By solving $\frac{c_2}{c_1} = \frac{(D_2 -x)}{\sqrt{D_1^2 + (D_2-x)^2}}$, we can find the distance $x$ of the pipeline to minimise construction costs. For example, if $D_1 = 10, D_2 = 20, c_1 = 5$ and $c_2 = 2$, then the solution gives us $ x \approx 15.64$. Plus that into the original cost function and we're golden. Not much else to it. \\
\newline \textbf{Objectives:} 
We want to finish the pipeline project for the least amount of money. And have no major delays, if possible. 
\newline \newline \textbf{Options:} 
\newline Option 1: Cheapest Path
\begin{itemize}
  \item Pros: Super cheap, which is great! 
  \item Cons: Might not be perfect.
\end{itemize} 
Option 2: Slightly more expensive path
\begin{itemize}
  \item Pros: Less complicated maths. Might avoid other problems.
  \item Cons: It's more expensive, so seems like a bad option.
\end{itemize} 
\textbf{Recommendation:} I say we just go with Option 1 and put the pipeline wherever is mathematically cheapest. Just do what the maths says. There’s no point in wasting more time or money looking into extra stuff. \\
\newline \textbf{Implementation:} Solve $\frac{c_2}{c_1} = \frac{(D_2 -x)}{\sqrt{D_1^2 + (D_2-x)^2}}$ and get the exact $x$. Then build the pipeline. Done. Success is measured by completion of the pipeline in the cheapest possible way.
\end{mdframed}

\newpage \noindent\textit{\textbf{Note to Instructor:} This last memo deliberately ignores critical ethical, legal, and environmental issues, provides minimal context, uses overly technical jargon, and lacks professionalism. In more detail, here is a list of shortcoming of the memo, giving roughly in order that they appear. This could be turned into a student exercise itself, by giving the students the `bad memo' and asking them to identify as many problems with it as they can find.}
\textit{\begin{itemize}
    \item The memo is not given any sort of description; what type of memo is it?
    \item Surnames are missing, as are basic descriptions of the parties, and contact details.
    \item Topic is completely lacking in detail; ``calculations'' could refer to any number of projects.
    \item Date is not given in full, and is ambiguous as to whether it is day/month or month/day format.
    \item The issue doesn't really specify an ethical challenge. It is packed with technical formulae, and leaves a lot for the reader to infer themselves. 
    \item Legal interests are completely vacuous. Vital and Important interests give minimal details.
    \item There is no mention of any outside interests.
    \item There no analysis of any ethical problem; just a description of the mathematics being used, and even that is presented in a very technical manner which would be inaccessible to many decision makers.
    \item The objectives just re-state the mathematical problem in simplistic terms, with no acknowledgement of any additional factors that may need to be considered.
    \item The options are presented in a loaded manner.
    \item The options lack detail, and any deep reflection on impacts to either the organisation or external parties.
    \item All that is ever considered is immediate cost to build, which is a superficial optimisation setup.
    \item The recommendation explicitly dismisses other potential concerns, and relies only on crude mathematical considerations.
    \item The implementation gives only details for what to do mathematically, and no details about the physical act of building the pipeline.
\end{itemize}
}

\fi

\newpage 
\begin{enumerate}[start=2]
    \item \textbf{Group project: How your science can end up in the newspapers.} \\
When you read about scientific research in the media, it has gone through several steps before the information has reached you. Consider David Spiegelhalter's \footnote{D. Spiegelhalter (2017). Trust in numbers. \textit{Journal of the Royal Statistical Society Series A: Statistics in Society}, 180(4), pp. 948-965. \url{https://academic.oup.com/jrsssa/article/180/4/948/7068387}} pipeline for how statistical results end up in the news:
\begin{figure}[H]
  \includegraphics[width=\linewidth]{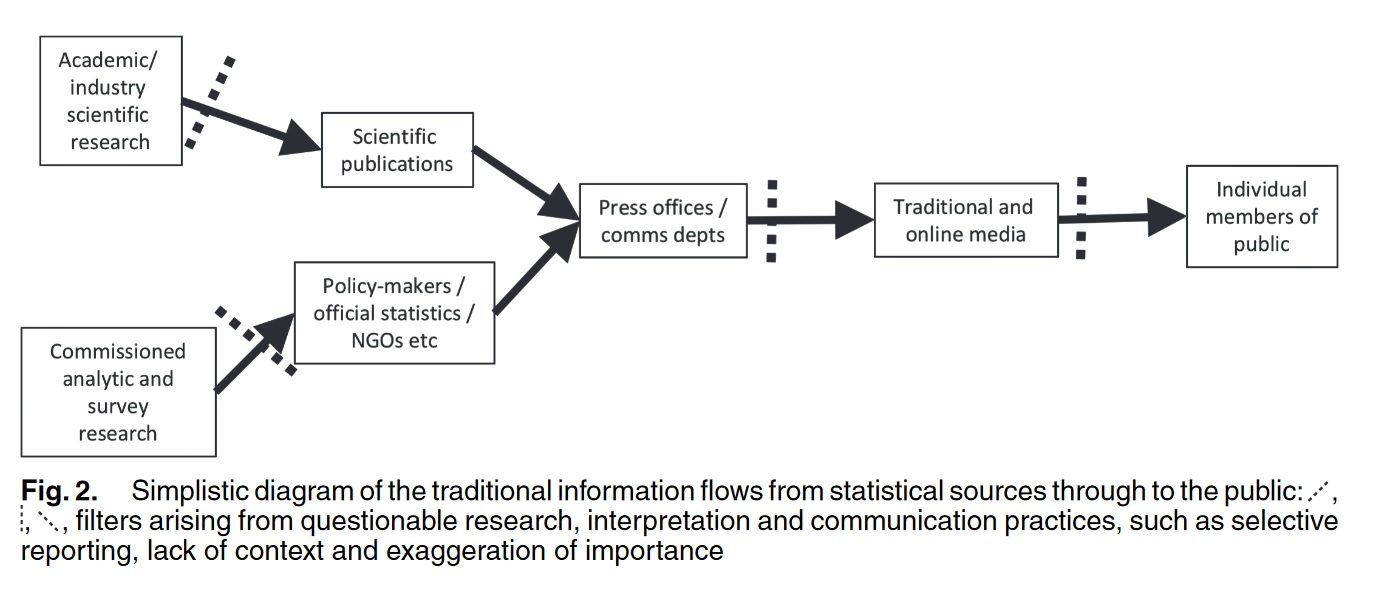}
  \caption{Pipeline of Science Communication (Spiegelhalter, 2017).}
  \label{fig:news}
\end{figure}

Step 1: Understanding the diagram
\begin{enumerate}
    \item Discuss the diagram in your group. 
    \item Try to understand the work, interests and standards of each entity in the diagram. Do you see any conflicts of interest arising between the different entities?
    \item Do you know what all the filters are and how they work? 
\end{enumerate}

Step 2: Analysing a case study
\begin{enumerate}
    \item Read the scientific paper [INSERT TITLE] and summarise its main results.
    \item Find and read newspaper articles about this research. Summarise how they presented the results to the public.
    \item Understand if the presentation by the news is accurate. If possible, try to identify where the problem started.
\end{enumerate}
This question is supposed to teach students that they might not always have control over how other people perceive their published work.

\textit{\textbf{Notes to Instructor: } It is best to select a research paper such that its (mathematical) results are not just suitable for your course, but for which students can find articles for the entire pipeline, including communication from the university's press office and a range of newspaper articles. This enables them to properly understand how research is publicised and broken down to fit different audiences.}

\end{enumerate}
\newpage
\section{Projects: Open-ended discussion questions}

The following questions, based on particular case studies, require a little more thought, discussion and some research. They are well suited for seminar-style discussions.

\begin{enumerate}
    \item This question is about the Diffie--Hellman key exchange protocol. 
    \begin{enumerate}
        \item Explain to a non-mathematical audience why the Diffie--Hellman algorithm is important.
    \end{enumerate}
    
        In 2015, a major flaw was discovered in the implementation of the Diffie--Hellman algorithm by Adrian et al.\footnote{D. Adrian et al. (2015). Imperfect Forward Secrecy: How Diffie--Hellman Fails in Practice. \textit{CCS '15: Proceedings of the 22nd ACM SIGSAC Conference on Computer and Communications Security}. \href{https://dl.acm.org/doi/10.1145/2810103.2813707   }{https://dl.acm.org/doi/10.1145/2810103.2813707  }}, which affected tens of thousands of popular websites.
        
        Read Adrian \textit{et al.} (2015) and pay particular attention to \S\S5--7. 
    \begin{enumerate}
        \setcounter{enumii}{1}
        \item What were the reasons for the flaw?
        \item Do you think the authors were right to publish this work?
    \end{enumerate}

\ifdefined\solution
    \solution This question deals with two main ideas:
\\1. What are the ethical issues associated with creating large-scale tools to break encryption en masse, especially when the tool itself is prohibitively expensive to produce (large lookup table), but comparatively easy to use and to lose/leak/steal (how ``safe'' can one keep such a tool)?
\\2. Regarding responsible disclosure, what are the ethical issues with exposing state surveillance capabilities, and does the fact that a state enabled the vulnerability in the first place justify full disclosure (which might incentivise other states to develop, or just outright steal, this technology)?

    \fi

    \newpage
    
    \item The statistician Ronald Fisher (1890--1962) made many important contributions to statistics, including the concepts of hypothesis testing and $p$-values. However, Fisher's work on statistics was partially motivated by a desire to justify his beliefs in parapsychology and racialism. This motivation does not necessarily invalidate Fisher's mathematical work. But how can the scientific community safeguard against misuse of statistical techniques? Start by reading about Fisher, his mathematics and views on eugenics and race.\footnote{W. Bodmer et al. (2021). The outstanding scientist, R.A. Fisher: his views on eugenics and race. \textit{Heredity} 126, pp. 565–576. \href{https://www.nature.com/articles/s41437-020-00394-6}{https://www.nature.com/articles/s41437-020-00394-6}}

    \ifdefined\solution
    \solution The consequences of \textit{Dual use technology} are well-explored  in science, and this is something that mathematicians should also start to be aware of.

    The example of Fisher's work here is a unique one, as it is the reverse of dual-use. \textit{Usually}, a new piece of technology is discovered/developed with benevolent intentions, and then later on is reappropriated/misappropriated for more nefarious use. In Fisher's case, the statistical tools he developed were for very harmful purposes, with almost no regard for ethics or morality. However, these same tools turned out to be very useful in doing good, helpful, and socially beneficial statistical analysis.

    This is an example where students might identify and highlight the good, positive uses of these mathematical tools, while at the same time acknowledging and scrutinising their problematic origins. Distinguishing between ``the scientist'' and ``the science'' is often a delicate topic, and one that this question allows students to explore in some depth.

    \fi

        \newpage
    \item Read about the Sally Clark case.\footnote{The corresponding Wikipedia article and press release by the Royal Statistical Society are a good start. \newline Wikipedia about Sally Clark: \href{https://en.wikipedia.org/wiki/Sally\_Clark}{https://en.wikipedia.org/wiki/Sally\_Clark} \newline Statement by the Royal Statistical Society: \href{http://www.inference.org.uk/sallyclark/RSS.html}{http://www.inference.org.uk/sallyclark/RSS.html}} Explain the chain of faulty reasoning that led to her wrongful conviction. Can you explain it to a layperson?

    The Royal Statistical Society has released several guidelines and reports\footnote{These are all available at: \href{https://rss.org.uk/membership/rss-groups-and-committees/sections/statistics-law/}{https://rss.org.uk/membership/rss-groups-and-committees/sections/statistics-law/}} regarding ``Statistics and the Law'', and in particular for the use of statistical in medical misconduct cases. These might be useful a useful starting point for your investigations.

In addition, the RSS released a report specifically dealing with ``uncertainty in medical `murder' cases''\footnote{Available at: \\\href{https://rss.org.uk/news-publication/news-publications/2022/section-group-reports/rss-publishes-report-on-dealing-with-uncertainty-i}{https://rss.org.uk/news-publication/news-publications/2022/section-group-reports/rss-publishes-report-on-dealing-with-uncertainty-i}}. Would it be sufficient to just send this to the legal team(s) and/or judiciary? What else would be right to do? And what is actually \textit{legal} to do? There are very strict rules on what external parties are, and are not, allowed to say or do in regards to criminal proceedings that are in progress.

    \ifdefined\solution
    \solution Sally Clark's conviction, based on several pieces of very bad mathematics, is an example of the power that mathematics has over people. Mathematical results are often seen as `neutral' or `unquestionable' -- and mathematicians are sometimes to blame for advertising their subject as that way. While that description might apply to high school geometry (most people's encounter with `real' mathematics), it certainly does not apply to the mathematics that is used in human contexts: like any other science, mathematics suffers from `garbage in, garbage out'. 
    
    Clark was not the only person to be wrongfully imprisoned by Roy Meadow's bad testimony: two other mothers, Trupti Patel and Angela Cannings, were convicted under similar circumstances. 
    
    \fi

\end{enumerate}

\newpage
\section{Projects: Longer Essay Topics}

\subsection{Introduction to essay tasks}

\subsubsection{Aims}

These  essays serve many purposes, and there are many new skills that you can learn, and knowledge that you can gain. Ultimately, what you get out of this is entirely up to you and what you choose to put in; there are no set expectations. As a guide, here are some things that you \textit{might} learn or gain from doing the essay:

\begin{enumerate}
  \setlength\itemsep{1em}
  
 \item Learning to use \LaTeX.
 \\Yes, this is important! You're mathematicians, or are in the process of becoming mathematicians. Our community communicates in documents written almost exclusively in \LaTeX. This is excellent practice of \LaTeX  formatting. Better to do it now in a non-examinable relaxed setting, than later in a masters or PhD thesis! For example, when you first try and do quotations marks like ``this'', you'll probably end up doing "this".
 
 \item Looking up research papers.
 \\It is actually a non-trivial task to find research papers and preprints online. Various publications are placed behind paywalls, and although you all have access to most of these (in theory) via a university account login, getting to them is never very easy. You'll also learn about arXiv, which is where mathematicians deposit preprints of papers before they go through the (long and arduous) refereeing and publication process. It is not uncommon for an arXiv preprint to go up online in year $N$, yet for the published paper to appear in year $N+3$.
 
 \item Writing prose.
 \\When was the last time you wrote an entire paragraph of text that didn't include any numbers, variables or subscripts? And what proportion of the sentences that you've written in the last 3 years have begun with words other than ``If, then, else''? If we want to communicate mathematically-based ideas to non-mathematicians, then we need to learn how to write prose! This is important: we can't send politicians or upper management lines of C++ code or complicated formulae and expect them to extract the same meaning that we do. We need to \textit{write} for them, and they simply won't read code or formulae.
 
 \item Looking into something with technical and non-technical aspects.
 \\Some of this will involve looking at the technical details of the work of various mathematicians; you're probably already quite competent at that. But you'll also need to do things like make judgement calls, interpret the actions, and analyse and predict things in a qualitative way (``The effects were harmful'', without having a precise definition of ``harmful''.) You'll need to discuss what the mathematicians did, and why it may have been a bad idea! You'll need to make an \textit{argument}, not a \textit{proof}. This is a different skill to what is usually taught in a mathematics degree, and now is as good a time as any to practise it. Remember, mathematics lecturers say ``Here are \textit{true things}'', while lecturers in all other faculties simply say ``Here are \textit{things}, and here is my argument as to why they are  worth considering''.
  
 \item Doing research on Ethics in Mathematics.
 \\This is \textbf{hard}. There are very few well-defined resources on this. So you may find it somewhat unproductive to simply Google ``Ethics in Mathematics''. You'll need to find and read the existing works, and then analyse them in your own mind (think of this as rotating them to get a different/better/more useful point of view), and then discuss them from that particular point of view. There are no right or wrong answers here. The challenge is that you need to find something sensible to study, and then figure out something sensible to say about it.

\end{enumerate}

 \subsubsection{Format and length}
 
 Obviously, this should be formatted in \LaTeX. Which style file you use is up to you (if you don't know what that means, then look it up). As for length: try not to write a phone book here. Perhaps aim for 3000--4000 words, and have a reasonable hard upper bound (e.g. 6000 words) that prevents you from going overboard with your ideas.

\textbf{Don't feel pressured to complete everything that is listed in an essay description; you can chop and choose as you like. So long as you present \textit{something} resembling a complete argument (so if an essay has 3--4 parts, perhaps you cover 1--2 of them). If you're feeling creative, you can add/modify part of an essay, provided it makes sense to do so.}
 
 Try and lay out your essay as you might lay out an undergraduate/masters essay, with a title, abstract, table of contents, introduction, conclusions, and proper (full) referencing. There are places where it may be necessary to insert some formatted mathematics; that's fine, and you can be the judge how much detail to put in. A couple of lines is ok, but a page of algebra should be relegated to an appendix (or removed entirely), as this is not a mathematics essay!

\newpage
 \subsection{Essay: Algorithmic fairness}
 
 With the rapid development and deployment of algorithmic processes in society, many people have raised concerns about potential bias in these algorithms, and the ways in which that bias may affect people's lives. Some mathematicians seem to have ``come to the rescue'', creating systems that identify and measure bias in algorithms, but this has spawned even more (technical and non-technical) debate about how reasonable the bias-checking processes are! The issue is complicated.

We start by looking at a particular example of algorithmic decision making: the prison sentencing algorithm known as COMPAS. You can find a reconstruction of the algorithm in the paper by Rudin et al. (2019) \textit{The age of secrecy and unfairness in recidivism
prediction}: \href{https://arxiv.org/pdf/1811.00731.pdf}{https://arxiv.org/pdf/1811.00731.pdf}.

\begin{center}
  \textit{1. Summarise what the COMPAS algorithm does, how it works, how it's derived, and how it's applied.}
\end{center}

The work of ProPublica was one of the first pieces to give a thorough analysis and testing of the COMPAS algorithm. It can be downloaded from:
\\ \href{https://www.propublica.org/datastore/dataset/compas-recidivism-risk-score-data-and-analysis}{https://www.propublica.org/datastore/dataset/compas-recidivism-risk-score-data-and-analysis}

\begin{center}
  \textit{2. Download the COMPAS data set and explore it. Which variables seem to be ethically relevant and which variables were not collected but could have given insight into ethical issues?}
\end{center}

You could then read the original ProPublica articles:
\\ \href{https://www.propublica.org/article/machine-bias-risk-assessments-in-criminal-sentencing}{https://www.propublica.org/article/machine-bias-risk-assessments-in-criminal-sentencing}
\\\href{https://www.propublica.org/article/how-we-analyzed-the-compas-recidivism-algorithm}{https://www.propublica.org/article/how-we-analyzed-the-compas-recidivism-algorithm}

\begin{center}
  \textit{3. Compare the insights of ProPublica with your own perspective on the data set and your own ideas on algorithmic fairness.}
\end{center}

There are several technical, and non-technical, papers which deal with fairness in algorithms. Here are a few that you can read over:
\begin{itemize}
 \item Cynthia Dwork et al., ``Fairness Through Awareness''. Technical paper giving a mathematical way of measuring if an algorithm is fair.

 \item Cynthia Dwork et al., ``Fairness Under Composition''. Technical paper showing how composition of fair algorithms need not be fair.
 
 \item Sorelle Friedler et al., ``On the (im)possibility of Fairness''. Technical paper showing that various mathematical definitions of ``fair'' are inconsistent.
 
 \item Michael Skirpan, ``The Authority of `Fair' in Machine Learning''. Non-technical paper discussing the  fairness of \textit{trying} to use an algorithm for a particular decision-making task.

 \item Alexandra Chouldechova et al., ``The Frontiers of Fairness in Machine Learning''. A semi-technical paper surveying existing works which discuss algorithmic bias and fairness.

 \item Reuben Binns, ``Fairness in Machine Learning: Lessons from Political Philosophy''.  A non-technical paper discussing how mathematicians are trying to re-invent the wheel when it comes to describing fairness, and that philosophers have already thought about this.
 
\end{itemize}
 
\begin{center}
  \textit{4. Try and make a summary of the algorithmic fairness debate, or at least certain parts of it.}
 \end{center}

 \begin{center}
  \textbf{Remember, you are free to pick and choose whichever aspects of the essay interest you the most.}
 \end{center}

\newpage
 \subsection{Essay: Mathematicians and cryptography}
 
 In recent years the work of mathematicians in state surveillance and cryptography has come under some degree of public scrutiny. This essay aims to explore some of that scrutiny, and look at the arguments that mathematicians, along with computer scientists, are making that might be changing the way mathematicians view their own involvement in cryptographic work.
 
 We take as our starting point the Snowden revelations of 2013. There it was revealed that the NSA, in collaboration with many other surveillance agencies collectively known as the ``Five eyes'' network, had engaged in systematic mass digital surveillance of practically the entire internet-connected world. This is a very very large topic, and not one that can be covered in depth over a matter of weeks. However, the technical details are still of interest. You can find a list of news articles here: \href{https://github.com/iamcryptoki/snowden-archive}{https://github.com/iamcryptoki/snowden-archive} , and the complete document archive here: \href{https://search.edwardsnowden.com}{https://search.edwardsnowden.com} . Thus, the first part of the essay is:
 \begin{center}
 \textit{1. Discuss some of the technical capabilities uncovered in the Snowden revelations, give technical details of how these were carried out, and discuss the sorts of technically-trained people who may have done this work.}
 \end{center}

 From this treasure trove of information spawned many discussions, and many pieces of serious academic work. An early discussion piece was written by Tom Leinster in 2014 titled ``Maths spying: The quandary of working for the spooks'', available here:
 \\\href{https://www.newscientist.com/article/mg22229660-200-maths-spying-the-quandary-of-working-for-the-spooks}{https://www.newscientist.com/article/mg22229660-200-maths-spying-the-quandary-of-working-for-the-spooks}
 \\There, he begins to discuss the merits and drawbacks of the involvement of mathematicians in such work. So the second part of the essay is:
 
 \begin{center}
  \textit{2. Analyse the points that Leinster raises in his piece, look at who else is making these points, and see if you can document a ``community'' of vocal objectors and their specific objections or arguments.}
 \end{center}

 A much more detailed academic paper written on this was the (now very well-known) work by Phillip Rogaway titled ``The Moral Character of Cryptographic Work'' in late 2015, available here:
 \\\href{https://eprint.iacr.org/2015/1162}{https://eprint.iacr.org/2015/1162}
 \\In it, he gives a history and overview of the role of cryptography and cryptographers in society; past, present, and future. This is a long piece, and covers many sub-points and case studies. So the third part of the essay is:
 
 \begin{center}
  \textit{3. Read through the piece by Rogaway and summarise his main arguments, both technical and non-technical. Perhaps choose one or two particular case studies that he mentions, and investigate them further.}
 \end{center}

 If choosing a case study from the piece by Phillip seems daunting, here is one I have pre-selected that involves a very technical academic paper by David Adrian \textit{et al.} from 2015 titled ``Imperfect forward secrecy: how Diffie-Hellman fails in practice'', available here:
 \\\href{https://cacm.acm.org/research/imperfect-forward-secrecy/}{https://cacm.acm.org/research/imperfect-forward-secrecy}
 
  \begin{center}
  \textit{4. Discuss the technical aspects covered in the piece by Adrian \textit{et al.}, as well as their practical consequences.}
 \end{center}

And finally, see if you can find out who the following quote is from: \textit{``The laws of mathematics are very commendable, but the only law that applies in Australia is the law of Australia.''}

 \begin{center}
  \textbf{Remember, you are free to pick and choose whichever aspects of the essay interest you the most.}
 \end{center}

\newpage
 \subsection{Essay: The Ofqual grading algorithm}

 In mid-2020 the UK Office of Qualifications and Examinations Regulation (Ofqual) produced an algorithm to standardise exam grades across England. Though Ofqual had been involved with standardisation algorithms in the past, this was the first time that \textit{all} student grades had been \textit{completely} predicted by the agency, done so in response to the cancellation of all A-level exams as a result of COVID-19. This algorithm was somewhat complicated, relied on many assumptions, took in various data sets, and was designed to achieve certain objectives. Ofqual have produced a handy 319 page report on the algorithm, ``Awarding GCSE, AS \& A levels in summer 2020: interim report'', which can be found here:
 \\ \href{https://www.gov.uk/government/publications/awarding-gcse-as-a-levels-in-summer-2020-interim-report}{https://www.gov.uk/government/publications/awarding-gcse-as-a-levels-in-summer-2020-interim-report}
 
 \begin{center}
 \textit{1. Analyse the Ofqual algorithm and discuss what it did, how it worked, and what it was designed to achieve. Comment on how various assumptions and parameters were made and set, how well it reflected what it was designed to achieve, and how reasonable those aims were.}
 \end{center}

  The public reaction to the grades produced by Ofqual was one of extreme criticism, and even the UK prime minister eventually referred to it as a ``mutant algorithm'' when he announced that all grades would revert back to those predicted for students by their teachers.
 
 \begin{center}
 \textit{2. Give a commentary on why the public reacted so badly to the output of Ofqual's algorithm. How did Ofqual defend it as ``fair'', and what public criticisms were made arguing that it was ``unfair''? How accessible to the public was their $319$ page report? What ultimately led to the UK government abandoning the algorithm?}
 \end{center}
 
 This algorithm was built at the request of the UK government. But ultimately, it was designed and implemented by mathematicians, tasked with ``solving'' a hard and important problem. The mathematical modelling process can sometimes be flawed, with rather bad consequences, as outlined in Chiodo and  M\"uller's \textit{Manifesto for the Responsible Development of Mathematical Works}, found here: \href{https://arxiv.org/abs/2306.09131}{https://arxiv.org/abs/2306.09131}. In particular, this article identifies 10 steps common to most mathematical development process (and in particular, modelling) which might lead to a problematic output.

 \begin{center}
 \textit{3. Explain what parts of the algorithm and its development process led to such an unacceptable output. What where the problems with the assumptions and data sets used by the mathematicians? Identify and comment on the mathematical and non-mathematical aspects of development you see as being problematic. You might use the 10-point list in the above manifesto as a starting point.}
 \end{center}

 Ultimately, the Ofqual algorithm had no impact on the eventual grades received by students. Thus, it achieved nothing, apart from creating a great deal of angst for the students for a week or so. 
 
  \begin{center}
 \textit{4.  What might the mathematicians at Ofqual, and the agency as a whole, done differently? What sorts of oversights or mistakes did they make, what pressures and constraints were they under, and was it even possible to produce a ``reasonable'' algorithm or was this outcome a forgone conclusion?}
 \end{center}

 \begin{center}
  \textbf{Remember, you are free to pick and choose whichever aspects of the essay interest you the most.}
 \end{center}

\newpage
 \subsection{Essay: Modelling COVID-19}
 
 In April 2020 the Virtual Forum for Knowledge Exchange in the Mathematical Sciences (V-KEMS) hosted a 2 day study group which undertook some mathematical modelling on how to re-open industry in UK after the national lockdown (to prevent COVID-19 transmission) ended. At the time it was still unknown how to safely open up shops, factories, etc, so that they could operate in a somewhat effective manner and yet still keep COVID-19 transmission rates low.  A summary document from that study group, titled ``Guiding Principles for Unlocking the Workforce - What Can Mathematics Tell Us?'' can be found here: 
 \\\href{https://gateway.newton.ac.uk/news/2020-05-12/10198}{https://gateway.newton.ac.uk/news/2020-05-12/10198}
 
  \begin{center}
 \textit{1. Summarise the aims and purpose of the study group, and who it comprised of. What problems they did they approach, why did they choose those problems, and what tools did they use to try and produce meaningful output.}
 \end{center}
 
 Understanding and modelling a new disease is a very difficult task, especially when one does not have precise or complete information about the spaces and places where disease transmission is of concern, such as in workplaces.
 
 \begin{center}
 \textit{2. Give an outline of some of the work done in this study group, and the preliminary results obtained. Comment on how the problems were set up and ``mathematised'', what assumptions were made, and what data was used to support the modelling.}
 \end{center}
 
 Often, mathematicians look to ``help'' in situations and scenarios where they believe their work might be of assistance. However, their perspectives and insights might not be sufficient to fully appreciate the problem at hand, yet their output is given great authority as it is presented as ``The mathematics \textit{says} X''. Some mathematicians will claim that ``Mathematical modelling \textit{always} improves a situation or problem.'' But The mathematical modelling process can sometimes be flawed, with rather bad consequences, as outlined in Chiodo and  M\"uller's \textit{Manifesto for the Responsible Development of Mathematical Works}, found here: \href{https://arxiv.org/abs/2306.09131}{https://arxiv.org/abs/2306.09131}. In particular, this article identifies 10 steps common to most mathematical development process (and in particular, modelling) which might lead to a problematic output.
 
 \begin{center}
 \textit{3. Comment on as many potential flaws or shortcomings as you can find in this modelling output, and/or in the process that led to it. How might the use of this work ``go wrong'' and potentially lead to harm? You might use the 10-point list in the above manifesto as a starting point.}
 \end{center}

 The same organisation (V-KEMS) ran another study group in June 2020, this time focusing on reopening UK universities in a COVID-secure manner. They produced a working paper titled ``Unlocking Higher Education Spaces - What Might Mathematics Tell Us?'', which can be found here: \href{https://gateway.newton.ac.uk/event/tgmw82}{https://gateway.newton.ac.uk/event/tgmw82}

 \begin{center}
 \textit{4. Summarise the objectives and output of this second report, and comment on any potential flaws or shortcomings you see in it. How does it differ from the first report, in terms of the issues raised in the 8-point list in the above SIAM article? If there are any noticeable differences, then what might have happened in the few months between the reports to effect such change?}
 \end{center}

 \begin{center}
  \textbf{Remember, you are free to pick and choose whichever aspects of the essay interest you the most.}
 \end{center}

\newpage
\subsection{Essay: Medical Statistics and Probability}

 Retinal detachment (RD) is a serious medical emergency in which the retina (the light-sensitive tissue at the back of the eye) separates from the eye. RD can lead to permanent blindness if not treated quickly. Myopia (short-sightedness) is a risk factor for RD. 
    
    An optician says: `Retinal detachment affects about $1$ in $10,000$ people, but $40\%$ of RD cases happen in people with severe myopia.'
    
    You have severe myopia. How worried should you be about RD? 
    
    Hints: 
    \begin{itemize}
        \item You will need supplementary information, such as  Haimann \textit{et al.} (1982) \\(\href{https://www.ncbi.nlm.nih.gov/pubmed/7065947}{https://www.ncbi.nlm.nih.gov/pubmed/7065947})
        \\and Kempen \textit{et al.} (2004) (\href{https://doi.org/10.1001/archopht.122.4.495    }{https://doi.org/10.1001/archopht.122.4.495    }).
        \item The optician has actually misrepresented the findings of Haimann \textit{et al.} (1982). Can you see why?
    \end{itemize}

    \ifdefined\solution
    \solution \textit{This problem gets students to further investigate what statistics actually mean, and to see that they can be very easily misrepresented.}
    
    The optician's statements by themselves do not tell you anything about your own probability of RD: you will need to know the prevalence of severe myopia in order to calculate your own probability, using Bayes' formula.
    
    Haimann \textit{et al.} (1982) report that $1$ in $10,000$ people is the \emph{annual incidence} of RD. The lifetime risk of RD is rather higher; a crude estimate is to multiply the annual incidence by the number of years one is alive.
    
    The statistics reported in those papers apply to very particular populations: for example, the subjects in Haimann \textit{et al.} (1982) were all from Iowa. 
    
    Myopia depends on both genetics and environment, so more information about incidence broken down by sex, ethnicity, age, \textit{etc.} would be needed to give a more detailed probability for an individual. 
    \fi

\newpage
 \subsection{Essay: A topic of your choosing!}
 
 Yes, that's right: you get to choose your own topic! This is a horribly double-edged sword, and not for the faint-hearted. You can choose any topic or point relating to Ethics in Mathematics that you have come across in the news, in the literature, or that you have even been part of yourself. It really is up to you.

 There is a course given at MIT on ``The Ethics and Governance of Artificial Intelligence'', and they have a very comprehensive reading list which might give you some starting points. It can be found here:
 \\ \href{https://www.media.mit.edu/courses/the-ethics-and-governance-of-artificial-intelligence/}{https://www.media.mit.edu/courses/the-ethics-and-governance-of-artificial-intelligence/}
 
 There is also the now-popular (though somewhat devoid of technical detail) book by Cathy O'Neil titled ``Weapons of math destruction'', which you can read through and look for inspiration.

 And there is the following ``Manifesto for the Responsible Development of Mathematical Works'', where you might take one of the sub-questions and reverse-engineer it to find examples of mathematical work where that particular aspect has gone terribly wrong:
 \\\href{https://arxiv.org/abs/2306.09131}{https://arxiv.org/abs/2306.09131}

 Whatever topic you choose, be sure to keep it as close to ethics in \textit{mathematics} as possible; it can be very easy and tempting to wander in to other neighbouring disciplines such as computer science, engineering and physics.

\newpage
\section{Handout: Further Reading}
\label{sec:further_reading}
\nocite{*}

Here we provide an (initial) reading list of papers, articles, books, etc on Ethics in Mathematics for students to read. We note that Allison Miller also maintains a substantive reading list on Ethics in Mathematics, which is available at:
\\\href{https://docs.google.com/document/d/1P4xAUa9qfEbeCjurGOdJH3dW2JErVRJZ1K9DLAT8uYE/edit}{https://docs.google.com/document/d/1P4xAUa9qfEbeCjurGOdJH3dW2JErVRJZ1K9DLAT8uYE/edit}

\printbibliography[keyword={basics},title={Neutrality and Politics of Mathematics}]
\printbibliography[keyword={teaching},title={Ethics and Mathematics Education}]
\printbibliography[keyword={community},title={The Mathematical Community}]
\printbibliography[keyword={justice},title={Ethnomathematics and Social Justice}]
\printbibliography[keyword={guides},title={Examples of Ethical Engagement}]
\printbibliography[keyword={society}, title={Mathematics, Science and Society}]

\end{document}